\documentclass{amsproc}

\newtheorem{theorem}{Theorem}[section]
\newtheorem{lemma}[theorem]{Lemma}
\newtheorem{corollary}[theorem]{Corollary}
\newtheorem{proposition}[theorem]{Proposition}

\theoremstyle{definition}
\newtheorem{definition}[theorem]{Definition}
\newtheorem{example}[theorem]{Example}

\theoremstyle{remark}
\newtheorem{remark}[theorem]{Remark}

\numberwithin{equation}{section}

\newcommand{\R}{\mathbb R}

\DeclareMathOperator{\ind}{ind} 
\DeclareMathOperator{\sign}{sign}

\usepackage{amsfonts,amsmath,amssymb}
\input xy
\xyoption{all}
\input epsf
\usepackage{amscd}

\newlength{\rig}

\newlength{\hei}

\newcommand{\inv}{{\rm inv}}

\DeclareMathOperator{\im}{Im} \DeclareMathOperator{\Vect}{Vect}
\DeclareMathOperator{\defind}{def-ind} \DeclareMathOperator{\Mod}{mod}

\newcommand{\func}[1]{{\rm #1} \,}
\newcommand{\limfunc}[1]{{\rm{#1}}}
\def\stackunder#1#2{\mathrel{\mathop{#2}\limits_{#1}}}%

\newcommand{\fgr}[3]{
\setlength{\rig}{-0.05\textwidth} \setlength{\hei}{7.8cm}
\begin{figure}
\rule{\rig}{0in} \epsfxsize=9cm \epsffile{#1} \caption{#3}\label{#2}
\end{figure}
}

\begin{document}

\title[Index defects in the theory of spectral problems]
{Index defects in the theory of spectral boundary value problems}

\author{Anton Savin}
\address{Independent University of Moscow,
Moscow 121002, Bolshoi Vlasevsky Pereulok, Dom 11, Russia}
\email{antonsavin@mtu-net.ru}
\thanks{The work was partially supported by the Russian
Foundation for Basic Research under grants Nos. 02-01-00118, 02-01-06515, and
00-01-00161, by Arbeitsgruppe Partielle Differentialgleichungen und Komplexe
Analysis, Institut f\"ur Mathematik, Universit\"at Potsdam and DFG}

\author{Boris Sternin}
\address{Independent University of Moscow,
Moscow 121002, Bolshoi Vlasevsky Pereulok, Dom 11, Russia}
\email{sternin@mail.ru}


\begin{abstract}
We study index defects in spectral boundary value problems for elliptic
operators. Explicit analytic expressions for index defects in various
situations are given. The corresponding topological indices are computed as
homotopy invariants of the principal symbol.
\end{abstract}

\maketitle
\tableofcontents

\section*{Introduction}

The classical Hirzebruch  formula
\begin{equation}
\label{hi1} \sign M=\int_M \widehat{L}(p_1,p_2,\ldots,p_k)
\end{equation}
expresses the signature of a closed oriented $4k$-dimensional manifold in terms
of its Pontryagin characteristic classes. From the  viewpoint of elliptic
operator theory,  formula~\eqref{hi1} expresses the index of a specific
elliptic operator (later called the Hirzebruch operator) via  stable homotopy
invariants of  its principal symbol.  (For the Hirzebruch operator, these
invariants coincide with the Pontryagin classes of the manifold.)

Unfortunately,  formula~\eqref{hi1} has no immediate
analog for manifolds with boundary: there are
examples showing that the signature of such a manifold
 cannot be expressed in terms of
Pontryagin classes.

In 1973, Hirzebruch~\cite{Hirz2} considered a class of manifolds with
boundary (arising from algebraic-geometric
considerations on Hilbert modular varieties) that have naturally defined
relative Pontryagin classes.  Although the right-hand side of~\eqref{hi1} makes sense in this case, the equality
in~\eqref{hi1}  fails. The
difference between the right- and left-hand sides  was called the \emph{signature defect},
and the problem was to compute  it,
i.e.,  find a function $f$ of the boundary of the manifold
such that the difference
$$
\sign M-f(\partial M)
$$
can be expressed via Pontryagin classes of the manifold, or, in the language of
elliptic theory, via the principal symbol of the Hirzebruch operator.

 Hirzebruch conjectured a formula for  $f(\partial M)$ and proved  it in a
number of examples.  A complete solution of the
signature defect problem was given later by Atiyah--Donnelly--Singer
\cite{ADS1} and M\"uller \cite{Mul5}.

The aim of the present survey  is to describe index defects for some natural
classes of general elliptic operators on manifolds with boundary. We consider
only boundary value problems. Note, however, that index defects  also occur in
completely different situations, e.g., for elliptic operators  in
pseudodifferential subspaces  of Sobolev spaces (rather than in  Sobolev spaces
themselves) on compact closed manifolds (see \cite{SaSt1,SaSt2}). Here we  deal
with two classes of operators important in applications, namely, operators
satisfying Gilkey's parity condition~\cite{Gil2}  and operators on
$\mathbb{Z}_n$-manifolds in the sense of Freed and Melrose~\cite{FrMe1}, and
give explicit index defect formulas in both cases.

\subsection{The classical theory: the
Atiyah--Singer and Atiyah--Bott index formulas}
\hspace{1mm}

\textbf{The Atiyah--Singer formula on closed manifolds.} Let $D$ be
an elliptic operator, say, in  Sobolev
spaces on a closed manifold $M$. It is well known that  $D$ is Fredholm. The celebrated Atiyah--Singer theorem \cite{AtSi1} gives a
topological formula for  the index $\ind D$ in terms of
the principal symbol $\sigma(D)$.  By applying the difference construction to the
principal symbol,  one obtains an element
$$
[\sigma(D)]\in K_c(T^*M)
$$
in the $K$-group with compact supports
of the cotangent bundle  $T^*M$. The
Atiyah--Singer formula  reads
\begin{equation}
\label{atsf}
\ind D = \ind_t[\sigma(D)],
\end{equation}
where $\ind_t[\sigma(D)]$ is a functional of the principal
symbol of the operator which can be written out in closed form. In other
words, the Atiyah--Singer formula expresses an analytic invariant of
the operator (the index) in terms of  topological
invariants of the principal symbol.

\textbf{The Atiyah--Bott index formula for boundary value problems.}  If the
boundary $\partial M$ is not empty, then the operator $D$ is no longer Fredholm
(one can show that it always has an infinite-dimensional kernel), and one
should equip it with boundary conditions to obtain a well-posed problem. The
classical boundary conditions are most natural.

A \textit{classical boundary value problem} is a system of equations of the
form\footnote{To simplify the presentation, we consider first-order operators
and  occasionally speak of functions instead of sections of vector bundles.}
\begin{equation}\label{cla}
\begin{cases}
Du=f,
\\
B\left(u\big|_{\partial M}\right)=g,
\end{cases}
\end{equation}
where $u$ and $f$  are functions on $M$ and $g$ is a function on $\partial M$.
The operator $B$ in the boundary condition is  a differential operator; it is
applied to the restriction of the unknown function to the boundary.

The ellipticity condition for problem~\eqref{cla} (see \cite{Hor3})
 can be stated in terms of the principal
symbols $\sigma(D)$ of the operator and
$\sigma(B)$ of the boundary condition. Atiyah and Bott \cite{AtBo2}
showed that the index theory  of classical boundary value
problems  is similar to
that of elliptic operators on closed manifolds. Namely,  problem~\eqref{cla} defines a difference element
$$
[\sigma(D,B)]\in K_c(T^*(M\setminus \partial M)),
$$
where $T^*(M\setminus \partial M)$ is the cotangent bundle over the
interior of  $M$, and the
index of the corresponding Fredholm operator  is given
by the formula
\begin{equation}
\label{mumq} \ind(D,B)=\ind_t[\sigma(D,B)],
\end{equation}
similar to~\eqref{atsf}.

However, the theory of classical boundary value problems has
an essential  drawback. \emph{For  some operators, there  are no
well-posed\footnote{That is, defining a Fredholm problem
 for the original
differential expression $D$  in suitable
spaces.} classical boundary conditions at all}!

Atiyah and Bott showed that the obstruction to the existence of
well-posed boundary conditions  is of topological
 nature and computed
it.  The obstruction
proves to be nonzero for most  geometric operators:
the Dirac operator, the Hirzebruch operator, and
the Cauchy--Riemann operator. In other words,  there are no Fredholm classical boundary value
problems for these operators.

 It is still possible to
sidestep the obstruction  and, in particular,
 equip the
above-mentioned operators with well-posed boundary
conditions. To this end,  one has to consider a
more general class of boundary value problems, namely,
 so-called \emph{problems in subspaces}, which
are described in the next subsection.

\subsection{Boundary value problems in subspaces. Spectral problems.}
A \emph{boundary value problem in subspaces} is a boundary value problem of the
form
\begin{equation}
\begin{cases}
Du=f,\\
B\left(u\big|_{\partial M}\right)=g, & g\in \im P,
\end{cases}
\end{equation}
where the right-hand side $g$ of the boundary condition  lies in  the range
$$
\im P\subset C^\infty(\partial M, G)
$$
of a pseudodifferential projection  operator
$$
P:C^\infty(\partial M, G)\to C^\infty(\partial M, G)
$$
in the function space on the boundary. This class of boundary value
problems was  introduced in~\cite{ScSS18} and further
 studied in~\cite{NScSS3,SaScS8,SaScS4}. In
particular, it was shown that the ellipticity
condition  can be stated in terms of the principal symbols of
$D$, $B$, and $P$, just as in the classical case. However,
from the topological point of view  these
 problems  are opposite
to classical boundary value problems.

 The two most important differences
are as follows.\vspace{1mm}

I. \emph{There exists a Fredholm
boundary value problem in subspaces for an arbitrary elliptic operator}.
 An example is given by the \emph{spectral
Atiyah--Patodi--Singer} boundary value problem \cite{APS1}\footnote{Atiyah,
Patodi, and Singer
 used only homogeneous boundary
conditions. However,
 problem~\eqref{apso} is
equivalent to the corresponding homogeneous problem as far
as the solvability and the index problem are concerned.}
\begin{equation}
\label{apso}
\begin{cases}
Du=f, \\
\Pi_+(A)u|_{\partial M}=g, & g\in \im \Pi_+(A).
\end{cases}
\end{equation}
Here $D$ is a
first-order operator assumed to have  the form
\begin{equation}
D\big|_{U_{\partial M}}\simeq \frac\partial{\partial t}+A
\label{omgo}
\end{equation}
in a collar neighborhood $U_{\partial M}$ of the boundary, where $A$ is  an
elliptic self-adjoint operator called the \emph{tangential operator}  of $D$,
and $\Pi_+(A)$ is the spectral projection of $A$ on $\overline{\R}_+$, i.e.,
the orthogonal projection on the subspace spanned by eigenvectors of $A$ with
nonnegative eigenvalues.
 Topologically,\footnote{But not
analytically!}  problem~\eqref{apso}  can in
essence be viewed as the general case of a problem in
subspaces, since an arbitrary boundary value problem can be reduced to
 a spectral problem by a stable homotopy (see
\cite{SaSt1,SaScS4}). Therefore, for simplicity we consider  only spectral problems~\eqref{apso} for operators $D$
satisfying~\eqref{omgo}.  By $\ind(D,\Pi_+(A))$ we
denote the index of
problem~\eqref{apso}.\vspace{1mm}

II. \emph{The index of a boundary value problem in subspaces is not
determined by the principal symbol of the operator} $D$.  To illustrate this, consider  a deformation
of  lower-order terms of $D$  such
that some eigenvalue of the tangential operator  changes its sign.  At this point, the spectral projection
$\Pi_+(A)$  experiences a jump,  so that the index of the problem may
change. On the  other hand, the index
remains constant as long as the  deformation produces
continuously varying spectral projections. Let us give a simple example.

Consider the zero-order deformation
$$
D_\tau=D+\tau\chi(t)
$$
of the Cauchy--Riemann operator
$$
D=\frac\partial{\partial t}+i\frac\partial{\partial \varphi}
$$
on the cylinder $\mathbb{S}^1\times[0,1]$, where $\chi(t)$ is a smooth function  such that
$\chi(0)=1$ and $\chi(1)=0$. The tangential operator  of the family $D_\tau$  depends on $\tau$
only on one of the  bases of the cylinder, namely, on
$\mathbb{S}^1\times\{0\}$, where it has the form
$$
A_\tau=i\frac\partial{\partial \varphi}+\tau.
$$
The eigenvalues of $A_\tau$ are given by the
formula $\tau+2\pi n, n\in\mathbb{Z}$.  As $\tau$  passes through zero, one of the eigenvalues changes
its sign, so that  the spectral projection
 undergoes a jump. The index  is
also discontinuous:
$$
\ind(D_\tau,\Pi_+(A_\tau))=
\begin{cases}
-2, & - 2\pi<\tau<0, \\
-1, & \tau =0, \\
0, & 0<\tau< 2\pi.
\end{cases}
$$

This example makes it clear that  one cannot
 obtain an index formula  similar
to~\eqref{mumq} for  spectral boundary value problems; in other words, a topological computation of the index in this
case is  impossible in principle.

\subsection{The index and the index defect for spectral problems}
The  aim of this survey is to show that in
many cases of interest one obtains a homotopy
invariant of the principal symbol of the operator by adding some analytic
invariant to the index of a spectral boundary value problem.

 The correction term
is naturally called  an \emph{index defect} of the problem,
 since it is  this term that restores
the homotopy invariance of the index.  It is natural to
 require that the correction term  be
determined solely by the structure of the operator in a
neighborhood of the boundary,  for  on closed manifolds  the
analytic index itself is homotopy invariant and  zero
can be taken  for the correction term.

Therefore, we introduce the following statement of the index defect
problem.\vspace{1mm}

\textbf{The index defect problem for spectral boundary value
problems}. Construct a functional $\defind (D)$ of elliptic operators $D$ on a manifold with
boundary such that
\begin{enumerate}
\item the sum $\ind (D,\Pi_+(A))+\defind D$  is a homotopy
invariant  of $D$;\footnote{It
readily  follows that the sum is a homotopy invariant of the
principal symbol of $D$.}

\item  $\defind (D)$ is determined solely by the
tangential operator $A$.
\end{enumerate}
A functional with these properties will be called an \textit{index
defect}.  Conditions~(1) and~(2) imply that the
\emph{index defect is determined by the spectral projection and is its
homotopy invariant}.

 Needless to say, to obtain an actual defect formula
(which is our main problem), we should compute the
homotopy invariant  in (1)
topologically, i.e., express it in the form
\begin{equation}
\ind (D,\Pi_+(A))+\defind D=\ind_t(\sigma(D)), \label{defo}
\end{equation}
where $\ind_t(\sigma(D))$  is a functional on the
set of homotopy classes of elliptic principal symbols.

In the remaining part of the introduction, we explain the main methods  that
can be used to define index defects and describe approaches to  the proof of
the corresponding index defect formulas~\eqref{defo}. However, prior to
proceeding to these topics,  we consider the following phenomenon of utmost
importance. \vspace{1mm}

\textbf{The obstruction to index defect formulas.} The desired index
defect formula~\eqref{defo} can be  viewed as a
decomposition of the index of a spectral boundary value
problem into a finite-dimensional contribution of the
principal symbol  and an (infinite-dimensional)
contribution of the tangential operator:
\begin{equation}
\label{muo} \ind (D,\Pi_+(A))=f_1(\sigma(D))+f_2(A),
\end{equation}
where the functional $f_1$ is a homotopy invariant of the principal symbol.

It turns out that  \textit{there is no index decomposition
of  the form~\eqref{muo} on the set of all elliptic
operators}  (see \cite{SaScS1}). Therefore,
index defect formulas and decompositions of the form \eqref{muo} can
 be sought only in some \emph{subsets}  of
the space of elliptic operators on a manifold with boundary.

The obstruction to the existence of decompositions \eqref{muo} was computed in
\cite{SaScS1}.  It is  the
one-dimensional cohomology class of the space of elliptic operators whose
value on a cycle is equal to the Atiyah--Patodi--Singer spectral flow of
the corresponding family of tangential operators.
There exists a decomposition~\eqref{muo}  on a subspace
$\Sigma$  of the space of elliptic operators if and only if the
restriction of  this cohomology class to $\Sigma$ is trivial.

The cited result, unfortunately, proves only the \emph{existence} of a
decomposition~\eqref{muo} and does not give a
satisfactory formula for  an index defect.  To
study  index defects, one has to use other methods.

\subsection{Approaches to the definition of index defects}

Let us briefly describe two methods useful  in defining
 index defects.\vspace{1mm}

{\bf A. The geometric index formula of Atiyah--Patodi--Singer.} In
1975, Atiyah--Patodi--Singer \cite{APS3} obtained  the
formula
\begin{equation}
\label{apsf} \ind (D,\Pi_+(A))=\int_M a(D)-\eta(A)
\end{equation}
for the index of spectral boundary value problems, where the
density $a(D)$ is determined by the coefficients of
$D$, just as in the case of closed manifolds. The new contribution to the index  is given by
the spectral $\eta$-invariant $\eta(A)$ of the tangential operator $A$.

This formula is often called a \emph{geometric index formula}, since for
geometric operators (the Hirzebruch, Dirac, Todd and Euler
operators) the integrand  on the right-hand side  is determined by the metric and coincides with the local
Atiyah--Singer density in  the case of closed
manifolds  (i.e., with the $L$-form for the
Hirzebruch operator, the $A$-form for the Dirac operator, etc.).

Unfortunately, the Atiyah--Patodi--Singer formula does not define an
index defect, since  neither of the terms
on the right-hand side is a homotopy invariant of
 $D$.  However, formula~\eqref{apsf}
can be used to define index defects as follows.

 We have already pointed out that  the index of the spectral boundary value problem
experiences jumps under homotopies of $D$.  However,
the sum
\begin{equation}
\label{suma} \ind (D,\Pi_+(A))+\eta(A)
\end{equation}
varies smoothly by the Atiyah--Patodi--Singer  theorem. Moreover, the
sum~\eqref{suma} is homotopy invariant if and only if for an arbitrary homotopy
$D_\tau$ (with parameter~$\tau$) in our class of operators the derivative
\begin{equation}
\label{defi2}
\frac d{d\tau}\int_Ma(D_\tau)
\end{equation}
is zero. Since the density
$a(D_\tau)$ is given by  a closed-form
expression involving the coefficients of $D_\tau$, one can
use  a detailed analysis of  $a(D_\tau)$ to construct classes of operators for
which the derivative is zero and hence the
$\eta$-invariant of the tangential operator is the desired index defect.
\vspace{1mm}

{\bf B. Operator algebras.}
Another method for defining  index defects relies on
$K$-theory of operator algebras.

As was mentioned already,  any index defect is determined by the spectral
projection $\Pi_+(A)$ and is its homotopy invariant.
 For simplicity, consider the class of matrix
projections
\begin{equation}
\Pi_+(A):C^\infty(\partial M, \mathbb{C}^N)\longrightarrow
C^\infty(\partial M, \mathbb{C}^N) \label{pro1}
\end{equation}
and assume that the matrix entries  lie in some algebra $\mathcal{A}$ of operators on the
boundary.

Then  index defects admit an alternative
description as  homotopy invariants of projections.

 The homotopy
classes of projections~\eqref{pro1} with entries in $\mathcal{A}$ (for
arbitrary $N$) generate  the $K$-\textit{group
$K_0(\mathcal{A})$  of the algebra}
$\mathcal{A}$ (e.g., see \cite{Bla1}).  An index defect functional defines a homomorphism
$$
d: K_0(\mathcal{A})\longrightarrow \mathbb{R}
$$
of the $K$-group  into real numbers.

Moreover, a simple computation shows that the homotopy invariance
 of the  sum
$\ind(D,\Pi_+(A))+d(\Pi_+(A))$ is equivalent
to either of the following two conditions for  the
functional $d$  (provided that $\mathcal{A}$ contains the ideal of
finite rank operators).

\begin{enumerate}
\item The functional $d$ is a dimension type invariant of
projections. More precisely, for  two arbitrary
projections
$$
P_{1,2},\quad \im P_1\subset \im P_2, \quad \dim\im (P_2-P_1)<\infty,
$$
differing by a finite rank projection $P_2-P_1$,
one has
$$
d(P_2)-d(P_1)=\dim\im (P_2-P_1).
$$

\item The functional $d$ defines a commutative diagram
\begin{equation}
\xymatrix{ \quad K_0(\mathcal{K})=\mathbb{Z}
\ar[d] \ar[rd]& \\
K_0(\mathcal{A})\ar[r]^{d} &\mathbb{R}, } \label{treug}
\end{equation}
where the  diagonal arrow is the natural embedding
 $\mathbb{Z}\subset\mathbb{R}$.
\end{enumerate}

\subsection{Examples}
\label{exi} Now we consider  some explicit
 implementations of the methods.  Each example starts with
a description of  the class of operators
for which the index defect is to be considered.

\begin{example}
\label{ex1} {\bf Problems in even subspaces} {\cite{SaSt1}. Consider operators
$D$ such that  the symbol of the tangential operator  is even  in the momentum
variables~$\xi$:
$$
\sigma(A)(x,\xi)=\sigma(A)(x,-\xi).
$$
This condition  singles out the subalgebra of operators with even principal
symbol in the algebra of pseudodifferential operators on the boundary. In this
subalgebra, consider the subalgebra $\mathcal{A}$ of zero-order operators.  If
$\partial M$ is \emph{even-dimensional}, then the vertical  arrow in
diagram~\eqref{treug} is a monomorphism, which implies the existence of the
desired functional $d$. One can also prove  that $d$ is unique under some
natural conditions.

Thus, there is a well-defined homotopy invariant $d(P)\in\mathbb{R}$ on the set
of  pseudodifferential projections $P$ with even principal symbol.  In contrast
with the index of elliptic operators,  this functional is not integer; it can
take arbitrary rational values whose denominators are powers of $2$ (dyadic
values) see \cite{SaSt7}.

Thus, the value of the dimension functional $d$ on the spectral
projection can be taken as the index defect for operators $D$
whose tangential operator $A$  has an even
principal symbol. In other words, the sum
$$
\ind (D,\Pi_+(A))+d(\Pi_+(A))
$$
is a homotopy invariant of $D$. This poses a problem of computing this
invariant in terms of the principal symbol of $D$. In  the next subsection, we explain
the main ideas underlying the computation of this invariant.}
\end{example}

\begin{example}\label{ex2}
{ {\bf Spectral problems on $\mathbb{Z}_n$-manifolds} \cite{SaSt10}. Consider a
manifold $M$ whose boundary is represented as the total space of a covering
$$
\pi:\partial M\longrightarrow X
$$
over  a smooth base $X$. Geometrically, such a manifold can be  viewed as a
smooth model of the singular space $\overline{M}^\pi$ (known as a
$\mathbb{Z}_n$\emph{-manifold}, where $n$ is the number of sheets of the
covering) obtained by identifying the points in  each fiber of $\pi$. A
neighborhood of a singular point looks like an open book with $n$ sheets (see
Fig.~\ref{fi13}), and $X\subset \overline{M}^\pi$ is the edge, where the sheets
meet. \fgr{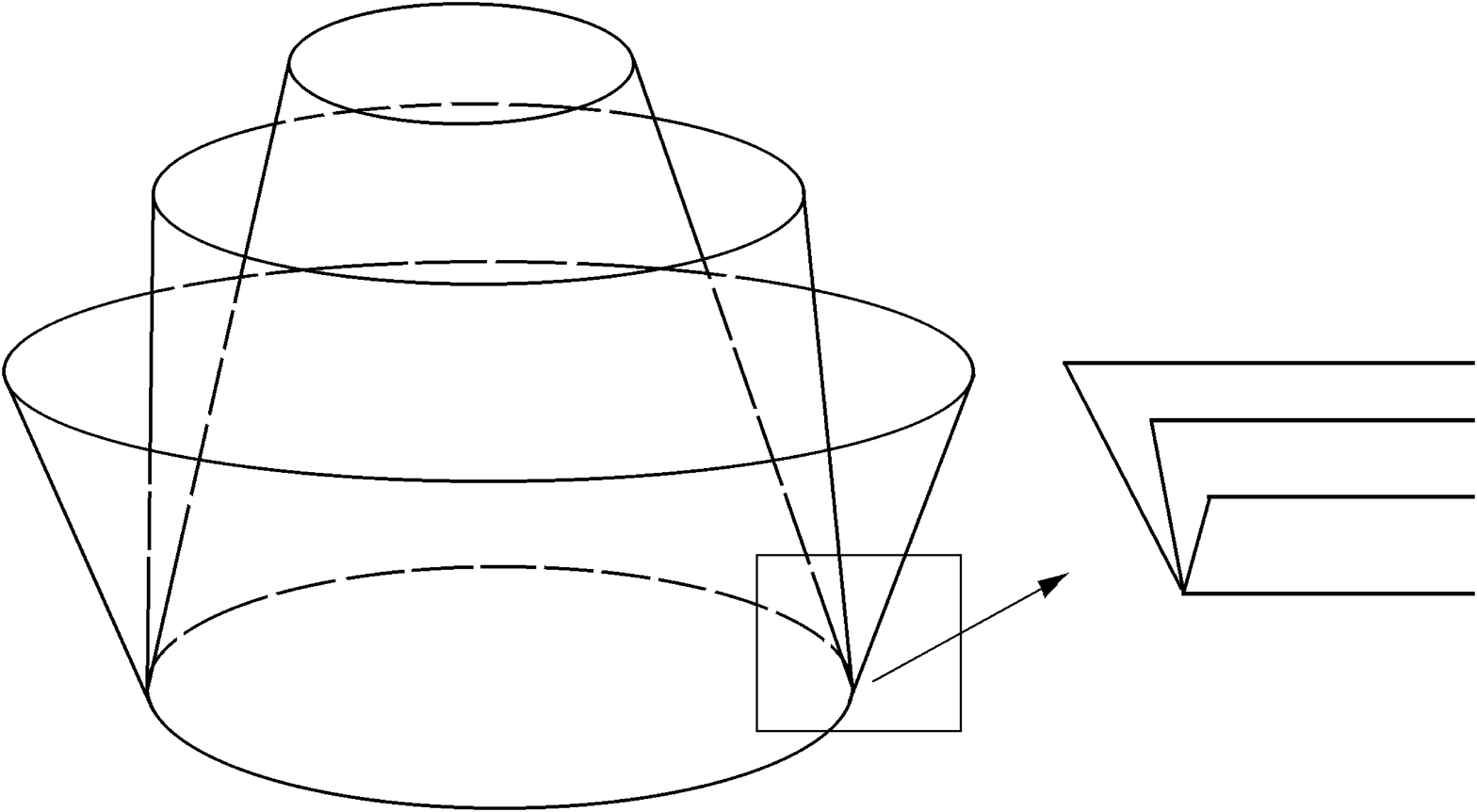}{fi13}{A singular space}

On $M$ we consider elliptic operators whose tangential operator $A$ is  the
lift by $\pi$ of some operator $A_0$ on the base of the covering.  (The lift is
well defined, since $\pi$ is a local diffeomorphism.)

For a trivial covering, this condition guarantees that the index of
the corresponding spectral boundary value problem viewed as
the $\Mod n$-residue
$$
\Mod n\text{-}\ind D\in \mathbb{Z}_n
$$
is a homotopy invariant of $D$. Freed and Melrose \cite{FrMe1} obtained
the $\Mod n$-\emph{index formula}
$$
\Mod n\text{-}\ind D=\ind_t\left[ \sigma \left( D\right) \right],
$$
where
\[
\left[ \sigma \left( D\right) \right] \in K_c\bigl( \overline{T^{*}M}^\pi
\bigr)
\]
is the difference element defined by the principal symbol of $D$ in the
$K$-group of the singular space $\overline{T^{*}M}^\pi $ obtained from the
cotangent bundle $T^{*}M$ by identifying the points in  each fiber of the
covering  $\partial T^*M\longrightarrow T^*X\times\mathbb{R}$,  by analogy with
the definition of $\overline{M}^\pi$.

However,  the $\Mod
n$-index \emph{is not} an invariant of the
principal symbol if the covering is nontrivial. The index defect in
this case turns out to be given by the  difference
$$
\defind (D)\overset{\operatorname{def}}=\eta \left( A\right)-n\eta(A_0) \in
\mathbb{R}/n\mathbb{Z}.
$$
An expression of this type is known as the \emph{relative
Atiyah--Patodi--Singer $\eta$-in\-vari\-ant} \cite{APS2} of an operator $A_0$
with coefficients in a flat bundle. The index defect problem in this case
 is to compute the homotopy invariant
$$
\inv D=\Mod n\text{-}\ind D+\defind D \in \mathbb{R}/n\mathbb{Z}
$$
as a residue modulo $n$. }
\end{example}

\subsection{Approaches to index defect formulas}
\label{cch6}

To state and prove an index defect formula, one
can apply all  methods that are useful in the proof of  ordinary index formulas. We content ourselves with
describing only two approaches.\vspace{1mm}

{\bf Method 1. Homotopy classification.}

Roughly speaking, the method consists of two steps.
\begin{enumerate}
\item First, one  carries out the homotopy
classification of elliptic operators $D$ to be considered,
or,  more technically,  computes the
group of stable homotopy classes of  these operators.

\item Second,  one  finds a generating set
 of  this group,  so that the proof of  an index
defect formula is reduced to its verification
 for the generators.
\end{enumerate}

This scheme goes back to the first proof of the Atiyah--Singer
formula~\cite{AtSi0}, where  elliptic operators are
classified (modulo stable homotopies) by  elements of the $K$-group
$K_c(T^*M)$ of the cotangent bundle at step~(1), the
Hirzebruch operator with coefficients in various vector bundles gives a
rational generating set of the $K$-group (on an orientable even-dimensional
manifold)  at step~(2), and finally cobordism theory is used to compare
the analytic and topological index of these geometric operators.\vspace{1mm}

{\bf Method 2. Poincar\'e duality.}

Another method  for proving index formulas is based on
 Poincar\'e duality in $K$-theory. Let us illustrate this method
 using the classical Atiyah--Singer theorem as an
example.

On a closed manifold, there is  an index homomorphism
\begin{equation}
\label{str} \ind_a:K_c(T^*M)\longrightarrow \mathbb{Z},
\end{equation}
which takes each element of $K$-theory  to the (analytic) index of the corresponding elliptic operator. On
the other hand,  Poincar\'e duality in $K$-theory gives
the pairing
\begin{equation}\label{alf3}
\begin{aligned}
\langle\,,\,\rangle:K_c(T^*M)\times K(M) & \longrightarrow \mathbb{Z},
\\
(x,y) &\longmapsto p_!(xy),
\end{aligned}
\end{equation}
which is nonsingular on the free parts of the groups. Here
$$
p_!:K_c(T^*M)\longrightarrow \mathbb{Z}
$$
is the direct image mapping induced by the projection
$p:M\to pt$ to  a one-point space. It follows  that the homomorphism~\eqref{str} can be
represented as  the pairing with
some element $y\in K(M)$; i.e.,
$$
\ind D=\langle[\sigma(D)],y\rangle
$$
for all elliptic operators $D$, where $y$ is uniquely determined by
$M$ modulo torsion. Therefore, to obtain an index formula, it
suffices to compute the element $y$. In these terms, the Atiyah--Singer
formula states that  one can take
$y=1\in K(M)$, the element corresponding to the trivial line
bundle.

\medskip

Let us show how one can apply these methods  to
find and prove index defect formulas.

\medskip

\noindent \textbf{Continuation of Example \ref{ex1} (an index
defect formula in even subspaces).}

 We consider operators $D$ with even principal
symbol of the tangential operator on an even-dimensional manifold
$M$.
 The homotopy classification of the
corresponding spectral boundary value problems  turns out to be isomorphic
(modulo 2-torsion) to  that of classical boundary value problems,
i.e., to the group $K_c(T^*(M\setminus
\partial M))$. Therefore, to obtain a topological formula for the
homotopy invariant
$$
\inv D=\ind (D,\Pi_+(A))+ d(\Pi_+(A) ),
$$
it suffices to generalize the Atiyah--Bott topological index~\eqref{mumq}
to  boundary value problems in even subspaces. Such a
generalization was obtained in~\cite{SaSt1}.
We point out that the topological index in this formula  proves to be  a half-integer, and a topological
consequence of this formula is  the half-integrality
 of the index defect. The
index defect formula has a number of applications. For
example, it enabled the authors~\cite{SaSt6} to solve  Gilkey's problem on the nontriviality of $\eta$-invariants
of even-order operators on odd-dimensional manifolds. \vspace{2mm}

\noindent{\textbf{Continuation of Example \ref{ex2} (an index
defect formula on $\mathbb{Z}_n$-manifolds).}

For a manifold  whose boundary is
an $n$-sheeted covering, the  sum
$$
{\inv}D\overset{\operatorname{def}}= \Mod n\text{-}\ind D+\defind
D
$$
can be  viewed as a homomorphism
\begin{equation}
\inv:K_c(\overline{T^*M}^\pi)\longrightarrow \mathbb{R}/n\mathbb{Z},
\label{stst}
\end{equation}

The  two main differences
between~\eqref{stst} and~\eqref{str}
are as follows:
\begin{enumerate}
\item now we use the group $\mathbb{R}/n\mathbb{Z}$ rather than
$\R$;
\item the space $\overline{T^*M}^\pi$ has singularities.
\end{enumerate}

To take account  of~(1), one should
 replace the classical Poincar\'e duality~\eqref{alf3} by
``Poincar\'e duality with coefficients," i.e.,
\emph{Pontryagin duality}
\begin{equation}
\label{alf4} \langle,\rangle:K_c(T^*M)\times
K(M,\mathbb{R}/\mathbb{Z})\longrightarrow \mathbb{R}/\mathbb{Z}.
\end{equation}
(See~\cite{SaSt6}; here $K(M,\mathbb{R}/\mathbb{Z})$ is the
$K$-group with coefficients in  the compact group
$\mathbb{R}/\mathbb{Z}$.)

 To tackle~(2), one needs an
extension of  Poincar\'e duality to singular manifolds like
${\overline{M}}^\pi$.  We point out
that duality in the classical sense may not be valid on a manifold
with singularities. However, the desired duality can be obtained in the framework of
Connes' noncommutative geometry~\cite{Con1}. A detailed
exposition  is given in  Appendix C,
 and now we describe only the main ideas involved.

To a singular $\mathbb{Z}_n$-manifold ${\overline{M}}^\pi$, one
 assigns a noncommutative $C^*$-algebra
 $\mathcal{A}_{M,\pi}$, which should be regarded as the
``function algebra" on
${\overline{M}}^\pi$.  Now Poincar\'e duality  gives a pairing
\begin{equation}
\label{bet2} \langle\,,\,\rangle:K_c(\overline{T^*M}^\pi)\times
K_0(\mathcal{A}_{M,\pi})\longrightarrow \mathbb{Z},
\end{equation}
where the second group is the $K$-group of  the $C^*$-algebra
$\mathcal{A}_{M,\pi}$. As a special case of this construction for
$\partial M=\varnothing$, one obtains the pairing \eqref{alf3},
since  in this case the algebra
$\mathcal{A}_{M,\pi}$ coincides with the algebra of continuous functions
on $M$ and the $K$-group of the algebra of continuous
functions on a space coincides with the $K$-group of the space.

Applying the  constructions~\eqref{alf4} and~\eqref{bet2}
to the mapping~\eqref{stst}, one  finds that
$$
\inv D=\langle[\sigma(D)],y\rangle\in \mathbb{Q}/n\mathbb{Z},
$$
just as in the case of a closed manifold, for some element
$$
y\in K_0(\mathcal{A}_{M,\pi},\mathbb{Q}/n\mathbb{Z})
$$
depending only on the manifold. The index defect theorem
for $\mathbb{Z}_n$-manifolds  in the
authors's paper~\cite{SaSt10} gives an explicit formula for this
element.}

\medskip

{\bf Acknowledgements.} The authors wish to express their keen gratitude to
Professor B.-W.~Schulze of Potsdam University, where the paper was written, for
his kind hospitality. We also thank V.~Nazaikinskii, who read a preliminary
version of the paper and made a number of important remarks.

The contents of the paper are as follows. In the first section, we
define  spectral boundary value problems and prove a theorem on
index decompositions.  Section~2 deals with
 index defects in even subspaces. In
Section~\ref{seczn}, we study index defects on
$\mathbb{Z}_n$-manifolds. The paper concludes with three
appendices, the first  dealing with
the Atiyah--Patodi--Singer $\eta$-invariant and the
remaining two  with Poincar\'e duality in $K$-theory on smooth
manifolds and $\mathbb{Z}_n$-manifolds.

\section{Spectral boundary value problems and their index}

\subsection{Atiyah--Patodi--Singer spectral boundary value problems}
We start  by introducing some notation. Let $M$ be a
smooth compact manifold with boundary and
$$
D:C^\infty(M,E)\longrightarrow C^\infty(M,F)
$$
an elliptic first-order differential  operator on
$M$ acting in the spaces of sections of vector bundles $E,F\in
\Vect(M)$.

 We choose some collar neighborhood
$$
U_{\partial M}\simeq \partial M\times [0,1)
$$
of the boundary $\partial M$. The normal coordinate on the
half-open interval $[0,1)$ is denoted by $t$. Then $D$ can
be represented in  this neighborhood
in the form
\begin{equation}
\label{deco2} D\big|_{U_{\partial M}}\simeq \frac\partial{\partial
t}+A(t)
\end{equation}
(here $\simeq$  stands for equality up to a vector bundle isomorphism) for a
smooth family $A(t)$ of elliptic first-order differential  operators on
$\partial M$. This representation can be obtained as follows.  In  the collar
neighborhood of the boundary, we take some isomorphisms
$\pi^*(\left.E\right|_{\partial M})\simeq \left. E\right|_{\partial
M\times[0,1)}$ and $\pi^*(\left.F\right|_{\partial M})\simeq \left.
F\right|_{\partial M\times[0,1)}$, where $\pi:\partial M\times [0,1)\to
\partial M$ is the natural projection. Then $D$ is isomorphic in
$U_{\partial M}$ to an operator in the
spaces
$$
C^\infty\Bigl(\partial M\times [0,1),\pi^*(\left.E\right|_{\partial M})
\Bigr)\longrightarrow C^\infty\Bigl(\partial M\times
[0,1),\pi^*(\left.F\right|_{\partial M})\Bigr),
$$
where the operator $\partial/\partial t$ is
well defined, and we obtain a decomposition
$$
D|_{U_{\partial M}}\simeq \Gamma(t)\frac\partial{\partial t}+A'(t),
$$
where $\Gamma(t)$ is a vector bundle homomorphism. By ellipticity,
$\Gamma(t)$ is an isomorphism, and we arrive
at~\eqref{deco2}.

 For simplicity, we also assume that
for small $t$ the family $A(t)$  is
independent of  $t$ and coincides with a
 given self-adjoint operator
$$
A:C^\infty(\partial M, E|_{\partial M})\longrightarrow
C^\infty(\partial M, E|_{\partial M})
$$
on the boundary. The operator $A$ is called the \emph{tangential
operator}  of $D$.

\begin{definition}{
The \emph{Atiyah--Patodi--Singer spectral boundary value problem} for an
elliptic operator $D$ is  the boundary value problem
\begin{equation}
\label{aps1}
\begin{cases}
Du=f, \\
\Pi_+(A)u|_{\partial M}=g, & g\in \im \Pi_+(A),
\end{cases}
\end{equation}
where $\Pi_+(A)$ is the nonnegative spectral projection of
$A$:
$$
\Pi_+(A)e_\lambda=
\begin{cases}
e_\lambda, & \lambda\ge 0, \\
0, & \lambda<0,
\end{cases}
$$
for any eigenvector $e_\lambda$ of $A$ with eigenvalue $\lambda$.
}
\end{definition}
For an arbitrary elliptic operator $D$, the spectral boundary value
problem defines a bounded Fredholm operator in the  spaces
\begin{equation}
\left({D},{\Pi_+(A)}\right):H^s(M,E)\longrightarrow H^{s-1}(M,F)\oplus
\overline{\im \Pi_+(A)}, \quad s>1/2,
\end{equation}
where $\overline{\im \Pi_+(A)}$ is the closure of the range of $\Pi_+(A)$ in  $H^{s-1/2}(\partial M, E|_{\partial M})$
(see \cite{APS1}).

As usual, the index $\ind(D,\Pi_+(A))$ of the spectral boundary
value problem
is independent of the Sobolev smoothness
exponent $s$ and can be computed in  spaces of smooth
functions.

Note that,  in contrast with the index of elliptic
operators on closed manifolds, the index of the Atiyah--Patodi--Singer
problem is not invariant  under homotopies of
 $D$. Indeed, consider  a smooth homotopy
$$
D_\tau:C^\infty(M,E)\longrightarrow C^\infty(M,F), \quad \tau\in [0,1],
$$
of elliptic operators, i.e. an elliptic operator family with coefficients
smoothly depending on $\tau$. Suppose that  an eigenvalue of the tangential
operator $A_{\tau}$ changes its sign at some point $\tau=\tau'$. Then the
corresponding spectral projection  experiences a jump, and consequently, the
spaces  in which the operator  acts change discontinuously. This intuitive
argument is stated in Theorem~\ref{sft} below in terms of the  spectral flow.

\subsection{The spectral flow}
\label{sflow1}

Consider a smooth family $\left\{ A_\tau\right\} _{\tau\in \left[
0,1\right] }$ of elliptic self-adjoint operators on a closed manifold and
assume that the operators at the endpoints $t=0$ and $t=1$ are
invertible.
\begin{definition}{\em
The { spectral flow $\limfunc{sf}\left\{ A_\tau\right\} _{\tau\in \left[
0,1\right] }$ \emph{of the family }$\left\{ A_\tau\right\} _{\tau\in \left[
0,1\right] }$ {\em is  the}
\emph{net number of eigenvalues of }$A_\tau$\emph{\
 that change their sign from minus to plus as the
parameter $\tau $ increases from $0$ to $1$ (see Fig.~\ref{fi14}). } }}
\end{definition}
Unfortunately, this definition does not makes sense for general families. The
spectral flow  of an arbitrary family $\{A_\tau\}$ can be defined by  different
methods (see \cite{Phil1}, \cite{Mel2}, \cite{Sal1}, \cite{BBW1}, \cite{DaZh3},
\cite{NScS5}, and other papers). For  example, one can slightly  deform the
straight line $\lambda=0$ in the $(\lambda,\tau)$-plane in a way such that the
spectral curves of the operators $A_\tau$ intersect the perturbed line
transversally. Then the spectral flow can be defined as the intersection number
of the perturbed line with the graph of the spectrum of the family. The desired
perturbation can be constructed explicitly as follows (e.g., see \cite{Mel2}).
\fgr{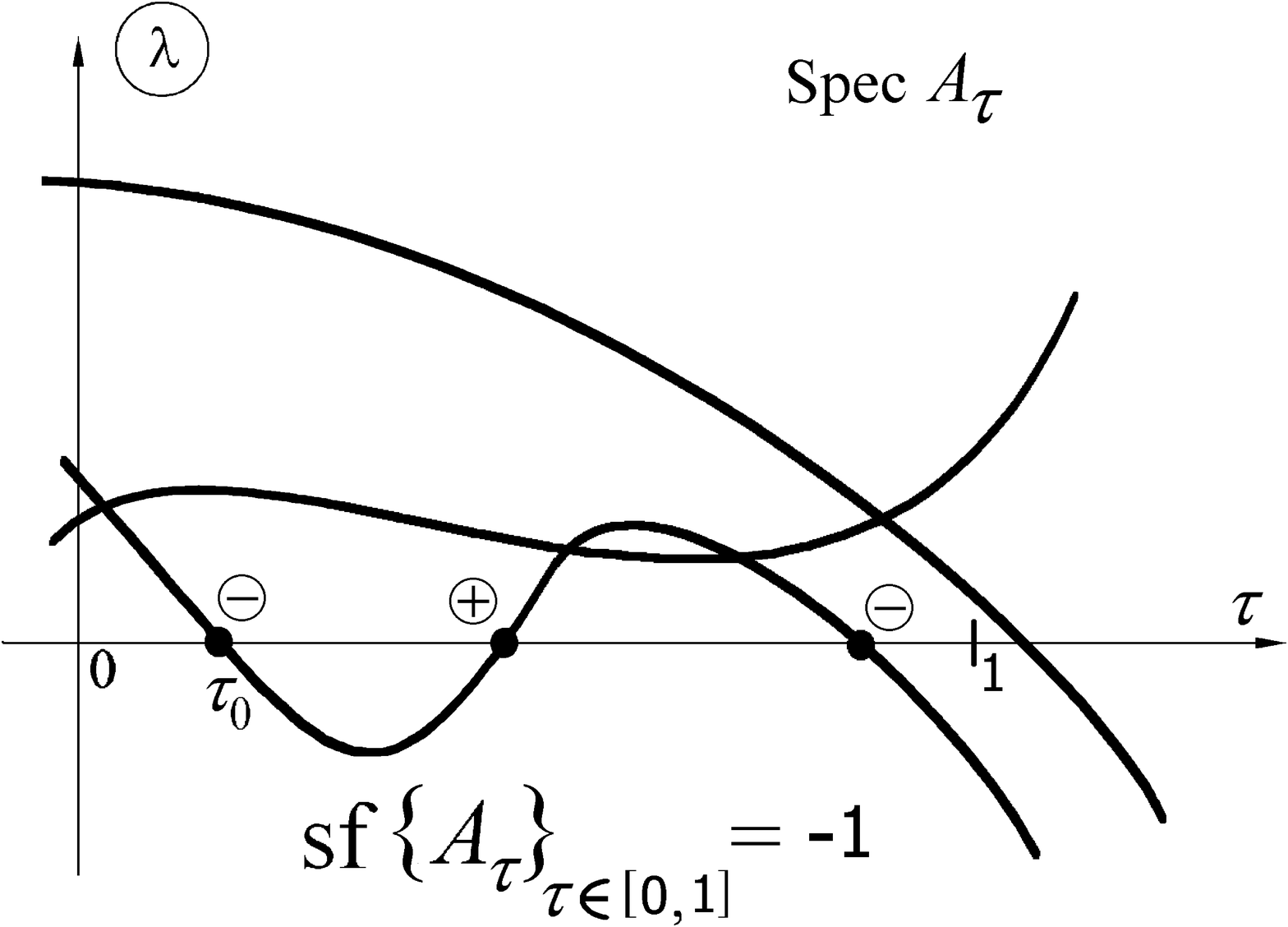}{fi14}{Spectral flow}

A continuous family of elliptic self-adjoint operators has the following
property. For an arbitrary $\tau'\in [0,1]$, there exists a number
$\lambda_{\tau'}$  that is not an eigenvalue of
$A_\tau$ in  an $\varepsilon$-neighborhood of
$\tau'$.
 This follows from the
fact that the spectrum of $A_\tau$ is discrete.

 Now we choose a finite partition
\begin{equation}
\label{dec1} 0=\tau_0<\tau_1<...<\tau_N=1
\end{equation}
of the interval $[0,1]$ and  numbers
$\{\lambda_i\}_{i=0,N-1}$,  referred to as
weights, such that $\lambda_i$ is not an eigenvalue of the
family $A_\tau$ on the interval $[\tau_i,\tau_{i+1}]$. We  also assume
that $\lambda_0=\lambda_N=0$.

\begin{definition}{
The \emph{spectral flow} of the family $\{A_\tau\}_{\tau\in[0,1]}$ is the
 number
\begin{equation}
{\rm sf}\{A_\tau\}_{\tau\in [0,1]}=\sum_{i=1}^{N-1}\ind
(\Pi_{\lambda_{i}}(A_{\tau_i}),\Pi_{\lambda_{i-1}}(A_{\tau_{i}})), \label{ala1}
\end{equation}
where $\Pi_{\lambda}(A)$ is the spectral projection of a self-adjoint
operator $A$ corresponding to eigenvalues greater than or equal to
$\lambda$ and
$$
\ind(\Pi_{\lambda}(A),\Pi_{\mu}(A))={\rm sgn}(\mu-\lambda) \dim
\Lambda_{\lambda,\mu},
$$
is the relative index of two projections. Here $\Lambda_{\lambda,\mu}$
 is the spectral subspace of $A$ corresponding to
eigenvalues in the interval
$[\min\{\lambda,\mu\},\max\{\lambda,\mu\})$.}
\end{definition}

One can show that the spectral flow~\eqref{ala1} is well defined,
i.e.,  is
independent of  the choice of the partition
$\{\tau_i\}$ and the weights $\{\lambda_i\}$.

Let us now return to our original problem and consider a homotopy $D_\tau,$
$\tau\in[0,1]$, of elliptic operators on a manifold with boundary.  The
corresponding family of tangential operators will be denoted by $\{A_\tau\}$.
We shall now give a formula for the difference of  indices at the endpoints of
the homotopy. It turns out that the difference is equal to the spectral flow of
the family of tangential operators.

\begin{theorem}
\label{sft} \emph{(The spectral flow theorem.)} One has
\begin{equation}
\label{sfff} \ind (D_0,\Pi_+(A_0))-\ind (D_1,\Pi_+(A_1))= {\rm
sf}\{A_t\}_{t\in[0,1]}.
\end{equation}
\end{theorem}

\begin{remark}
This result does not contradict the homotopy invariance of the index,
since  we
change not only the operator~$D$, but also the spaces where
the index is computed.
\end{remark}

\begin{proof}

1.  One can readily obtain~\eqref{sfff} if the family $A_\tau$ is invertible
for all $\tau$. Indeed, in this case the right-hand side of~\eqref{sfff}
is zero by the definition of the spectral flow. Let us
verify that the indices on the left-hand side
are  equal. To this end, we note that the family of
nonnegative spectral projections is smooth. Consider the Cauchy problem
\[
\stackrel{.}{U}_\tau=\left[ \stackrel{.}{\Pi }_{+}\left( A_\tau\right) ,\Pi
_{+}\left( A_\tau\right) \right] U_\tau,\quad U_0=\operatorname{Id}.
\]
 One can readily  verify that
 the solution $U_\tau$ is a unitary elliptic operator
 specifying an isomorphism
$$
U_\tau:\im \Pi_+(A_0)\longrightarrow \im\Pi_+(A_\tau)
$$
of the subspaces defined by the pseudodifferential
projections $\Pi_+(A_0)$ and $\Pi_+(A_\tau)$. The composition
$$
\left(
\begin{array}{cc}
1 & 0 \\
0 & U^{-1}_\tau
\end{array}
\right)(D_\tau,\Pi_+(A_\tau))^t: C^\infty(M,E) \longrightarrow
C^\infty(M,F)\oplus \im \Pi_+ (A_0)
$$
has the same index as the  original problem
$(D_\tau,\Pi_+(A_\tau))$. On the other hand, the composition acts in
 spaces independent of $t$. Thus, by the standard index
invariance, its index does not change.

2. Now consider  the general case
in which the operators of the tangential family $A_\tau$ may
 be noninvertible. To this end, we choose some partition
\eqref{dec1}. Using the argument  in the previous part of
the proof,  on the first  interval
$[\tau_0,\tau_1]$ we obtain
$$
\ind (D_0,\Pi_{\lambda_{0}}(A_{\tau_0}))= \ind
(D_{\tau},\Pi_{\lambda_{0}}(A_{\tau})), \quad \tau \in [\tau_0,\tau_1].
$$
Considering this equation for $\tau=\tau_1$ and replacing the projection
$\Pi_{\lambda_0}(A_{\tau_1})$ by $\Pi_{\lambda_1}(A_{\tau_1})$ (they differ by a finite-dimensional projection), we
obtain\footnote{To  justify this and subsequent index
computations,  one uses the fact that
the relative index $\ind(P,Q)$ of projections  coincides
with the index of the Fredholm operator $Q:\im P\to \im Q.$}
$$
\ind (D_0,\Pi_{\lambda_{0}}(A_{\tau_0}))= \ind
(D_{\tau_1},\Pi_{\lambda_1}(A_{\tau_1}) )+ \ind
(\Pi_{\lambda_1}(A_{\tau_1}),\Pi_{\lambda_{0}}(A_{\tau_{1}})).
$$
A similar modification of $\ind (D_{\tau_1},\Pi_{\lambda_1}(A_{\tau_1}))$
at $\tau_2$ gives
\begin{multline*}
\ind (D_0,\Pi_{\lambda_{0}}(A_{\tau_0}))=  \ind
(D_{\tau_2},\Pi_{\lambda_2}(A_{\tau_2}) ) \\
+\ind (\Pi_{\lambda_1}(A_{\tau_1}),\Pi_{\lambda_{0}}(A_{\tau_{1}})) +
\ind(\Pi_{\lambda_2}(A_{\tau_2}),\Pi_{\lambda_{1}}(A_{\tau_{2}})).
\end{multline*}
Proceeding similarly at the subsequent points $\tau_3,\tau_4,\ldots$, we
obtain the desired  equation
$$
\ind (D_0,\Pi_+(A_0))=\ind (D_1,\Pi_+(A_1))+
{\rm sf}\{A_\tau\}_{\tau\in[0,1]}.
$$
\end{proof}

\subsection{A theorem on index decompositions}

 It follows from the spectral flow theorem
in the previous subsection  that the index of the
spectral boundary value problem $(D,\Pi_+(A))$ is uniquely determined by
the principal symbol $\sigma \left( D\right) $ and the tangential
operator $A$. There arises a natural question: Is it possible
to decompose the index of the spectral boundary value problem
 into the sum
\begin{equation}
\limfunc{ind}\left( D,\Pi_+(A)\right) =f_1\left( \sigma \left( D\right)
\right) +f_2\left( A\right)  \label{labb7}
\end{equation}
of a homotopy invariant $f_1\left( \sigma \left( D\right) \right) $ of
the principal symbol of the operator and a functional $f_2\left(A\right)$
of the tangential operator? If this representation is possible,
 how to find it?

We  shall sometimes refer to Eq.~\eqref{labb7} as
 an ``index decomposition."

\begin{remark}
{The functional $f_2$ is not  a homotopy invariant of
the tangential operator in general. Indeed, the set of operators
with  a  given principal symbol
$\sigma(D)$ contains operators with  arbitrary index. Thus,
$f_2$ can take infinitely many values. A more precise analysis shows that
$f_2$ is a homotopy invariant of the corresponding spectral projection.}
\end{remark}

It turns out that there is an \emph{obstruction} to the index
decomposition.

Indeed, suppose that \eqref{labb7} is valid. Consider a homotopy $D_\tau$
of elliptic operators such that the homotopy of their tangential
operators $A_\tau$ is periodic: $A_0=A_1$. We claim that in this
case the indices of the spectral boundary value problems for $D_0$ and
$D_1$ are equal. Indeed,
$$
f_1(\sigma(D_0))=f_1(\sigma(D_1)),
$$
since the symbols are homotopic,
and
$$
f_2(A_0)=f_2(A_1),
$$
since the tangential operators $A_0$ and $A_1$ coincide by  assumption.

On the other hand, by virtue of the spectral flow theorem, the difference
of  indices at the endpoints of the homotopy is equal to the
spectral flow of the periodic family of tangential operators. Thus, we
obtain the following result.

\begin{proposition}\label{lem1}
If the index decomposition is valid, then for an arbitrary homotopy
of tangential operators $A_\tau$, $\tau\in\mathbb{S}^1$, one has
\begin{equation}
{\rm sf}\{A_\tau\}_{\tau\in \mathbb{S}^1} =0. \label{labb8}
\end{equation}
\end{proposition}

It is well known that there exist periodic families of elliptic self-adjoint
operators with  nontrivial spectral flow~(\ref{labb8}). (Simple examples
can be found in~\cite{Sav1}.) Therefore, \emph{the index
decomposition~\eqref{labb7} does not exist on the space of all elliptic
operators.}

In other words, to  achieve~\eqref{labb7}, one has to consider  subspaces
rather than the entire space of elliptic operators.  Using this idea,  one can
prove a result similar to Proposition~\ref{lem1}, where one considers only
homotopies of tangential operators within some  given class of operators.

Namely, let $\Sigma$ be a subspace  of the space of all elliptic Hermitian
symbols acting in the restriction of the bundle $E$ to the boundary.  In the
space ${\rm Ell}(M,E,F)$ of elliptic operators on $M$ acting between the spaces
of sections of vector bundles $E$ and $F$, we consider the subspace ${\rm
Op}(\Sigma)$ of operators such that  the symbols of the corresponding
tangential operators belong to $\Sigma$:
$$
{\rm Op}(\Sigma)=\Bigl\{D\in{\rm Ell}(M,E,F)\;\;\bigl|\;\;\ \sigma(A)\in
\Sigma\Bigr\}.
$$

We  assume that $\Sigma$ is nondegenerate in the following
sense: the natural mapping ${\rm Op}(\Sigma)\to \Sigma$ taking an
elliptic operator on  the manifold with boundary to the
principal symbol of its tangential operator is surjective. In other
words, every element of $\Sigma$ can be realized as the symbol of
the tangential operator for some elliptic operator on $M$.
\begin{definition}{\em
{\em The class ${\rm Op}(\Sigma)$ \emph{admits an index
decomposition} if there exist two functionals }
$$
f_{1,2}: {\rm Op}(\Sigma)\longrightarrow \mathbb{R}
$$
{\em such that}
\begin{enumerate}
{\em \item the first functional is a homotopy invariant of the principal symbol
of the operator, i.e., $f_1(D)=f_1(\sigma(D))$; \item the second functional is
determined by the tangential operator, i.e.,\\ $f_2(D)=f_2(A)$; \item  for
$D\in {\rm Op}(\Sigma)$, one has
$$
\ind (D,\Pi_+(A))=f_1(\sigma(D))+f_2(A).
$$
}
\end{enumerate}
}
\end{definition}
We  shall state
 a \emph{necessary and sufficient condition} for the
existence of a decomposition \eqref{labb7} in terms of the following condition
on the class $\Sigma$ of symbols.
\begin{definition}
{ The class ${\rm Op}(\Sigma)$ is said to be \emph{admissible} if
for an arbitrary periodic family $\{A_\tau\}_{\tau\in \mathbb{S}^1}$ of
elliptic self-adjoint operators on $\partial M$ one has
$$
{\rm sf}\{A_\tau\}_{\tau\in \mathbb{S}^1}=0
$$
provided that $\sigma(A_\tau)\in \Sigma$ for all $\tau$.}
\end{definition}

\begin{theorem}\label{thh1}
\emph{(The index decomposition theorem.)} The  class ${\rm Op}(\Sigma)$ admits
an index decomposition if and only if  it is admissible.
\end{theorem}

\begin{proof}
Necessity  can be proved  by analogy with Proposition \ref{lem1}. The proof  of sufficiency can be found
in~\cite{SaScS1}.
\end{proof}

The admissibility condition can be verified effectively.  Indeed, the principal
symbol of an elliptic self-adjoint operator $A$ on a manifold $X$ (in our
case, $X=\partial M$) defines an element
$$
[\sigma(A)]=[\im \Pi_+\sigma(A)]\in K^0(S^*X)
$$
in the $K$-group, where $\im \Pi_+\sigma(A)\in{\rm Vect}(S^*X)$ is
the subbundle generated by the positive spectral subspaces of the
principal symbol $\sigma(A)$ on the cosphere bundle $S^*X$ (with respect
to some Riemannian metric). Then the spectral flow of a periodic family
$A=\left\{ A_t\right\} _{t\in \Bbb{S}^1}$ of elliptic self-adjoint
operators can be computed by the Atiyah--Patodi--Singer formula
 \cite{APS3}
\begin{equation}
{\rm sf} \left\{ A_t\right\} _{t\in \Bbb{S}^1} =\left\langle
\limfunc{ch}\left[ \sigma(A)\right] \cup \limfunc{Td}\left(
T^{*}X\otimes \Bbb{C}
\right) ,\left[ S^{*}X\times \Bbb{S}^1\right] \right\rangle . \label{omg1}
\end{equation}
Here $\limfunc{ch}\left[ \sigma(A)\right] \in H^{ev}\left( S^{*}X\times
\Bbb{S}^1\right) $ is the Chern character of the element
\[
\left[\sigma(A) \right] =\left[ \im\Pi_+\sigma(A) \right] \in K^0\left(
S^{*}X\times \Bbb{S}^1\right)
\]
defined by the principal symbol of the family, $\limfunc{Td}$ is the
Todd class, and $\langle,\rangle$ is the pairing
 between homology  and cohomology.
\begin{remark}
{The obstruction to the index decomposition given in  Theorem \ref{thh1} has
the following cohomological interpretation. Note that the spectral flow of a
periodic family of tangential operators with symbols in $\Sigma$ defines a
homomorphism
\begin{align*}
{\rm sf}:\pi_1(\Sigma) & \longrightarrow \mathbb{Z},
\\
\{\sigma(A_t)\}_{t\in \mathbb{S}^1} & \longmapsto\left\langle
\limfunc{ch}\left[ \sigma(A)\right] \cup \limfunc{Td}\left( T^{*}X\otimes
\Bbb{C}
\right) ,\left[ S^{*}X\times \Bbb{S}^1\right] \right\rangle,
\end{align*}
of the fundamental group into integers. It vanishes on commutators.
Therefore, by the well-known isomorphism
$H_1(X)=\pi_1(X)/[\pi_1(X),\pi_1(X)]$, the spectral flow defines a
cohomology class
$$
[{\rm sf}]\in H^1(\Sigma,\mathbb{R}).
$$
Now the admissibility condition is equivalent to the vanishing of this
cohomology class. }
\end{remark}

Theorem \ref{thh1} shows that if  the
integral~(\ref{omg1}) is zero for an arbitrary periodic homotopy in
some class $\Sigma$, then the corresponding class of spectral
boundary value problems is admissible and the index admits a
decomposition.

 Now consider examples  in which this
condition is satisfied.

\subsection{Examples}

\hspace{1mm}

\textbf{1. Even subspaces.} { Consider the antipodal involution }
\[
\alpha :S^{*}\partial M\longrightarrow S^{*}\partial M,\quad \alpha
\left( x,\xi \right) =\left( x,-\xi \right)
\]
{of the cosphere bundle of the boundary $\partial M$. For a vector bundle
$E\in{\rm Vect}( M),$  we consider elliptic Hermitian
symbols }
\[
a:\pi^*E|_{\partial M}\longrightarrow \pi^* E|_{\partial M}, \quad \pi:
S^*\partial M\to
\partial M,
\]
{invariant under the involution}
\begin{equation}
a(x,-\xi)=a(x,\xi). \label{labb9}
\end{equation}
{In this case, the contributions to~\eqref{omg1} from antipodal points $\pm
\xi $ in the integral over the cosphere $S_x^{*}\partial M$
corresponding to any $x\in
\partial M$ cancel
provided that $\partial M$ is odd-dimensional. Here we have used the fact that
}$\alpha ${\ preserves (or reverses) the orientation of }$S^{*}\partial M$
depending on the parity of the dimension of $\partial M.${\ Thus,
Eq.~(\ref{labb8}) is satisfied  in the case of an \emph{even-dimensional}
manifold $M$ for operators whose tangential operators have even principal
symbols (see Eq.~\ref{labb9}), and the index decomposition for the
corresponding spectral boundary value problems  is possible. In the following
section, we obtain the index defect formula for this case.}\vspace{1mm}

\textbf{2. Odd subspaces.} { For an \emph{odd-dimensional}
manifold $M$, the antipodal involution }$\alpha ${ preserves
the orientation of }$S^{*}\partial M${. In this case, one should
consider odd symbols $a$ antiinvariant under $\alpha$:}
\[
a(x,-\xi)=-a(x,\xi).
\]
A computation shows that the contributions of  antipodal points to~\eqref{omg1}
cancel modulo a form lifted from the base of the cosphere bundle. Therefore,
the integral is  zero, and this class of boundary value problems has an index
decomposition. This ``odd" case is largely analogous to the even case. Some new
phenomena appear in this case. We touch on this theory only briefly at the end
of  Subsection~\ref{jaja2}. For a detailed exposition, we refer the reader to
\cite{SaSt2}.\vspace{1mm}

\textbf{3. Manifolds  whose boundary is a covering. } Suppose that the boundary is a smooth $n$-sheeted covering
$$
\pi:\partial M\longrightarrow X.
$$
We consider  the class of operators adapted to this
structure in the sense that their
tangential operators  are lifted from the base of
the covering. The  lift is defined by the local
diffeomorphism $\pi$.

Let us compute the spectral flow of  a periodic family of
tangential operators $\{A_\tau\}_{\tau\in \mathbb{S}^1}$. By assumption, this
family is the  lift of some family
$\{{\widetilde{A}}_\tau\}_{\tau\in \mathbb{S}^1} $ of elliptic self-adjoint
operators from the base $X$.  Since
formula~\eqref{omg1} is local, we obtain
$$
{\rm sf}\{A_\tau\}_{\tau\in \mathbb{S}^1} = n{\rm
sf}\{\widetilde{A}_\tau\}_{\tau\in \mathbb{S}^1}\in n\mathbb{Z}.
$$
This is zero as a residue modulo $n$. Therefore, the assumption of Theorem
\ref{thh1} is satisfied modulo $n$, and the index of the corresponding boundary
value problems as a residue modulo $n$ admits the desired
decomposition.\footnote{ We use a theorem similar to Theorem \ref{thh1}, where
one considers the index modulo $n$ and the spectral flow modulo $n$; this
result  can be proved by the same method.} The index defect in this situation
will be studied in Section~\ref{seczn}.

\section{Index defects for problems with parity conditions}
In this section, we describe index defects for spectral boundary
value problems in even subspaces (see Example~\ref{ex1} in the
introduction). The methods of $K$-theory of operator algebras are very
effective in this case.  In the framework of this
approach, the index defect appears naturally as a dimension type
functional of subspaces  determined by spectral
projections.

In this section, we first define a dimension-type functional of subspaces and
then prove the defect formula.

\subsection{The dimension functional for even subspaces}
\begin{definition}{\em
\label{defa}\emph{A pseudodifferential operator }$${P}:C^\infty \left(
X,E\right)\to C^\infty(X,E) $$\emph{is said to be} even \emph{if its principal
symbol }$\sigma(P)$\emph{\ is invariant under the antipodal involution:}
\begin{equation}
\sigma(P)=\alpha ^{*}\sigma(P),\qquad \alpha :S^{*}X\rightarrow S^{*}X,\;\;\;
\alpha \left( x,\xi \right) =\left( x,-\xi \right) . \label{omg}
\end{equation}
}
\end{definition}

\begin{proposition}
\label{prok} Let $A$ be an even-order elliptic self-adjoint
differential operator. Then the spectral projection $\Pi_{+}(A)$ is
even.
\end{proposition}

\begin{proof}
The principal symbol of a differential operator of order $n$ has the
property
\[
\alpha ^{*}\sigma \left( A\right) =\left( -1\right) ^{n}\sigma \left( A\right).
\]
Since
$$
\Pi_+(A)=\frac{A+|A|}{2|A|}
$$
(for an invertible $A$), we see that the principal symbol of
$\Pi_+(A)$ is given by
$$
\sigma(\Pi_+(A))=\Pi_+(\sigma(A)),
$$
where $\Pi_+(\sigma(A))$ is the orthogonal projection on the
nonnegative spectral bundle of the symbol $\sigma(A)$. (Here we use
 the following result  due to Seeley
\cite{See5}: the principal symbol of a function of a self-adjoint
operator is equal to the same function of the symbol.)

The last formula  gives the desired
equality
\[
\sigma\left(\Pi_{+}\left( A\right)\right) =\alpha
^{*}\sigma\left(\Pi_{+}\left( A\right)\right)
\]
for even-order operators.
\end{proof}

We denote the set of even pseudodifferential projections on $X$ by
$\widehat{\limfunc{Even}}\left( X\right)$, and the Grothendieck group of the
semigroup of homotopy classes of even projections  will be denoted by $K\bigl(
\widehat{\limfunc{Even}}\left( X\right) \bigr)$. Let us give an alternative
description of this group.

Consider the algebra $\Psi_{ev }(X)$ of scalar even pseudodifferential
operators of order zero. The Grothendieck group of the semigroup of
homotopy classes of projections in matrix algebras over $\Psi_{ev}(X)$ is
denoted by $K_0(\Psi_{ev}(X))$ and called the \emph{even $K$-group of
the algebra} $\Psi_{ev}(X)$ (e.g., see \cite{Bla1}).
\begin{lemma}
One has
$$
K_0(\Psi_{ev }(X)) \simeq K\bigl( \widehat{\limfunc{Even}} \left( X\right)
\bigr).
$$
\end{lemma}

\begin{proof}
The proof is immediate from the definitions of the Grothendieck group $K\bigl(
\widehat{\limfunc{Even}} \left( X\right) \bigr)$ and the $K$-group of  an
algebra.
\end{proof}

Let us compute the $K$-group  of $\Psi_{ev }(X)$.
\begin{theorem}\label{tma}
 If $X$ is odd-dimensional, then there is an isomorphism
\begin{equation}
\left(\mathbb{Z}\oplus K(X)\right)\otimes \mathbb{Z}\left[\frac 1
2\right]\longrightarrow K_0({\Psi}_{ev }(X))\otimes \mathbb{Z}\left[\frac 1
2\right]. \label{dima}
\end{equation}
Here the mapping takes an integer $k\in \mathbb{Z}$ to a projection of
rank $k$ and a vector bundle $E\in\Vect(X)$ to a projection
defining $E$ as a subbundle  of some trivial bundle. By
$\Bbb{Z}\left[ 1/2\right] $ we denote the ring of dyadic
rationals.
\end{theorem}

\begin{proof}
1.  Let $\overline{\Psi}_{ev}(X)$ be the
$C^*$-algebra obtained as closure  of  ${\Psi}_{ev}(X)$ in the $L^2$-norm.  The closure does not change the $K$-group:
$$
K_0(\Psi_{ev }(X))\simeq K_0(\overline{\Psi}_{ev }(X)).
$$

2. The algebra $\overline{\Psi}_{ev }(X)$  contains
the ideal $\mathcal{K}$ of compact operators, and one has
the exact sequence of algebras
\begin{equation}
0\longrightarrow \mathcal{K} \longrightarrow \overline{\Psi}_{ev }(X)
\longrightarrow C(P^*X)\longrightarrow 0 \label{into}.
\end{equation}
Here the projection takes each operator to its principal symbol. We also  treat
even symbols as continuous functions on the projectivization
$P^*X=S^*X/\mathbb{Z}_2$ of the cosphere bundle. The sequence is well defined
by virtue of the well-known norm estimates
$$
\inf_{K\in \mathcal{K}} \|D+K\|_{L^2(X)\to L^2(X)}=\sup_{(x,\xi)\in
S^*X}|\sigma(D)(x,\xi)|.
$$
Furthermore, the sequence \eqref{into}
induces the sequence
$$
\to K_1(C(P^*X))\stackrel{\ind}\longrightarrow
K_0(\mathcal{K})\longrightarrow K_0(\overline{\Psi}_{ev
}(X))\stackrel{\rm smbl}\longrightarrow K_0(C(P^*X))\to
K_1(\mathcal{K})\to
$$
of $K$-groups, which is in our case reduced to
\begin{equation}
\label{pluss}K^1(P^*X)\stackrel{\ind}\longrightarrow \mathbb{Z}\longrightarrow
K_0(\overline{\Psi}_{ev }(X))\stackrel{\rm smbl}\longrightarrow K(P^*X)\to 0.
\end{equation}
Here we have substituted the well-known computations
$K_0(\mathcal{K})=\mathbb{Z}$ and
$K_1(\mathcal{K})=0$ and replaced the $K$-groups of the
function algebra  $C(P^*X)$ by the topological
$K$-groups of  $P^*X$.

Let us describe the mappings in~\eqref{pluss}. The first mapping takes an
elliptic even symbol to the index of the corresponding operator. The second
mapping takes a positive integer $k$ to a projection of rank $k$. (Such a
projection is even, since its symbol is zero.) Finally, the mapping ${\rm
smbl}$ takes a pseudodifferential projection $P$ to the range $\im
\sigma(P)\in\Vect(P^*X)$ of its principal symbol
 treated as a vector bundle over $P^*X$.

For the existence of a functional $d$ making the diagram
\begin{equation}
\xymatrix{ \quad K_0(\mathcal{K})=\mathbb{Z}
\ar[d] \ar[rd]& \\
K_0({\overline{\Psi}}_{ev }(X))\ar[r]^{\qquad d} &\mathbb{R} }
\end{equation}
commute, it is necessary that the vertical arrow be a monomorphism
or, equivalently, the index mapping in \eqref{pluss} be zero. This
condition is satisfied if $X$ is odd-dimensional. Indeed, it is well
known (e.g., see \cite{Pal1}) that the index of operators with even
principal symbol is zero on such manifolds. Therefore, the sequence can
finally be rewritten as
\begin{equation}\label{exa5}
0\longrightarrow \mathbb{Z}\longrightarrow K_0(\overline{\Psi}_{ev
}(X))\longrightarrow K(P^*X)\longrightarrow 0.
\end{equation}

3. Let us slightly simplify this sequence  further. To
this end,  we note that the natural projection
$p:P^{*}X\rightarrow X$ for an odd-dimensional $X$ induces an isomorphism in
$K$-theory modulo $2$-torsion. More precisely, the following
 is valid.
\begin{proposition}
\emph{\cite{Gil7}}\label{pro7} The projection $P^*X\to X$ induces an
isomorphism
\[
p^{*}:K^*\left( X\right) \otimes \Bbb{Z}\left[ \frac 12\right] \longrightarrow
K^*\left( P^{*}X\right) \otimes \Bbb{Z}\left[ \frac 12\right].
\]
\end{proposition}

\begin{proof}[Sketch of Proof] One first verifies the statement over a point
$x\in X$: here $K^1(\mathbb{RP}^{2n})=0$ and $\widetilde
K^0(\mathbb{RP}^{2n})=\mathbb{Z}_{2^n} $, which shows that
$$
p^{*}:K^*\left(\{x\}\right) \otimes \Bbb{Z}\left[ \frac 12\right]
\longrightarrow K^*\left( P^{*}_xX\right) \otimes \Bbb{Z}\left[ \frac 12\right]
$$
is an isomorphism. Now the isomorphism for the entire space can be proved by
the Mayer--Vietoris principle (\cite{BoTu1}).
\end{proof}

Taking a tensor product of the sequence \eqref{exa5}  by
the ring of dyadic rationals (which does not
violate the exactness), we obtain  the exact
sequence
$$
0\longrightarrow \mathbb{Z}\left[\frac 1 2\right]\longrightarrow K_0(\overline{\Psi}_{ev
}(X))\otimes \mathbb{Z}\left[\frac 1 2\right]\longrightarrow K(X)\otimes
\mathbb{Z}\left[\frac 1 2\right]\longrightarrow 0.
$$

4. This exact sequence has a splitting mapping
\[
K^0\left( X\right) \otimes \mathbb{Z}\left[ \frac 12\right] \longrightarrow
K\bigl( \widehat{\limfunc{Even}}\left( X\right) \bigr) \otimes \mathbb{Z}\left[
\frac 12\right]
\]
taking a vector bundle to  a projection on the space
of its sections in $C^\infty(X,\mathbb{C}^N)$ for
 sufficiently large  $N$.

5. The splitting obviously gives the desired isomorphism \eqref{dima}.
\end{proof}

The first component of the isomorphism \eqref{dima} will be called the
\emph{dimension functional} for even pseudodifferential
projections. It is useful to restate
Theorem~\ref{tma} in the following way.
\begin{theorem}
\emph{\cite{SaSt1} (A theorem on the dimension
functional.)}\label{thmm1} On an odd-dimensional manifold $X$,
there exists a unique group homomorphism
\[
d:K\bigl( \widehat{\limfunc{Even}}\left( X\right) \bigr) \longrightarrow
\Bbb{Z}\left[ \frac 12\right]
\]
with the property
\begin{equation}
d\left( P\right) ={\rm rank} P \label{labb10}
\end{equation}
{for a finite rank projection }$P$ under the normalization
\[
d\left( P_F \right) =0,
\]
for all projections $P_F:C^\infty(X,E)\longrightarrow C^\infty(X,E)$ on
the space of sections of a subbundle $F\subset E$. Moreover, for an
arbitrary even projection $P$ and for a sufficiently large $N$ the
projection $2^NP$ is homotopic to  the direct sum of a
projection on the space of sections of a subbundle and a finite rank
projection.
\end{theorem}

\subsection{The index defect formula}
\label{jaja2}

Let us apply the dimension functional $d$ of even projections
to the theory of spectral boundary value problems.

{\bf 1. The dimension $d$ as an index defect.} We shall consider elliptic
operators $D$ with even tangential operator $A$. In this case, the spectral
projection $\Pi_+(A)$ is even as well (see Proposition \ref{prok}).
\begin{remark}
Clearly, this condition cannot be satisfied for first-order differential
operators. Thus, in this subsection we consider  the more general class of
operators that are pseudodifferential far from the boundary and  have the
form~\eqref{deco2} with a first-order pseudodifferential operator $A$ on
$\partial M$ near the boundary. For this class of operators, the spectral
boundary value problems are  also well defined and have the Fredholm property.
Operators of this form are considered in detail in~\cite{Hor3}.
\end{remark}

\begin{lemma}
The sum
$$
\ind(D,\Pi_+(A))+d(\Pi_+(A))
$$
is a homotopy invariant of
the operator $D$.
\end{lemma}

\begin{proof}
Consider a homotopy $D_\tau,$ $\tau\in[0,1],$ and  the corresponding homotopy
$A_\tau$ of tangential operators.  By the homotopy invariance of the dimension
functional and  property~\eqref{labb10}, it follows that the variation of the
dimension is equal to the spectral flow
$$
d(\Pi_+(A_1))-d(\Pi_+(A_0))={\rm sf} \{A_\tau\}_{\tau\in [0,1]}.
$$
(The equality can be obtained  by analogy with the
proof of the theorem on the spectral flow with
the use of a partition of the  interval $[0,1]$ and
some choice of weights on the  intervals of the
partition.) A similar formula holds for the variation of the index:
$$
\ind(D_1,\Pi_+(A_1))-\ind(D_0,\Pi_+(A_0))=-{\rm sf} \{A_\tau\}_{\tau\in [0,1]}.
$$
Combining the two expressions, we obtain the desired homotopy invariance
of the sum  indicated in the lemma.
\end{proof}

This homotopy invariant  will be denoted by
$$
\inv D\stackrel{\operatorname{def}}=\ind(D,\Pi_+(A))+d(\Pi_+(A)).
$$
Let us give a topological formula for this analytic invariant. This
 can be done by generalizing the Atiyah--Singer topological
index to the  case of spectral boundary value
problems in even subspaces.

\textbf{2. The topological index.} It turns out that the principal
symbol $\sigma(D)$ has a natural extension to the double
\[
2M=M\bigcup\nolimits_{\partial M}M
\]
of the manifold $M$. The double is obtained by gluing two copies of $M$ along
the common boundary.

To construct the desired extension,  we
take the symbol $\sigma(D)$ on the first copy of $M$ and the
symbol $\alpha ^{*}\sigma \left( D\right)$  on the second copy.
Here
$$
\alpha:S^*M\longrightarrow S^*M
$$
is the antipodal involution. Then in a
neighborhood of the boundary the symbols $\sigma(D)$ and $\alpha
^{*}\sigma(D)$  have the form
$$
i\tau+a(x,\xi) \quad \text{and} \quad -i\tau+a(x,\xi),
$$
respectively. They are
taken into one another  by the coordinate
transformation
$$
x\to x, \quad t\to -t
$$
in a neighborhood of the boundary. Therefore,
taken together, they define an elliptic symbol $\sigma(D)\cup \alpha ^{*}\sigma(D)$ on the double of $M$. The
difference element of this symbol  will be denoted by
$$
[\sigma(D)\cup \alpha ^{*}\sigma(D)]\in K_c(T^*2M).
$$
Then the \emph{topological index} of $D$ is defined as half the
Atiyah--Singer topological index of  this
element on the double:
$$
\ind_t[\sigma(D)]\stackrel{\operatorname{def}}= \frac{1}2\ind_t [\sigma(D)\cup
\alpha ^{*}\sigma(D)].
$$

{\bf 3. The index defect formula.}
\begin{theorem}
{\em \cite{SaSt1}} \label{thbvp} Let $D$ be an elliptic operator on an
even-dimensional manifold with even tangential operator $A$. Then
\begin{equation}
\inv D = \ind_t[\sigma(D)]. \label{quest}
\end{equation}
\end{theorem}

\begin{proof}
For sufficiently large $N$, Theorem \ref{thmm1}  in the previous subsection
gives  a homotopy $\pi_\tau$, $\tau\in[0,1]$, of the direct sum of $2^N$ copies
of the symbol $\sigma(\Pi_+(A))$  to a projection on some subbundle $E_0\subset
2^NE|_{\partial M} $. One can lift this homotopy of projections to a homotopy
$\sigma(D_\tau)$ of elliptic symbols with the properties
$$
\sigma(D_0)=2^N\sigma(D), \quad \sigma(\Pi_+(A_\tau))=\pi_\tau.
$$
By the homotopy invariance of both sides of Eq.~\eqref{quest}, it suffices to
prove the equality only for the symbol $\sigma(D_1)$ obtained at the end of the
homotopy. This symbol in a neighborhood of the boundary depends only on the
absolute value of the covector $\xi$. Thus, one can consider an elliptic
operator $D_1$ with this symbol and the corresponding spectral boundary value
problem $(D_1,\Pi_+)$ such that the spectral projection $\Pi_+$ is actually a
projection on the space of sections of the subbundle $E_0$. The spectral
boundary value problem in this case is a classical boundary value problem; the
index defect (i.e. the dimension functional) is zero, and Eq.~\eqref{quest}
follows from the Atiyah--Bott formula.
\end{proof}

\begin{remark}
(On the dimension functional of odd subspaces.) A similar
functional was constructed in \cite{SaSt2} on the space of ``odd
projections" with principal symbols satisfying
$$
\alpha^*\sigma(P)=1-\sigma(P).
$$
Such projections  arise as spectral projections of
odd-order differential operators (cf.\ Proposition \ref{prok}). However,
the methods of $K$-theory cannot be applied directly to odd projections,
since odd symbols do not form an algebra. Moreover, the geometry of odd
projections differs from the geometry of even ones: for
example, an odd projection on a manifold of dimension $2k$ can
 act  only in the space of sections of a vector
bundle  whose dimension  is a multiple of~$2^{k-1}$ (see~\cite{SaSt2}). In the
cited paper, we constructed a dimension functional and proved the defect
formula. Let us only mention  that the topological index in the odd case
 can be obtained if on the double $2M$ one considers the
symbol $\sigma(D)\cup \alpha ^{*}\sigma(D)^{-1}$.
\end{remark}

\begin{remark}
One can show \cite{SaScS4} that the value $d(\Pi_+(A))$ of the dimension
functional and the element $[\sigma(D)\cup \alpha^*\sigma(D)]$ form a
complete system of stable homotopy invariants of the spectral boundary
value problem $(D,\Pi_+(A))$, i.e., classify these problems modulo stable
homotopy.
\end{remark}

\subsection{The dimension functional and the $\eta $-invariant\label{p82}}

The dimension functional for even pseudodifferential projections was
constructed in Theorem \ref{thmm1}. More precisely, the theorem claims
the existence of the functional. In this subsection, we address the
question of an explicit analytic expression for this functional. It turns
out that such a description can be given in terms of the spectral
Atiyah--Patodi--Singer $\eta$-invariant. The reader can find some details
about the $\eta$-invariant in Appendix A.

Gilkey \cite{Gil2} noted that the $\eta$-invariant is rigid within the class of
differential operators satisfying the \emph{parity condition}
\begin{equation}
{\rm ord}A+{\dim} X\equiv 1 ({\rm mod}\, 2), \label{para1}
\end{equation}
which relates the order of the operator to the dimension of the manifold.
Rigidity is understood in the sense that the fractional part of the
$\eta$-invariant is not only spectrally invariant but also homotopy
invariant.

Condition \eqref{para1} coincides with the conditions under which there
exists a functional $d$ of the corresponding spectral projections.
(Recall that in the previous subsections this functional was considered
for even projections on odd-dimensional manifolds and odd projections on
even-dimensional manifolds.) This is not a mere coincidence. In fact, the
$\eta$-invariant of an arbitrary elliptic self-adjoint operator $A$
coincides with the value of the dimension functional on the spectral
projection $\Pi_+(A)$.

\begin{theorem}{\em \cite{SaSt1,SaSt2}}
If an elliptic self-adjoint differential operator satisfies the parity
condition~\emph{(\ref{para1})}, then
\[
\eta \left( A\right) =d\left( \Pi_{+}(A)\right).
\]
\end{theorem}

\begin{proof}
For simplicity, we consider only opperators of even order $2l$.

1. It was shown in Theorem \ref{thmm1} that the dimension $d$ of even
projections can be defined as the unique homomorphism
\[
d:K\bigl( \widehat{\limfunc{Even}}\left( X\right) \bigr) \longrightarrow
\Bbb{R}
\]
with the following two properties: (a) $ d\left( P\right) ={\rm rank} P$
for a finite rank projection $P$, (b) $d\left( P_F \right)=0$ for
projections on the spaces of sections of vector bundles $F\in \Vect(X)$.
Therefore, to prove the theorem it suffices to show that the
$\eta$-invariant defines a similar homomorphism and enjoys the same
properties.

2. Thus, we should consider the $\eta$-invariant for general even
pseudodifferential operators. Unfortunately, the rigidity of the
$\eta$-invariant is lost in this class, since lower-order terms of the
operator contribute to the $\eta$-invariant. However, the invariance
still holds (see Appendix A) in the class of
$\mathbb{R}_*$-\emph{invariant} operators, i.e., operators with complete
symbol having an asymptotic expansion
$$
{\rm smbl} A\sim a_{2l}(x,\xi)+a_{2l-1}(x,\xi)+a_{2l-2}(x,\xi)+\ldots,
$$
where each homogeneous term $a_k(x,\xi)$ of order $k$ has the following
parity with respect to the momentum variables $\xi$:
$$
a_k(x,-\xi)=(-1)^ka_k(x,\xi).
$$
For $\Bbb{R}_{*}$-invariant elliptic self-adjoint operators, the
spectral projection $\Pi_{+}\left( A\right) $ is either even or odd
according to the parity of the order of the operator.

Let us verify property (a) for the $\eta$-invariant (as a homotopy
invariant of the spectral projection). For a nonnegative Laplacian
$\Delta $ with $k$-dimensional kernel, we obtain
\[
\dim \im \Pi_+\left( -\Delta^l \right) =k.
\]
On the other hand, the $\eta $-invariant can be expressed in terms of
the $\zeta$-invariant as
\[
\eta \left( -\Delta^l \right) =k-\zeta \left( \Delta^l \right) .
\]
But the $\zeta $-invariant of the differential operator $\Delta^l $ is
zero in odd dimensions (see the remark following Theorem~\ref{droba} in
Appendix~A). This proves property (a):
\[
\eta \left( -\Delta^l \right) =\dim \im\Pi_{+}\left( -\Delta^l \right)
=d(\Pi_+(-\Delta^l)) .
\]
Similarly, one obtains (b):
\[
\eta \left( \Delta^l \right) =0,\quad \im \Pi_{+}\left( \Delta^l \right)
=C^\infty \left( M,E\right) .
\]
By the characterization of the dimension functional, this proves that the
$\eta $-invariant of an operator $A$ coincides with the dimension $d$ of
the corresponding spectral projection.
\end{proof}

\section{Index defects on twisted $\mathbb{Z}_n$-manifolds}
\label{seczn}

In this section, we consider another geometric situation in which index
defects naturally arise. This situation comes from elliptic theory on
manifolds with singularities. More precisely, we consider so-called
twisted $\mathbb{Z}_n$-manifolds. Such manifolds look like a book with
$n$ sheets near the singularity; the singular set is the edge where the
sheets meet.

The index theorem was stated and proved by Freed and Melrose
\cite{FrMe1} for the special case in which a neighborhood of the singular
set consists of $n$ sheets globally. Note that in this case it is natural
to view the index as a residue modulo $n$, and the index theorem computes
this index-residue.

We consider general twisted $\mathbb{Z}_n$-manifolds. In this case, the
index-residue is no longer a homotopy invariant and an index defect arises.
Intuitively, the proof of the index defect formula follows the scheme explained
in Subsection~\ref{cch6}, so here we give only the main steps of the proof. We
first define a special element in the $K$-theory of a $\mathbb{Z}_n$-manifold
in Subsection~\ref{s61}. Then the topological index is defined as the pairing
with this element. The last subsection contains an application of the defect
formula.

\subsection{Twisted $\mathbb{Z}_n$-manifolds and elliptic operators}
\label{zns}

\hspace{2mm}

\textbf{1. Twisted $\mathbb{Z}_n$-manifolds}
\begin{definition}A{\em
\label{d6} {twisted }$\Bbb{Z}_n$-manifold, \emph{where $n$ is a positive
integer, is a smooth compact manifold $M$ with boundary }$\partial
M$\emph{\ equipped with the structure of the total space of an
$n$-sheeted covering}
\begin{equation}\label{cova}
\pi:\partial M\longrightarrow X
\end{equation}
\emph{over a smooth base $X$ (see Fig.~\ref {fi12}).}}
\end{definition}
\fgr{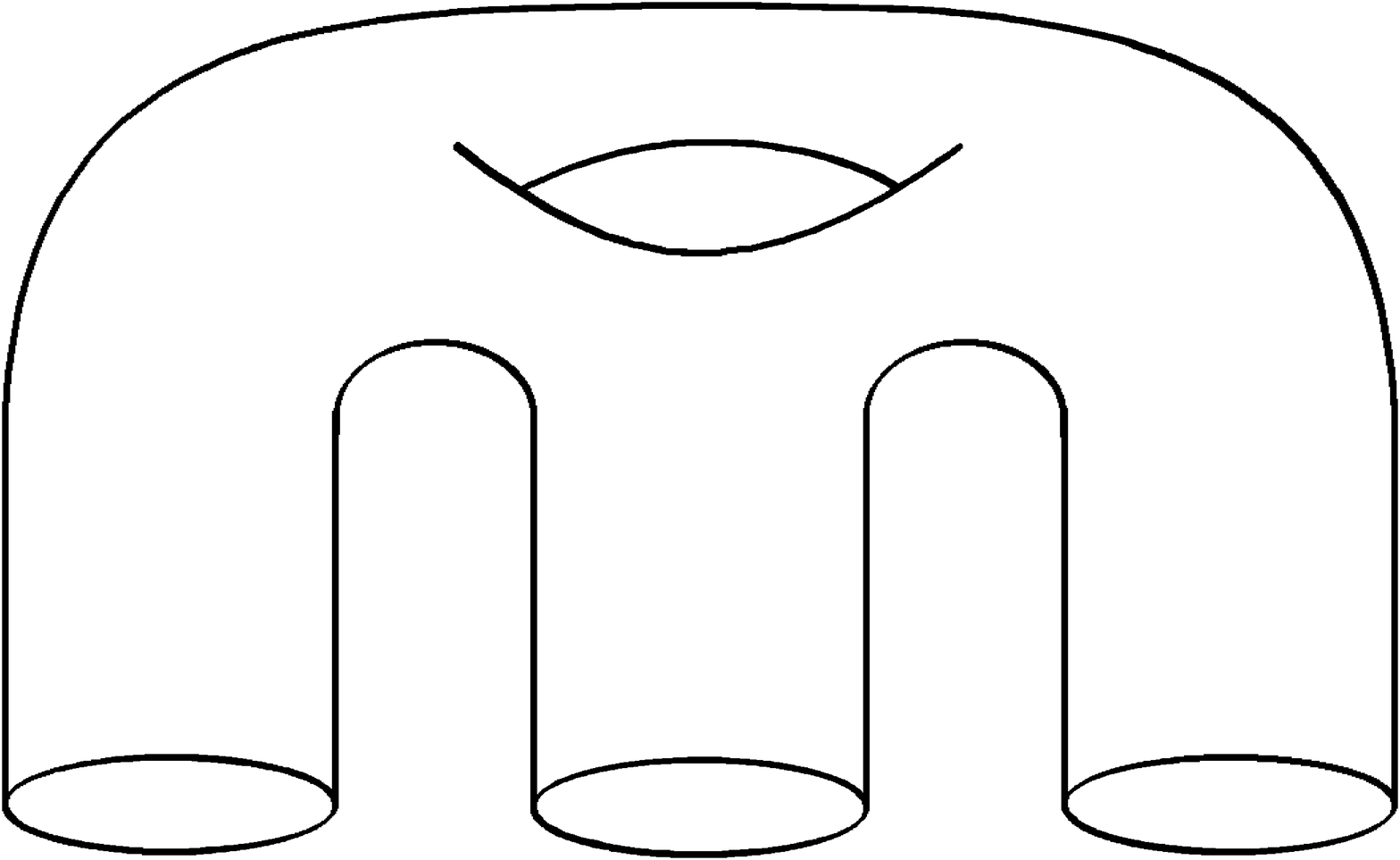}{fi12}{$\Bbb{Z}_3${-manifold}}

Geometrically, a twisted $\Bbb{Z}_n$-manifold $(M,\pi)$ naturally defines
the singular space
\[
\overline{M}^\pi=\left. M\right/\sim,
\]
obtained by the identification of points in each fiber of the covering (see
Fig.~\ref{fi13}); the corresponding equivalence relation $\sim$ is defined as
\begin{equation}
\label{eka} x\sim y \Longleftrightarrow x=y \mbox{ or } \{x,y\in \partial M
\mbox{ and } \pi(x)=\pi(y) \}.
\end{equation}

Sullivan \cite{Sul1} introduced the notion of a $\Bbb{Z}_n$-manifold.
These manifolds correspond to the structure of a \emph{trivial} covering
$\pi$. One of the motivations showing the interest in such manifolds is
that (in the orientable case) the manifold $\overline{M}^\pi$ has a
fundamental class in homology with coefficients in $\Bbb{Z}_n$:
\[
\bigl[ \overline{M}^\pi\bigr] \in H_m\bigl( \overline{M}^\pi,\Bbb{Z}_n\bigr)
,\qquad m=\dim M.
\]
These manifolds with singularities were also used to give a geometric
realization of bordism theory with coefficients in the group
$\mathbb{Z}_n$. Further developments in this direction can be found,
e.g., in~\cite{Bot2}.

For brevity, we frequently omit the word ``twisted" in what follows.

We point out that while most technical constructions of the theory deal
with the smooth model $(M,\pi)$, the most important results, e.g., the
defect formula in Subsection \ref{defzn}, are stated in terms of the
singular space $\overline{M}^\pi$ itself.

{\bf 2. Natural mappings associated with coverings.} Let us recall some
natural mappings induced by the covering~\eqref{cova}. First, one has the
\emph{direct image} mapping
\[
\pi _{!}:\limfunc{Vect}\left( \partial M\right) \longrightarrow
\limfunc{Vect}\left( X\right)
\]
taking a vector bundle $E\in \limfunc{Vect}\left(\partial M\right) $ to the
bundle
\[
\pi _{!}E\in \limfunc{Vect}\left( X\right) ,\quad \left( \pi _{!}E\right)
_x=C^\infty \left( \pi ^{-1}\left( x\right) ,E\right) ,\quad x\in X.
\]
This definition leads to a natural isomorphism
$$
\beta _E:C^\infty \left(\partial M
,E\right) \stackrel{\simeq }{\to} C^\infty \left( X,\pi _{!}E\right)
$$
of section spaces on the total space and the base of the covering. For example,
scalar functions on the total space correspond to sections of the bundle
$\pi_!1\in\Vect(X)$ on the base. (Here $1$ stands for the trivial line bundle.)

This isomorphism enables one to transfer operators acting on $\partial M$
to $X$ and vice versa. More precisely, for an operator
\[
P:C^\infty \left( \partial M,E\right) \longrightarrow C^\infty \left(
\partial M,E\right)
\]
on $\partial M$, by $\pi _{!}P$ we denote its \emph{direct image} given
by the formula
\[
\pi _{!}P=\beta _EP\beta _E^{-1}:C^\infty \left( X,\pi _{!}E\right)
\longrightarrow C^\infty \left( X,\pi _{!}E\right) .
\]
One also has the \emph{inverse image}
\begin{equation}
\pi ^{!}P^{\prime }=\beta _E^{-1}P^{\prime }\beta _E:C^\infty \left(
\partial M,E\right) \longrightarrow C^\infty \left( \partial M,E\right) \label{nona}
\end{equation}
of operators $P^{\prime } $ on $X$.

On twisted $\mathbb{Z}_n$-manifolds, we consider the following class of
operators adapted to the structure of the singularity.

\textbf{3. The class of operators}. For a pair $(M,\pi)$, we consider
elliptic differential operators
\[
D:C^\infty \left( M,E\right) \longrightarrow C^\infty \left( M,F\right) ,
\]
which in a neighborhood of the boundary are lifted from the base of the
covering. More precisely, we suppose that they satisfy the following
assumption.

\emph{Assumption 1}. The restrictions of $E$ and $F$ to the boundary are
lifted from the base:
\[
\left. E\right| _{\partial M}\simeq \pi ^{*}E_0,\quad \left. F\right|
_{\partial M}\simeq \pi ^{*}F_0,\quad \quad E_0,F_0\in
\limfunc{Vect}\left( X\right) ;
\]
moreover, the corresponding isomorphisms are given. For some operator
$$
D_0:C^\infty \left( X\times \left[ 0,1\right) ,E_0\right) \longrightarrow
C^\infty \left( X\times \left[ 0,1\right) ,F_0\right)
$$
on the cylinder with base $X$, the direct image $\left( \pi \times
1\right) _{!}D$ in a neighborhood of the boundary can be inserted in the
commutative diagram
\begin{equation}
\xymatrix{ C^\infty \left( X\times \left[ 0,1\right) ,(\pi\times 1)
_{!}E\right) \ar[d]_\simeq \ar[r]^{\left( \pi \times 1\right) _{!}D}& C^\infty
\left(X\times \left[ 0,1\right) ,(\pi\times 1) _{!}F\right) \ar[d]^\simeq
\\
C^\infty \left( X\times \left[ 0,1\right) ,E_0\otimes \pi _{!}1\right)
\ar[r]^{D_0\otimes 1_{\pi_!1}} & C^\infty \left( X\times \left[0,1\right)
,F_0\otimes \pi _{!}1\right) . } \label{pdown}
\end{equation}
Here we suppose that we have chosen a diffeomorphism
\begin{equation}
U_{\partial M}\simeq \left[ 0,1\right) \times \partial M \label{rava}
\end{equation}
of a collar neighborhood of the boundary and extended the projection $\pi$ to
it as $\pi\times 1$. Then $D_0\otimes 1_{\pi_!1}$ is the operator $D_0$ with
coefficients in the flat bundle $\pi _{!}1$ (e.g., see \cite{APS3}).

We also assume for simplicity that the operator $D_0$ in a neighborhood of the
base of the cylinder $X\times \left[ 0,1\right)$  has the form \eqref{deco2},
i.e.,
\[
\left. D_0\right| _{X\times \left[ 0,\varepsilon \right) }=\Gamma
\left( \frac \partial {\partial t}+A_0\right) ,
\]
where $\Gamma $ is some vector bundle isomorphism and the tangential
operator $A_0$ is a first-order elliptic self-adjoint operator on $X$.

It follows from this assumption that in a neighborhood of the boundary the
operator $D$ has the form
\[
\frac \partial {\partial t}+\pi ^{!}\left( A_0\otimes 1_{\pi_!1}\right)
\]
(up to a vector bundle isomorphism); i.e., the tangential operator $A$ is
equal to $\pi^{!}\left( A_0\otimes 1_{\pi_!1}\right)$.

We note that the classical operators satisfy Assumption 1 in an
appropriate geometric setting. For example, the Hirzebruch operator on an
oriented even-dimensional Riemannian manifold $M$ satisfies the
assumption if (1) in a collar neighborhood of the boundary we take a
metric pulled back from the cylinder $X\times[0,1]$; (2) the covering
$\pi:\partial M\to X$ is oriented. Similar statements hold also for the
Dirac and Todd operators; we leave them to the reader. \vspace{1mm}

{\bf 4. The difference element.} Assumption 1 is closely related to the
above-discussed manifolds with singularities. Indeed, note that the total
space $T^*M$ of the cotangent bundle is a (noncompact)
$\Bbb{Z}_n$-manifold. It follows that the principal symbol of an elliptic
operator $D$ defines a difference element
\[
\left[ \sigma \left( D\right) \right] \in K_c\bigl( \overline{T^{*}M}^\pi\bigr)
\]
in the $K$-group, since the commutative diagram \eqref{pdown} shows that the
restriction of $\sigma(D)$ to the boundary is isomorphic to a symbol lifted
from the base of the covering $T^*M|_{\partial M}\to T^*X\times \mathbb{R}$.

\subsection{The Freed--Melrose index theorem modulo $n$}

This subsection deals with the index theory on $\mathbb{Z}_n$-manifolds
corresponding to trivial coverings $\pi$; i.e., the boundary of the
corresponding smooth manifold $M$ is a disjoint union of $n$ copies of
the base of the covering.

{\bf 1. The index modulo $n$.} It turns out that in this case operators
satisfying Assumption 1 have a nontrivial homotopy invariant. Namely, let
$$
\Mod n\text{-}\ind D\in \mathbb{Z}_n
$$
be the index of the spectral boundary value problem $(D,\Pi_+(A))$
treated as a residue modulo $n$.

\begin{proposition}
The index-residue $\Mod n\text{-}\ind D$ is a homotopy invariant of the
operator $D$.
\end{proposition}

\begin{proof}
Consider a homotopy $\{D_t\}_{t\in [0,1]}$. By the spectral flow theorem,
$$
\ind(D_1,\Pi_{+}(A_1))-\ind(D_0,\Pi_{+}(A_0))=-{\rm sf} \{A_t\}_{t\in[0,1]}.
$$
On the other hand, by assumption, the family $A_t$ of tangential
operators is the direct sum of $n$ copies of the family $\{A_{0,t}\}$ of
tangential operators on $X$. Therefore,
$$
{\rm sf} \{A_t\}_{t\in[0,1]}=n{\rm sf} \{A_{0,t}\}_{t\in[0,1]}\equiv 0 ({\rm
mod} n).
$$
This shows that the spectral flow vanishes as a residue modulo $n$ and proves
the homotopy invariance of the index-residue.
\end{proof}

This homotopy invariant index residue was computed in terms of the
principal symbol by Freed and Melrose. In Subsection~\ref{defzn}, we
obtain a more general formula, and now we briefly recall the
Freed--Melrose formula, just to make the exposition complete.

{\bf 2. The Freed--Melrose theorem.} Consider the category of
$\Bbb{Z}_n$-manifolds with morphisms given by embeddings that take the
boundary to the boundary and the fibers of the coverings to the fibers.
The direct image mapping in $K$-theory extends to this category. More
precisely, for an embedding $f$ of a pair $(M,\pi)$ in $(N,\pi_N)$ there
is a direct image mapping
\[
f_{!}:K_c\bigl( \overline{T^{*}M}^\pi\bigr) \longrightarrow K_c\bigl(
\overline{ T^{*}N}^{\pi_N}\bigr) .
\]
On the other hand, there is a universal space in which one can embed an
arbitrary $\mathbb{Z}_n$-manifold corresponding to a \emph{trivial} covering.
The universal space can be defined as follows. From $\Bbb{R}^L$ we cut away the
union of $n$ disjoint discs. We denote the resulting manifold with boundary by
$M_n$. It can be viewed as a $\Bbb{Z}_n$-manifold, since its boundary consists
of $n$ diffeomorphic spheres. (The diffeomorphisms are given by translations.)
One can readily compute the $K$-group of the cotangent bundle of this space:
\[
K_c\bigl( \overline{T^{*}M_n}^{\pi'}\bigr) \simeq \Bbb{Z}_n.
\]

Freed and Melrose proved the following index theorem.

\begin{theorem}
\emph{\cite{FrMe1} } One has
\[
\func{mod}n\text{-}\limfunc{ind}D=f_{!}\left[ \sigma \left( D\right) \right] ,
\]
where the direct image
\[
f_{!}:K_c\bigl( \overline{T^{*}M}^\pi\bigr) \longrightarrow K_c\bigl(
\overline{ T^{*}M_k}^{\pi'}\bigr) \simeq \Bbb{Z}_n
\]
is induced by an embedding
\[
f:M\longrightarrow M_n.
\]
\end{theorem}
The proof models the $K$-theoretic proof of the Atiyah--Singer index theorem
based on embeddings. The cornerstone of the proof is the statement that the
analytic index is preserved under the direct image mapping, that is, the
diagram
$$
\xymatrix{ K_c(\overline{T^*M}^\pi)
\ar[rr]^{f_!} \ar[dr]_{\func{mod}n\text{-}\limfunc{ind}} & &
K_c(\overline{T^*N}^{\pi_N})\ar[ld]^{\func{mod}n\text{-}\limfunc{ind}} \\
&\mathbb{Z}_n }
$$
commutes for an embedding of $M$ in $N$.

\subsection{The index defect problem on twisted $\mathbb{Z}_n$-manifolds}

In contrast to the case of $\mathbb{Z}_n$-manifolds corresponding to
trivial coverings, considered in the previous subsection, the index
modulo $n$ is no longer a homotopy invariant if the covering is
\emph{nontrivial}. By way of illustration, consider the following simple
example.

On the cylinder $\mathbb{S}^1\times [T_1,T_2]$ with coordinates
$(\varphi,t)$, consider the scalar elliptic operator
$$
D=\frac\partial{\partial t}+ \left(-i\frac\partial{\partial\varphi}+t\right).
$$
Clearly, the tangential operator is lifted from the base of any of the
coverings
$$
\pi_n:\mathbb{S}^1\to \mathbb{S}^1,\quad \varphi\mapsto n\varphi.
$$
On the other hand, if $T_1$ or $T_2$ passes through zero, then the index of the
corresponding spectral boundary value problem jumps by $1$ and is not constant
as a residue modulo $n$ for any $n\ne 1$.

{\bf 1. The index defect.} The example shows that the index-residue is
not homotopy invariant and cannot be computed topologically. The
following theorem gives the index defect in this situation.
\begin{theorem}\cite{SaSt10}
On a twisted $\mathbb{Z}_n$-manifold, the difference
$$
\eta \left( A\right) -n\eta \left( A_0\right) \in \mathbb{R}/n\mathbb{Z}
$$
is an index defect. Here $n$ is the number of sheets of the covering, and
$\eta \left( A\right)$ and $\eta \left( A_0\right)$ are the spectral
Atiyah--Patodi--Singer $\eta $-invariants of the tangential operators $A$
and $A_0$, respectively. In other words, the sum
\begin{equation}
\func{mod}n\text{-}\limfunc{ind}\left( D,\Pi _{+}(A)\right) +\eta \left(
A\right) -n\eta \left( A_0\right) \in \Bbb{R}/n\Bbb{Z}, \label{moda1}
\end{equation}
is a homotopy invariant of $D$.
\end{theorem}

\begin{proof}
Consider a homotopy $D_t$ of operators. The corresponding homotopy of
tangential operators will be denoted by $A_t$ and the homotopy of
operators on the base $X$ by $A_{0,t}$.

As $t$ varies, the index $\ind(D_t,\Pi_+(A_t))$ can jump. The geometric
Atiyah--Patodi--Singer index formula (see Appendix A) shows that the sum
$$
\limfunc{ind}\left( D_t,\Pi _{+}(A_t)\right) +\eta \left( A_t\right)
$$
is a smooth function of the parameter. Furthermore, the derivative of the sum
with respect to $t$ is a local invariant; namely, it is equal to the integral
over the manifold of an expression determined by the complete symbol of the
family $A_t$ of tangential operators. On the other hand, the complete symbols
of $A_t$ and $A_{0,t}$ locally coincide by Assumption l. Therefore, the jumps
in the sum
\begin{equation}
\limfunc{ind}\left( D_t,\Pi _{+}(A_t)\right) +\eta \left( A_t\right) -n\eta
\left( A_{0,t}\right) \label{bet3}
\end{equation}
of three terms come only from the third term. These jumps, however, are
multiples of $n$. Hence, reducing \eqref{bet3} modulo $n$, we obtain a quantity
independent of $t$. This proves the theorem.
\end{proof}

{\bf 2. A relation to known invariants.} We denote the homotopy invariant
constructed in the theorem by
$$
\inv D\stackrel{\operatorname{def}}=\func{mod}n\text{-}\limfunc{ind}\left(
D,\Pi _{+}(A)\right) +\eta \left( A\right) -n\eta \left( A_0\right) \in
\Bbb{R}/n\Bbb{Z}.
$$
In some special cases, it can be reduced to the following invariants.
\begin{enumerate}
\item {\em \emph{For a trivial covering, we have an isomorphism
$$
A\simeq \stackunder{n\;\;{\text{copies}}}{\underbrace{A_0\oplus A_0\oplus
...\oplus A_0}}.
$$
Hence $\eta(A)=n\eta(A_0)$, and $\inv D$ is thereby reduced to the
Freed--Melrose} $\func{mod}n$\emph{-index} \emph{\cite{FrMe1}}
\[
\func{mod}n\text{-}\limfunc{ind}(D,\Pi_+)\in \Bbb{Z}_n\subset \Bbb{R}/n\Bbb{Z}.
\]}

\item {\em \emph{\ On the other hand, if we consider only the fractional part
}$\left\{ \,\right\} $\emph{\ of (\ref{moda1}), then the index vanishes and we
obtain the Atiyah--Patodi--Singer} relative $\eta $-invariant
\emph{\cite{APS2,APS3}}
\[
\left\{ \eta \left( A_0\otimes 1_{\pi _{!}1}\right) -n\eta \left( A_0\right)
\right\} \in \Bbb{R}/\Bbb{Z}
\]
\emph{of the operator }$A_0$\emph{\ with coefficients in the flat bundle
}$\pi _{!}1\in \limfunc{Vect}\left( X\right) .$ }
\end{enumerate}

{\bf 3. The invariant $\inv D$ as an obstruction}. Let us give an
interpretation of the invariant $\inv$ as an obstruction.

Suppose that the manifold $M$ itself is the total space of a covering
$\widetilde{\pi }$ over some base $Y$ and the induced covering at the
boundary coincides with $\pi $:
\[
\begin{array}{ccc}
\partial M & \subset & M \\
\pi \downarrow \;\; & & \;\;\downarrow \widetilde{\pi } \\
X & \subset & Y.
\end{array}
\]

\begin{proposition}
\label{pro2} If an operator
\begin{equation}
\label{opa} D:C^\infty \left( M,E\right) \rightarrow C^\infty \left( M,F\right)
\end{equation}
satisfying Assumption {\rm 1} is the lift of some operator $D_0$ from
$Y$, then
\[
{\inv}D=0.
\]
\end{proposition}

\begin{proof}
By the Atiyah--Patodi--Singer formula, the sum
\[
\limfunc{ind}\left( D,\Pi _{+}(A)\right) +\eta \left( A\right)
\]
can be expressed as the integral over the manifold of an expression
depending on the complete symbol of $D$. On the other hand, the
operators $D$ and $D_0$ on the total space and on the base locally
coincide. Thus, we obtain
\[
\limfunc{ind}\left( D,\Pi _{+}(A)\right) +\eta \left( A\right) =n\left(
\limfunc{ind}\left( D_0,\Pi _+(A_0)\right) +\eta \left( A_0\right) \right) .
\]
Transposing the term $n\eta \left( A_0\right) $ to the left-hand side,
we obtain $\inv D=0$.
\end{proof}
\begin{remark}
{We point out that an arbitrary operator $D$ is not induced by some
elliptic operator $D_0$ on the base in general. Proposition~\ref{pro2}
gives a necessary condition for the existence of such operators: if $D_0$
exists, then $\inv D$ is zero.}
\end{remark}

The principal symbol of $D$ defines an element in the $K$-group of the
space $\overline{T^*M}^\pi$. To associate some topological index with
this element, in the following subsection we use Poincar\'e duality on
manifolds with singularities of the above-described type. For the
reader's convenience, the survey contains Appendices B and C, where we
discuss Poincar\'e duality on smooth manifolds and on singular
$\mathbb{Z}_n$-manifolds.

\subsection{The element of $K$-theory defined by a
manifold whose boundary is a covering} \label{s61} In this subsection, we
show that a pair $(M,\pi)$ consisting of a compact manifold $M$ and a
covering $\pi:\partial M\to X$ defines a special element in the $K$-group
$$
K_0(\mathcal{A}_{M,\pi},\mathbb{Q}/n\mathbb{Z})
$$
with coefficients in the group $\mathbb{Q}/n\mathbb{Z}$. Here
$\mathcal{A}_{M,\pi}$ is the $C^*$-algebra of the twisted
$\mathbb{Z}_n$-manifold. This algebra is defined in Appendix C. To define
our element, we first need a geometric description of elements of the
$K$-group. The coefficient group $\mathbb{Q}/n\mathbb{Z}$ looks slightly
complicated, and hence we first consider a geometric realization of the
$K$-groups with finite coefficient group $\mathbb{Z}_n$. After this, we
return to the complicated coefficients.

\textbf{1. $K$-theory with coefficients $\mathbb{Z}_n$.} There are several
equivalent approaches to the definition of these $K$-groups (e.g., see
\cite{Bla1}). We use the definition coming from topological $K$-theory. Namely,
the group $\mathbb{Z}_n$ has the corresponding \emph{Moore space}. It is the
2-dimensional complex $M_n$ obtained from the unit disk by the identification
of points on its boundary under the natural action of the group $\Bbb{Z}_n$:
\[
M_n=\left. \left\{ \left. z\in \Bbb{C}\;\right| \;\left| z\right| \leq
1\right\} \right/ \left\{ e^{i\varphi }\sim e^{i\left( \varphi +\frac{2\pi k}
n\right) }\right\}.
\]
In particular, for $n=2$ this space is the real projective plane
$M_2=\Bbb{RP}^2$.

The $K$-group with coefficients $\mathbb{Z}_n$ for an algebra
$\mathcal{A}$ can be defined as
\begin{equation}
K_0\left( \mathcal{A},\Bbb{Z}_n\right) =K_0\bigl( \widetilde{C} _0\left(
M_n,\mathcal{A}\right) \bigl) , \label{nesht1}
\end{equation}
where $\widetilde{C}_0\left( M_n,\mathcal{A}\right) $ is the algebra of
$\mathcal{A}$-valued continuous functions on the Moore space vanishing at
a marked point $pt\in M_n$. This definition is similar to the definition
of the topological $K$-groups with coefficients of a space $Y$:
$$
K^{0}\left( Y,\Bbb{Z}_n\right) =K^{0}\left( Y\times M_n,Y\times pt\right).
$$

\begin{proposition}
\label{zkoko} The group $K_0\bigl( \widetilde{C}_0\left( M_n,
\mathcal{A}_{M,\pi }\right) \bigr) $ is isomorphic to the group of stable
homotopy classes of triples $\left( E^{\prime },F^{\prime },\sigma ^{\prime
}\right),$ where $E^{\prime },F^{\prime }\in \limfunc{Vect}\left( M\times
M_n\right) $ and
\[
\sigma ^{\prime }:\pi _{!}\left. E^{\prime }\right| _{\partial
M}\longrightarrow \pi _{!}\left. F^{\prime }\right| _{\partial M}
\]
is an isomorphism over the product $X\times M_n.$ Here trivial triples
are understood as triples induced by a vector bundle isomorphism defined
on $M\times M_n$.
\end{proposition}
This result is a generalization of Lemma \ref{lemc2} to families. The
proof is similar to that of the lemma.

The triples introduced in this proposition can sometimes be written out
starting from explicit geometric data on the $\mathbb{Z}_n$-manifold
with the use of the following proposition.
\begin{proposition}
\label{zkk} A triple $\left( E,F,\sigma \right),$ where
\begin{equation}
E\in \limfunc{Vect}\left( M\right) ,\;F\in \limfunc{Vect}\left( X\right)
,\qquad \pi _{!}\left( \left. E\right| _{\partial M}\right) \stackrel{\sigma
}{\simeq }kF, \label{modnn1}
\end{equation}
and $\sigma $ is an isomorphism, defines an element of
$K_0\left( \mathcal{A}_{M,\pi },\Bbb{Z}_k\right) .$
\end{proposition}

\begin{proof}
On the Moore space $M_k$, we take a line bundle $\varepsilon $ such that
$\left[\varepsilon\right]-1\in\widetilde{K}\left(M_k\right)\simeq\Bbb{Z}_k$
is the generator of the reduced group. We also choose a trivialization
$\rho :k\varepsilon \rightarrow \Bbb{C}^k$.

Now consider the triple $\left( E,F,\sigma \right) $. With it we
associate the following element (in the sense of Proposition
\ref{zkoko}):
\[
\left[ E\otimes \varepsilon ,E,\sigma ^{\prime }\right] \in K_0\bigl(
\widetilde{C}_0\left( M_k,\mathcal{A}_{M,\pi }\right) \bigr) ,
\]
where the isomorphism $\sigma ^{\prime }$ is defined as the composition
\begin{equation}
\pi _{!}\left. E\right| _{\partial M}\otimes \varepsilon \stackrel{\sigma
\otimes 1}{\rightarrow }kF\otimes \varepsilon \simeq F\otimes k\varepsilon
\stackrel{1\otimes \rho }{\rightarrow }F\otimes \Bbb{C}^k\simeq kF\stackrel{
\sigma ^{-1}\otimes 1}{\rightarrow }\pi _{!}\left. E\right| _{\partial M}.
\label{longo1}
\end{equation}
\end{proof}

\textbf{2. $K$-theory with coefficients $\mathbb{Q}/n\mathbb{Z}$.} Now we
treat $\Bbb{Q}/n\Bbb{Z}$ as the direct limit of finite groups
$\Bbb{Z}_{nN}$ corresponding to the embeddings
\[
\Bbb{Z}_{nN} \subset \Bbb{Q}/n\Bbb{Z}, \quad x \mapsto x/N.
\]
Then the $K$-groups with coefficients $\Bbb{Q}/n\Bbb{Z}$ can also be
defined as the direct limit
\[
K_0\left( \mathcal{A}_{M,\pi },\Bbb{Q}/n\Bbb{Z}\right) =\stackunder{
\longrightarrow }{\lim }K_0\left( \mathcal{A}_{M,\pi },\Bbb{Z}_{nN}\right).
\]

{\bf 3. The $K$-theory element of a $\mathbb{Z}_n$-manifold.} The bundle
$\pi_!1\in\Vect(X)$ is flat, and therefore, for large $N$ there exists a
trivialization
\[
N\pi _{!}1\stackrel{\alpha }{\simeq }\Bbb{C}^{nN}.
\]
Then the triple $[\mathbb{C}^N,\mathbb{C},\alpha]$ defines an element of
$K_0(\mathcal{A}_{M,\pi},\mathbb{Z}_{Nn})$ by virtue of Proposition~\ref{zkk}.
By letting $N\to\infty$, we obtain an element of the $K$-group with
coefficients $\mathbb{Q}/n\mathbb{Z}$.

Unfortunately, this construction is ambiguous in the choice of the
trivialization $\alpha$, and different trivializations give different
elements.

It turns out that there is a canonical choice of $\alpha$. Namely, for
the covering $\pi:\partial M\to X$ consider the mapping
$$
f:X\longrightarrow BS_n
$$
to the classifying space $BS_n$ of the permutation group on $n$ elements. (Here
we treat the bundle $\pi$ as an associated bundle of the principal
$S_n$-bundle; then $f$ is the classifying mapping for this principal bundle.)
Moreover, $\pi_!1\in \Vect(X)$ is the pullback of the universal bundle
$\gamma_n\in \Vect_n(BS_n)$ over the classifying space. Assume that the range
of $f$ is contained in a finite skeleton $(BS_n)_{N'}$ of the classifying
space. The space $BS_n$ enjoys the Mittag-Leffler condition (see \cite{Ati6}):
for given $N'$, there exists an $L\geq 0$ such that
\begin{equation}
\func{Im}\!\left[ K^1\!\left( \left( BS_n\right) _{N+L}\right) \!\rightarrow\!
K^1\!\left( \left( BS_n\right) _N\right) \right] \!=\!\func{Im}\!\left[
K^1\!\left( \left( BS_n\right) _{N+M+L}\right) \!\rightarrow\! K^1\!\left(
\left( BS_n\right) _N\right) \right] \label{mile}
\end{equation}
for all nonnegative $M.$

Now let us choose an $N$ such that the sum of $N$ copies of the
restriction of the universal bundle $\gamma_n $ to $\left( BS_n\right)
_{N^{\prime }+L}$ is trivial. We take some trivialization
\[
N\gamma_n \stackrel{\alpha ^{\prime }}{\simeq }\Bbb{C}^{nN}.
\]
Then on $M$ we choose the induced trivialization
$$
\alpha=f^*\alpha'.
$$
One can show that for this special trivialization the $K$-theory element
defined by the triple $(\mathbb{C}^N,\mathbb{C},\alpha)$ is independent
of the choice of $\alpha'$ and $f$. We denote this element by
$$
[\widetilde{\pi_!1}]\in K_0\left( \mathcal{A}_{M,\pi },
\Bbb{Q}/n\Bbb{Z}\right).
$$

\subsection{The index defect formula}
\label{defzn}

\hspace{1mm}

{\bf 1. The general formula.}

\begin{theorem}
\label{defect2}{\em \cite{SaSt10}} For an elliptic operator $D$ on a
twisted $\mathbb{Z}_n$-manifold $(M,\pi)$, one has
$$
\inv D=\bigl\langle [\sigma(D)] ,\bigl[ \widetilde{\pi _{!}1}\bigr]
\bigr\rangle,
$$
where $\left\langle ,\right\rangle $ is Poincar\'e duality with
coefficients\textup:
\begin{equation}
\left\langle ,\right\rangle :K_c^0\bigl( \overline{T^{*}M}^\pi \bigr) \times
K_0\left( \mathcal{A}_{M,\pi },\Bbb{Q}/n\Bbb{Z}\right) \longrightarrow
\Bbb{Q}/n\Bbb{Z}. \label{proda2}
\end{equation}
\end{theorem}

\begin{proof}[Sketch of Proof]
It is well known in index theory how to compute fractional homotopy
invariants (see \cite{APS3}). The idea is to express the fractional
invariant as an index of some family of elliptic operators. Then it
suffices to apply the Atiyah--Singer formula for families. Let us use
this idea in the present situation.

1. Consider the family
\[
D^{*}\oplus \left( D\otimes 1_\varepsilon \right) :C^\infty \left( M,F\oplus
E\otimes \varepsilon \right) \longrightarrow C^\infty \left( M,E\oplus F\otimes
\varepsilon \right)
\]
of first-order elliptic operators on $M$ parametrized by the Moore space
$M_{nN}$. (The number $N$ will be chosen later on.) Here $D^{*}$ is the
adjoint operator, and $D\otimes 1_\varepsilon $ stands for operator $D$
with coefficients in the bundle $\varepsilon\in\Vect(M_{nN}) $. On the
other hand, consider the direct sum of $N$ copies of this family. It
turns out that the sum admits well-posed \emph{classical boundary
conditions}. Indeed, for sufficiently large $N$ the direct sum of $N$
copies of the flat bundle $\pi _{!}1$ is trivial. Let us choose some
trivialization
\begin{equation}
N\pi _{!}1\stackrel{\alpha }{\simeq }\Bbb{C}^{nN}. \label{triva1}
\end{equation}
Then we obtain an isomorphism
\[
\pi _{!}\left( N\left. E\right| _{\partial M}\right) \simeq \pi _{!}N\otimes
E_0\stackrel{\alpha \otimes 1}{\longrightarrow }\Bbb{C}^{nN}\otimes E_0
\]
on the base of the covering and a similar isomorphism
\[
\pi _{!}\left( N\left. E\right| _{\partial M}\right) \otimes \varepsilon
\simeq \pi _{!}N\otimes E_0\otimes \varepsilon \stackrel{\alpha \otimes 1}{
\longrightarrow }\Bbb{C}^{nN}\otimes E_0\otimes \varepsilon \simeq
nN\varepsilon \otimes E_0\stackrel{\rho \otimes 1}{\longrightarrow }\Bbb{C}
^{nN}\otimes E_0.
\]
The corresponding isomorphisms of the spaces of sections are denoted by
\[
B_1=\alpha \otimes 1:C^\infty \left( X,\pi _{!} \left(N\left. E\right|
_{\partial M}\right)\right) \stackrel{}{\longrightarrow }C^\infty \left(
X,\Bbb{C}^{nN}\otimes E_0\right) ,
\]
\[
B_2=(\rho \otimes 1)\left( \alpha \otimes 1\right) :C^\infty \left( X,\pi _{!}
\left(N\left. E\right| _{\partial M}\right)\otimes \varepsilon\right) \stackrel{}{
\longrightarrow }C^\infty \left( X,\Bbb{C}^{nN}\otimes E_0\right) .
\]
In this notation, we define the following family of boundary value
problems:
\begin{equation}
\left\{
\begin{array}{l}
\begin{array}{cc}
ND^{*}u=f_1,\quad & N\left( D\otimes 1_\varepsilon \right) v=f_2,
\end{array}
\vspace{1mm} \\
B_1\beta _E\left. u\right| _{\partial M}+B_2\beta _E\left. v\right| _{\partial
M}=g,\qquad g\in C^\infty \left( X,\Bbb{C}^{nN}\otimes E_0\right) .
\end{array}
\right. \label{nelin1}
\end{equation}
One can show that this family consists of elliptic boundary value
problems. By $\Phi_\alpha(D)$ we denote the family (\ref{nelin1}) for
some given trivialization $\alpha$.

2. By virtue of the embedding $\mathbb{Z}_{nN}\subset
\mathbb{Q}/n\mathbb{Z}$, we can treat the index-residue
$\ind\Phi_\alpha(D)\in \widetilde{K}(M_n)\simeq \mathbb{Z}_{nN}$ of the
family as a fractional rational number. The main step in the proof of the
defect formula is the following result.
\begin{lemma}
\label{lem3}One has
$$
\inv D=\ind \Phi_\alpha (D)\in \mathbb{R}/n\mathbb{Z}
$$
\emph{(}provided the trivialization $\alpha$ is chosen as in Subsection
\emph{\ref{s61}).}
\end{lemma}

\begin{proof}[Proof of the Lemma] There is a linear ellipticity-preserving
homotopy  between (\ref{nelin1}) and the boundary value problem
\[
\left\{
\begin{array}{l}
ND^{*}u=f_1,\qquad N\left( D\otimes 1_\varepsilon \right) v=f_2,\vspace{1mm}
\\
B_1\beta _E\Pi _{-}(A)\left. u\right| _{\partial M}+B_2\beta _E\Pi
_{+}(A)\left. v\right| _{\partial M}=g,\quad g\in C^\infty \left(
X,\Bbb{C}^{nN}\otimes E_0\right)
\end{array}
\right.
\]
for the same operator $ND^{*}\oplus N\left( D\otimes 1_\varepsilon
\right)$. By the Agranovich--Dynin formula \cite{AgDy1}, the index of the
family $\Phi_\alpha ( D) $ is the sum
\begin{multline}
\ind\Phi_\alpha ( D)=N \ind(D,\Pi_+(A))([\varepsilon]-1)  \\ +\ind\bigl(
N\func{Im}\Pi _{-}\left( A\right) \oplus N\func{Im}\Pi _{+}\left( A\right)
\otimes \varepsilon \stackrel{B_1+B_2}{\longrightarrow }C^\infty \left( X,
\Bbb{C}^{nN}\otimes E_0\right) \bigr) \label{opcl1}
\end{multline}
of the index of the family of spectral boundary value problems for
$ND^{*}$ and $N\left( D\otimes 1_\varepsilon \right)$ and the index of an
operator family on the boundary. Let us compute the index of that family.

To this end, we decompose the space of right-hand sides of the family into the
direct sum
\[
C^\infty \left( X,\Bbb{C}^{nN}\otimes E_0\right) \simeq nN\func{Im}\Pi
_{-}\left( A_0\right) \oplus nN\varepsilon \otimes \func{Im}\Pi _{+}\left(
A_0\right).
\]
Here the isomorphism is defined by the formula
\[
nN\func{Im}\Pi _{-}\left( A_0\right) \oplus nN\varepsilon \otimes \func{Im}
\Pi _{+}\left( A_0\right) \stackrel{1+\left( \rho \otimes 1\right) }{
\longrightarrow }C^\infty \left( X,\Bbb{C}^{nN}\otimes E_0\right) .
\]
This permits us to express the index of the family $B_1+B_2$ in
Eq.~(\ref{opcl1}) as
\[
\limfunc{ind}\bigl( N\func{Im}\Pi _{+}\left( A\right) \stackrel{\Pi _{+}\left(
A_0\right) \beta _E}{\longrightarrow }nN\func{Im}\Pi _{+}\left( A_0\right)
\bigr) \left( \left[ \varepsilon \right] -1\right) \in \widetilde{K}\left(
M_{nN}\right) .
\]
Finally, pushing down the space $\func{Im}\Pi _{+}\left( A\right) $ to
the base of the covering, we reduce the index to the form
\[
\limfunc{ind}\bigl( N\func{Im}\Pi _{+}\left( \pi _{!}A\right) \stackrel{\Pi
_{+}\left( A_0\right) }{\longrightarrow }nN\func{Im}\Pi _{+}\left( A_0\right)
\bigr) \left( \left[ \varepsilon \right] -1\right) .
\]
The index of an elliptic operator (not a family!) in this formula can be
expressed in terms of $\eta$-invariants by the Atiyah--Patodi--Singer ``index
formula for flat bundles" \cite{APS3}
\begin{multline}
\limfunc{ind}\bigl( N\func{Im}\Pi _{+}\left( \pi _{!}A\right) \stackrel{\Pi
_{+}\left( A_0\right) }{\longrightarrow }nN\func{Im}\Pi _{+}\left( A_0\right)
\bigr)  \\
=N\eta \left( A\right) -nN\eta \left( A_0\right) +\left\langle \left[ \sigma
\left( A_0\right) \right] ,\left[ \pi _{!}1\right] \right\rangle , \label{apps}
\end{multline}
where the brackets stand for the pairing
\begin{equation}
\left\langle ,\right\rangle :K^1_c\left( T^{*}X\right) \times K^1\left( X,\Bbb{
Q}\right) \longrightarrow \Bbb{Q} \label{pai2}
\end{equation}
of the difference element of the elliptic self-adjoint operator $A_0$
\[
\left[ \sigma \left( A_0\right) \right] \in K^1\left( T^{*}X\right)
\]
and the element $\left[ \pi _{!}1\right] \in K^1\left( X,\Bbb{Q}\right) $
corresponding to the trivialized flat bundle $N\pi _{!}1$. (More details on
this formula can be found in \cite{Gil1}.)

Substituting \eqref{apps} into Eq. \eqref{opcl1} and transposing the
$\eta$-invariants to the left-hand side, we obtain
$$
\inv D =\ind \Phi_\alpha(D)+\left\langle \left[ \sigma \left( A_0\right)
\right] ,\left[ \pi _{!}1\right] \right\rangle.
$$

2) It remains to show that for a special choice of the
trivialization~(\ref{triva1}) the last term in (\ref{apps}) is zero.

Indeed, consider the classifying mapping
$$
f:X\rightarrow \left( BS\right) _{N^{\prime }}.
$$
The computation of the pairing (\ref{pai2}) can be moved to the classifying
space:
\begin{equation}
\left\langle \left[ \sigma \left( A_0\right) \right] ,\left[ \pi _{!}1\right]
\right\rangle =\left\langle f_{!}\left[ \sigma \left( A_0\right) \right]
,\left[ \gamma \right] \right\rangle ,\quad \left[ \pi _{!}1\right]
=f^{*}\left[ \gamma \right] \in K^1\left( X,\Bbb{Q}\right), \label{al1}
\end{equation}
where $[\gamma_n]\in K^1((BS_n)_{N'})\otimes\Bbb{Q}$ is the element
defined by the trivialized flat bundle $N\gamma_n$. The embedding $\left(
BS_n\right) _{N^{\prime }}\subset \left( BS_n\right) _{N^{\prime
}+L^{\prime }}$ induces the commutative diagram
$$
\begin{array}{ccc}
K^1_c\!\left( T^{*}\!\left( BS_n\right) _{N^{\prime }}\right)\! \times\!
K^1\!\left( \left( BS_n\right) _{N^{\prime }+L^{\prime }}\!,\mathbb{Q}\right)
 &\!\!\rightarrow\!\!&\!\!K^1_c\!\left( T^{*}\!\left( BS_n\right) _{N^{\prime
 }}\right)\!\times\! K^1\!\left(
\left( BS_n\right) _{N^{\prime }}\!,\mathbb{Q}\right)
 \\
\downarrow& &\downarrow\\
K^1_c\!\left( T^{*}\!\left(\! BS_n\right) _{N^{\prime }+L^{\prime }}\right)
\!\times\! K^1\!\!\left( \left( BS_n\right) _{N^{\prime }+L^{\prime
}}\!,\mathbb{Q}\right)
 \!\!\!&\!\!\rightarrow\!\!& \mathbb{Q}.
\end{array}
$$
It follows from this diagram that the pairing (\ref{al1}) gives zero,
since the range of the mapping
$$
K^1_c\left( T^{*}\left( BS_n\right) _{N^{\prime }}\right)\to K^1_c\left(
T^{*}\left( BS_n\right) _{N^{\prime }+L^{\prime }}\right)
$$
is contained in the torsion subgroup.

Therefore, the expression for $\limfunc{ind}\Phi_\alpha (D) $ is reduced
to the desired relation
\[
\limfunc{ind}\Phi _\alpha(D) =\inv D .
\]
This proves Lemma \ref{lem3}.
\end{proof}

3. To complete the proof of the theorem, it suffices to relate the
families index of $\Phi_\alpha(D)$ to the Poincar\'e duality pairing. In
the index theory of classical boundary value problems, there is a
well-known operation of ``order reduction," which reduces elliptic
boundary value problems to zero-order operators (see \cite{Hor3} or
\cite{SaScS4}). This operation preserves the index. Applying this
operation to the family $\Phi_\alpha(D)$, one can show that the result
is a family of admissible operators in the sense of Subsection \ref{noni}
of Appendix C. A computation shows that this family of admissible
operators coincides with the family corresponding to the product of the
elements
$$
[\sigma(D)]\in K_c(\overline{T^*M}^\pi),\quad [\widetilde{\pi_!1}]\in
K_0(\mathcal{A}_{M,\pi},\mathbb{Q}/n\mathbb{Z}).
$$
Since the Poincar\'e duality pairing of two elements is defined as the
index, we obtain the desired formula for the index of the problem as the
Poincar\'e pairing with coefficients:
$$
\ind \Phi_\alpha(D)=\bigl\langle [\sigma(D)] ,\bigl[ \widetilde{\pi
_{!}1}\bigr] \bigr\rangle.
$$
Together with the equality in Lemma \ref{lem3}, this completes the proof
of the theorem.
\end{proof}

\begin{remark}
{ Using the results of Appendix C, one can compute Poincar\'e duality
topologically for a regular covering $\pi$. In this case, we obtain an index
defect formula in topological terms.}
\end{remark}

{\bf 2. The index defect in the $G$-equivariant case.} In a number of
cases, the invariant $\inv$ can be computed effectively with the use of
Lefschetz theory. Suppose we are given an action of a finite group $G$ on
$M$ such that the action is free on the boundary $\partial M$. As a
covering $\pi$, we take the natural projection to the quotient space:
$$
\pi:\partial M\to \partial M/G.
$$
Note that we do not require that the action be free in the interior of $M$.

Consider a $G$-equivariant elliptic operator $D$ on $M$. By $L\left( D,g\right)
\in \Bbb{C}$, $g\in G$, we denote the usual contribution to the Lefschetz
formula (see \cite{Don3}) of the fixed point set of  an element $g$.
\begin{proposition}
One has
\begin{equation}
\inv D\equiv -\sum_{g\neq e}L\left( D,g\right) \text{ }
\left( \func{mod}n\right) . \label{lef1}
\end{equation}
\end{proposition}

\begin{proof}
Consider the equivariant index $\limfunc{ind}_g\left( D,\Pi _{+}\right) $ of
the spectral boundary value problem and the equivariant $\eta $-function\ (see
\cite{Don3}) of the tangential operator $A$ on the boundary.

By $\left( D,\Pi _{+}\right) ^G$ and $A^G$ we denote the restrictions of the
corresponding operators to the spaces of $G$-invariant sections. Clearly, $A^G$
is equivalent to $A_0$ on the base of the covering. On the other hand, one can
express the usual invariants in terms of their equivariant counterparts:
\[
\limfunc{ind}\left( D,\Pi _{+}\right) ^G=\frac 1{\left| G\right| }\sum_{g\in
G}\limfunc{ind}_g\left( D,\Pi _{+}\right), \qquad \eta \left( A^G\right) =\frac
1{\left| G\right| }\sum_{g\in G}\eta \left( A,g\right).
\]
This is easy to check with the use of character theory. By virtue of these
expressions, we can rewrite $\inv D$ as
\[
\inv D=\limfunc{ind}_e\left( D,\Pi _{+}\right) -\sum_{g\neq e}\eta \left(
A,g\right) .
\]
The $\eta $-invariants here can be expressed by the equivariant
Atiyah--Patodi--Singer formula (see \cite{Don3}):
\[
-\eta \left( A,g\right) =\limfunc{ind}_g\left( D,\Pi _{+}\right) -L\left(
D,g\right) .
\]
Thus, we obtain
\[
\inv D=\left| G\right| \limfunc{ind}\left( D,\Pi _{+}\right) ^G-\sum_{g\neq
e}L\left( D,g\right) .
\]
This gives the desired equation~(\ref{lef1}).
\end{proof}

\subsection{An application to $\eta$-invariants\label{exa4}}

The index defect formula enables one to express the fractional part
of the $\eta $-invariant in the following situation.

Let $M$ be an even-dimensional spin manifold such that its boundary is a
covering with spin structure coinciding with that induced from the base. Let us
also choose a vector bundle $E\in \limfunc{Vect}\left( M\right) $ that at the
boundary is the lift of some bundle $E_0\in\Vect(X)$.  In a collar neighborhood
of the boundary, we take a product metric on $M$ that is the lift of some
metric from the base. Finally, we choose a similar connection in $E$.

\begin{proposition}
The Dirac operator $D_M$ on $M$ with coefficients in $E$ satisfies the
assumptions of Theorem \emph{\ref{defect2}}, and one has the formula
\[
\left\{ \eta \left( D_X\otimes 1_{E_0}\right) \right\} =\frac 1n\biggl(
\int\limits_M \widehat{A}\left( M\right) \limfunc{ch}E-\left\langle \left[
\sigma \left( D_M\right) \right] ,\left[ \widetilde{\pi _{!}1}\right]
\right\rangle \biggr) \in \Bbb{R}/\Bbb{Z},
\]
for the fractional part of the $\eta$-invariant of the self-adjoint Dirac
operator $D_X$ with coefficients in the bundle $E_0$ on the base of the
covering.
\end{proposition}

\begin{proof}
This formula follows from the defect formula, where the index of the
spectral boundary value problem is expanded by the Atiyah--Patodi--Singer
formula
\[
\limfunc{ind}\left( D_M,\Pi _{+}\right) =\int\limits_M\widehat{A}\left(
M\right) \limfunc{ch}E-\eta \left( D_{\partial M}\right) .
\]
\end{proof}

\section{Appendix A. The Atiyah--Patodi--Singer $\eta$-invariant}
In this appendix, we give a brief overview of the spectral
$\eta$-invariant. Most of the results were proved in the original paper
\cite{APS3}, and there is also a very stimulating exposition in
\cite{BBW1}. Therefore, here we either omit the proofs or only indicate
the main idea.
\subsection{The geometric index formula and the $\eta$-invariant}

Atiyah--Patodi--Singer \cite{APS1} gave a formula for the index of
spectral boundary value problems for geometric first-order operators.
Namely, using the heat equation method \cite{ABP1}, they obtained the
formula
$$
\ind (D,\Pi_+(A))=\int_X a(D)-\eta(A)
$$
for the index of the spectral boundary value problem on a manifold $X$
for an operator having the form \eqref{deco2} near the boundary. The
first contribution is defined by the constant term $a(D)$ in the local
asymptotic expansion of the heat kernel
$$
{\rm tr}(e^{-tD^*D}(x,x))-{\rm tr}(e^{-tDD^*}(x,x))
$$
as $t\to 0$. This term is determined just as in the case of operators on closed
manifolds as some algebraic expression in the coefficients of the operator. The
new feature of the spectral boundary value problem is the so-called
$\eta$\emph{-invariant} of the tangential operator $A$. Let us recall its
definition.

Let $A$ be an elliptic self-adjoint operator of some positive order on a closed
manifold $M.$ The $\eta $\emph{-function} of $A$ is defined by the formula
\[
\eta \left( A,s\right) =\sum_{\lambda _j\in \limfunc{Spec}A,\lambda _j\neq 0}
\limfunc{sgn}\lambda _i\left| \lambda _i\right| ^{-s}\equiv \limfunc{Tr}
\left( A\left( A^2\right) ^{-s/2-1/2}\right) .
\]
It is analytic in the half-space $\func{Re}s>\dim M/\limfunc{ord}D$
(where the series is absolutely convergent). { This spectral function is
a generalization of the $\zeta $-function
\[
\zeta \left( A,s\right) =\sum_{\lambda _j\in \limfunc{Spec}A}\left| \lambda
_i\right| ^{-s}
\]
of positive definite elliptic operators. } By analogy with the $\zeta
$-invariant
\[
\zeta \left( A\right) =\frac 12\zeta \left( A,0\right),
\]
it is natural to introduce the following definition.
\begin{definition}The $\eta$\emph{-invariant} of $A$ is the number
\begin{equation}
\eta \left( A\right) =\frac 12\left( \eta \left( A,0\right) +\dim \ker A\right)
\in \Bbb{R}. \label{omo}
\end{equation}
\end{definition}
Of course, for this definition to make sense, it is necessary to have an
analytic extension of the $\eta $-function to the point $s=0 $. Analytic
methods show that the $\eta$-function extends meromorphically to the
entire complex plane, possibly with a simple pole at the origin.
Atiyah--Patodi--Singer \cite{APS3} for odd-dimensional manifolds and
Gilkey \cite{Gil4} for even-dimensional manifolds proved, using global
topological methods, that the residue at the point $s=0$ is nevertheless
zero. Thus, the $\eta $-function is holomorphic at $s=0$ and the
$\eta$-invariant is well defined.

The $\eta$-invariant is by definition only a spectral invariant, and it
can vary for a deformation of the operator. Consider an example.

\begin{example}
\label{exo}On the circle of length $2\pi $ with coordinate $\varphi $,
consider the operator
\[
A_t=-i\frac d{d\varphi }+t.
\]
Here $t$ is some real constant. Let us compute the $\eta $-invariant.
Since the spectrum is given by the lattice $t+\Bbb{Z}$ (with period
one), it follows that the $\eta$-invariant is a periodic function of the
parameter $t$. Thus, we can suppose that $0<t<1$. Collecting the
eigenvalues in pairs, we can rewrite the $\eta $-function as
\[
\eta \left( A_t,s\right) =\sum_{n\geq 1}\left[ \left( n+t\right)
^{-s}-\left( n-t\right) ^{-s}\right] +t^{-s}.
\]
Let us show that this series absolutely converges on the real line for
$s>0$ and compute the limit as $s\rightarrow +0.$ By the Taylor formula,
we have
\[
\left[ \left( n+t\right) ^{-s}-\left( n-t\right) ^{-s}\right]
=-2tsn^{-s}+O\left( \frac s{n^{2+s}}\right) .
\]
Thus, as $s\rightarrow +0$ we obtain
\[
\sum_{n>1}\left[ \left( n+t\right) ^{-s}-\left( n-t\right) ^{-s}\right]
\sim -2ts\sum_{n\geq 1}n^{-s}\sim -2ts\int_1^\infty x^{-s-1}dx=-2t,
\]
and for the $\eta $-invariant we have
\[
\eta \left( A_t\right) =\frac{\eta \left( A_t,0\right) +\dim \ker A_t}
2=\frac 12-\left\{ t\right\} ,
\]
where $\left\{ \,\right\} \in \left[ 0,1\right) $ is the fractional part of a
number. Hence, for a smooth family $A_t$ of elliptic operators the
corresponding family of $\eta$-invariants is only piecewise smooth. In
addition, the jumps (which are integral!) occur for parameter values such that
some eigenvalue changes its sign. Let us also mention the half-integer
parameter  values $t\in \Bbb{Z}+1/2,$ where the spectrum of $A_t$ is symmetric
with respect to the origin. For these values, the $\eta $-function vanishes
identically. In this case, one says that the operator $A_t$ has a spectral
symmetry, and the $\eta$-invariant is regarded as a {\em measure of spectral
asymmetry} of the operator.
\end{example}

\subsection{The derivative of the $\eta$-invariant}

It turns out that the $\eta$-invariant is a piecewise smooth function of the
parameter in the general case as well. More precisely, the following result is
valid.

\begin{proposition}
\emph{\cite{APS3}}\label{pr16} For a smooth family $\{A_t\}_{t\in
\left[0,1\right]} $ of elliptic self-adjoint operators, the following
assertions hold.
\end{proposition}

\begin{enumerate}
\item \emph{If the operator }$A_{t_0}$\emph{\ is invertible, then the function
}$\eta \left( A_t\right) $\emph{\ is smooth in a neighborhood of }$t_0 $\emph{
and its derivative can be expressed in terms of the derivative of the }$\zeta
$\emph{-function of the auxiliary family}
$$ B_{t,t_0}=\left| A_{t_0}\right| +\left(
t-t_0\right) \left.\left(\frac d {d\tau}A_\tau\right)\right|_{\tau=t_0}.
$$
\emph{Namely,}
\begin{equation}
\left. \frac d{dt}\eta \left( A_t\right) \right| _{t=t_0}=\left. \frac
d{dt}\zeta \left( B_{t,t_0}\right) \right| _{t=t_0}. \label{deri}
\end{equation}

\item \emph{In the general case, }$\eta \left( A_t\right) $ \emph{is
piecewise smooth and can be decomposed in the form}
\begin{equation}
\begin{array}{c}
\eta \left( A_{t^{\prime }}\right) -\eta \left( A_0\right)
=\int\limits_0^{t^{\prime }}\omega \left( t_0\right)
dt_0+\limfunc{sf}\left(
A_t\right) _{t\in \left[ 0,t^{\prime }\right] }, \\
\\
\omega \left( t_0\right) =\left. \frac d{dt}\zeta \left(
B_{t,t_0}\right) \right| _{t=t_0}\in C^\infty \left[ 0,1\right] ,
\end{array}
\label{deka3}
\end{equation}
\emph{as the sum of a smooth function of the parameter and a piecewise constant
function $\limfunc{sf}$ called the spectral flow} (\emph{see Subsection
}\ref{sflow1})\emph{. Here}
\[
B_{t,t_0}=\left| A_{t_0}\right| +P_{\ker A_{t_0}}+\left( t-t_0\right)
\left.\left(\frac d {d\tau}A_\tau\right)\right|_{\tau=t_0},
\]
\emph{and }$P_{\ker A_{t_0}}$\emph{\ is the orthogonal projection on the
kernel of }$A_{t_0}.$
\end{enumerate}

\begin{corollary}
The fractional part
\[
\left\{ \eta \left( A_t\right) \right\} \in \Bbb{R}/\Bbb{Z\simeq S}^1
\]
of the $\eta$-invariant is a smooth function of the parameter $t$ for a smooth
family $A_t.$
\end{corollary}

In the general case, the $\eta$-invariant is not homotopy invariant and can
take arbitrary real values.

\subsection{The homotopy invariance of the $\eta$-invariant}

It turns out, however, that in some special operator classes the
$\eta$-invariant possesses homotopy invariance. To this end, it is
necessary that the two components in the decomposition (\ref{deka3})
vanish. The easiest way to eliminate the second component, i.e. the
spectral flow, is to consider only the fractional part $ \left\{ \eta
\left( A\right) \right\}$ of the $\eta$-invariant. To obtain the
vanishing of the second component, it is convenient to use the formula
for the derivative on the left-hand side in Eq. (\ref{deri}). Seeley
\cite{See5} proved (see also \cite{Agr3} and \cite{GrSe2}) that the value
of the $\zeta$-function at the origin can be expressed via the principal
symbol of the operator. Explicitly, for a positive self-adjoint operator
$A$ with principal symbol having the asymptotic expansion
\[
\sigma \left( A\right) \sim a_m+a_{m-1}+a_{m-2}+...,
\]
the $\zeta $-invariant is computed by the following procedure. Let us define
the symbols $b_{-m-j},$ $j\geq 0,$ by the recursion relations
\begin{multline}
b_{-m-j}\left( x,\xi ,\lambda \right) \left( a_m\left( x,\xi \right)
-\lambda \right) \\
+\sum\limits_{\substack{
k+l+\left| \alpha \right| =j, \\
l>0} }\frac 1{\alpha !}\left( -i\partial _\xi \right) ^\alpha b_{-m-k}\left(
x,\xi ,\lambda \right) \left( -i\partial _x\right) ^\alpha a_{m-l}\left( x,\xi
\right) =0. \label{reku}
\end{multline}
The symbols in (\ref{reku}) depend on the coordinates $x$, momenta
$\xi$, and additionally on the parameter $\lambda$. In this notation, the
$\zeta$-invariant is given by
\begin{equation}
2\zeta \left( A\right) =
\frac
1{\left( 2\pi \right) ^{\dim
M}\limfunc{ord}A}\int\limits_{S^{*}M}dxd\xi \int\limits_0^\infty
b_{-\dim M-\limfunc{ord}A}\left( x,\xi ,-\lambda \right) d\lambda .
\label{apofeoz}
\end{equation}
Note the following properties of this formula.

\begin{enumerate}
\item (Locality.) For two locallyisomorphic operators $A$ and $A^{\prime }$,
their $\zeta $-in\-var\-iants coincide:
\[
\zeta \left( A\right) =\zeta \left( A^{\prime }\right).
\]
(Operators are said to be \emph{locally isomorphic} if their complete
symbols coincide in a neighborhood of every point of the manifold in some
coordinate system for some trivialization of vector bundles.)

\item (Homogeneity.) The terms $b_j$ are positively homogeneous
functions:
\[
b_j\left( x,t\xi ,t^m\lambda \right) =\left( -1\right) ^jb_j\left(
x,\xi ,\lambda \right) ,\qquad t>0.
\]
\end{enumerate}

These properties enable one to find classes of operators where the
derivative of the $\eta$-invariant in Eq. (\ref{deri}) is zero. Let us
introduce some of the known classes.

First, we define a class of pseudodifferential operators that generalize
differential operators.

\begin{definition}{\em
\emph{A classical pseudodifferential operator }$A$\emph{\ with complete
symbol }
\[
\sigma \left( A\right) \sim a_m+a_{m-1}+a_{m-2}+...
\]
\emph{is said to be }$\Bbb{R}_{*}$-invariant\emph{ if the components of
its complete symbol are homogeneous functions}
\[
a_j\left( x,t\xi \right) =t^ja_j\left( x,\xi \right) ,\qquad t\in \Bbb{R}_{*},
\]
\emph{with respect to the group of nonzero real numbers }$\Bbb{R}_{*}.$ }
\end{definition}

To define the second class of operators, recall that a \emph{flat bundle}
$\gamma \in \limfunc{Vect}\left( M\right) $ is a vector bundle with locally
constant (i.e. constant on connected subsets) transition functions. For an
operator
\[
A:C^\infty \left( M,E\right) \longrightarrow C^\infty \left( M,F\right),
\]
one can define an \emph{operator with coefficients in the flat bundle}, denoted
by
\[
A\otimes 1_\gamma :C^\infty \left( M,E\otimes \gamma \right) \longrightarrow
C^\infty \left( M,F\otimes \gamma \right).
\]
This is locally isomorphic to the direct sum of $\dim \gamma $ copies of $A.$
One can globally define this operator by gluing the local complete symbols of
$A$ with the use of a partition of unity. We  also require that the transition
functions of a flat bundle be unitary.

\begin{example}
On the circle $\mathbb{S}^1$ with coordinate $\varphi$, consider the vector
bundle $\gamma $ with the transition function $e^{2\pi it}.$ Then the operator
\[
-i\frac d{d\varphi }+t
\]
from Example \ref{exo} is locally isomorphic to $-i\frac d{d\varphi
}\otimes 1_\gamma .$ The isomorphism
\[
e^{-ti\varphi }\left( -i\frac d{d\varphi }\right) e^{it\varphi
}=-i\frac d{d\varphi }+t
\]
is given by the trivialization $e^{it\varphi }$ of $\gamma $.
\end{example}

\begin{theorem}
\label{droba} {\em \cite{APS3}, \cite{Gil2}} The fractional part of the
$\eta$-invariant is homotopy invariant in the following two classes of elliptic
self-adjoint operators:
\end{theorem}

\begin{enumerate}
\item \emph{the class of direct sums}
\[
A\otimes 1_\gamma \oplus \left( -\dim \gamma A\right)
\]
\emph{with a given flat bundle }$\gamma\in\Vect(X) ${;}

\item \emph{the class of }$\Bbb{R}_{*}$\emph{-invariant operators if the
following parity condition is satisfied }:
\begin{equation}
\dim A+\limfunc{ord}M\equiv 1\left( \func{mod}2\right) . \label{para}
\end{equation}
\end{enumerate}

\begin{proof}
$1)$ The main idea of the proof is to use the locality of the $\zeta
$-invariant. More precisely, for the fractional part of $\eta$ one has
\[
\left\{ \eta \left( A\otimes 1_\gamma \oplus \left( -nA\right) \right)
\right\} =\left\{ \eta \left( A\otimes 1_\gamma \right) \right\}
-\left\{ n\eta \left( A\right) \right\} ,\qquad n=\dim \gamma .
\]
Then for a smooth homotopy $A_t$ this gives (see Eq. (\ref{deri}))
\[
\frac d{dt}\left\{ \eta \left( A_t\right) \right\} =\frac d{dt}\zeta \left(
B_t\right) ,\quad B_t=\left| A_{t_0}\right| +\left(
t-t_0\right)\left.\left(\frac d {d\tau}A_\tau\right)\right|_{\tau=t_0}.
\]
A similar formula
\[
\frac d{dt}\left\{ \eta \left( A_t\otimes 1_\gamma \right) \right\}
=\frac d{dt}\zeta \left( B_t^{\prime }\right)
\]
is valid for the derivative of the operators with coefficients in the
flat bundle. Recall that $A\otimes 1_\gamma $ and $nA$ are locally
isomorphic. Therefore, the positive definite operators $B_t$ and
$B_t^{\prime }$ are locally isomorphic as well. Therefore, the locality
of the $\zeta $-invariant gives the desired relation
\[
\frac d{dt}\left\{ \eta \left( A_t\otimes 1_\gamma \right) \right\}
=\frac d{dt}\left\{ n\eta \left( A_t\right) \right\} .
\]
This proves the homotopy invariance.

$2)$ Consider first the case of even-order operators. It is easy to see
that for a homotopy $A_t$ the corresponding positive definite operators
$B_t$ (see Eq.~(\ref{deri})) are also $\Bbb{R}_{*}$-invariant. By
induction, this gives the $\Bbb{R}_{*}$-homogeneity
\begin{equation}
b_j\left( x,-\xi ,\lambda ,t\right) =\left( -1\right) ^jb_j\left( x,\xi
,\lambda ,t\right) \label{odno}
\end{equation}
of the coefficients corresponding to $B_t$ (see the recursion
relations (\ref{reku})). Hence,
\[
\frac d{dt}\left\{ \eta \left( A_t\right) \right\} =\frac d{dt}\zeta \left(
B_t\right) =\limfunc{Const}\frac d{dt}\biggl( \int\limits_{S^{*}M}dxd\xi
^{\prime }\int\limits_0^\infty b_{-\dim M-\limfunc{ord}A}\left( x,\xi ,-\lambda
,t\right)d\lambda \biggr)  .
\]
Using the homogeneity (\ref{odno}) and the assumption that the manifold is
odd-di\-men\-sion\-al, we see that the integrand $b_{-\dim
M-\limfunc{ord}A}\left( x,\xi ,-\lambda ,t\right) $ is an odd function on the
sphere $S_x^{*}M.$ Therefore, the integral is zero, and we obtain the desired
relation
\[
\frac d{dt}\left\{ \eta \left( A_t\right) \right\} =0.
\]

For odd-order operators, one obtains the different homogeneity
\[
b_j\left( x,-\xi ,\lambda ,-t\right) =\left( -1\right) ^{j+1}b_j\left(
x,\xi ,\lambda ,t\right) .
\]
Substituting these homogeneous functions into the expression for the
$\zeta $-invariant, we obtain
\[
\frac d{dt}\zeta \left( B_t\right) =\frac d{dt}\zeta \left( B_{-t}\right) .
\]
This gives the desired homotopy invariance of the fractional part of the
$\eta$-invariant:
\[
\frac d{dt}\left\{ \eta \left( A_t\right) \right\} =\frac d{dt}\zeta
\left( B_t\right) =0.
\]
\end{proof}

\begin{remark}
Note that in the proof we also obtained the vanishing of both the
derivative and the $\zeta $-invariant itself for
$\Bbb{R}_{*}$-invariant even-order operators on odd-dimensional
manifolds.
\end{remark}

\section{Appendix B. Elliptic operators and Poincar\'e duality.
Smooth theory}

In this appendix, we show that Poincar\'e duality in $K$-theory on
smooth manifolds can naturally be described in terms of elliptic
operators. We consider both closed manifolds and manifolds with boundary.
Using the Poincar\'e isomorphism, we construct Poincar\'e duality as a
nonsingular pairing. In this context, the Atiyah--Singer index theorem
can be used to make the pairing effectively computable. The topics
covered in this appendix can also be found in the recent book
\cite{HiRo1}. Our approach is closer to differential equations.

\subsection{The Poincar\'e isomorphism on a closed manifold}\label{ss51}

\hspace{1mm}

{\bf 1. Atiyah's generalized elliptic operators.} It is well known that for a
sufficiently nice topological space $X$ (e.g., a finite $CW$-complex) there is
a pairing
\begin{equation}
\label{kog1} H_i(X,\mathbb{Z})\times H^i(X,\mathbb{Z})\longrightarrow
\mathbb{Z}
\end{equation}
of homology and cohomology groups. (The pairing is nondegenerate on the free
parts of the groups.) A similar pairing can be constructed in $K$-theory on a
smooth closed manifold $M$. Namely, we consider the pairing
\begin{equation}
\label{kog2}
\langle,\rangle:K^0_c(T^*M)\times K^0(M)\longrightarrow \mathbb{Z}
\end{equation}
taking the difference element $[\sigma(D)]\in K^0_c(T^*M)$ of an elliptic
operator
$$
D:C^\infty(M,E)\longrightarrow C^\infty(M,F)
$$
and a vector bundle $G\in\Vect(M)$ to the index of $D$ with coefficients
in $G$.\footnote{ Recall that an operator with coefficients in a vector
bundle has the principal symbol
$$
\sigma(D)\otimes 1_G:\pi^*(E\otimes G)\longrightarrow \pi^*(F\otimes G).
$$
Any operator with this symbol is denoted by $D\otimes 1_G$.} The pairing
\eqref{kog2} is nondegenerate on the free parts of the groups. (This can be
proved by using the Atiyah--Singer index formula and by passing with the use of
the Chern character to cohomology.) Comparing \eqref{kog1} and \eqref{kog2}, we
can make a guess that the \emph{$K$-homology} groups of $M$ can be defined in
terms of topological $K$-theory:
\begin{equation}
\label{muma} K_0(M)\stackrel{\operatorname{def}}= K^0_c(T^*M).
\end{equation}
Unfortunately, this definition has a significant drawback, since the right-hand
side of the formula does not make sense for a singular space $M$. Nevertheless,
the right-hand side of \eqref{muma} can be defined for an arbitrary compact
space provided that we interpret $K_c^0(T^*M)$ as the group of stable homotopy
classes of elliptic operators on $M$.

Atiyah \cite{Ati4} suggested the following abstract notion of an
elliptic operator.
\begin{definition}
{ A \emph{generalized elliptic operator} over a compact space $X$ is a triple
$(D,H_1,H_2)$, where
$$
D:H_1\longrightarrow H_2
$$
is a Fredholm operator acting in the Hilbert spaces $H_1,H_2$, these
spaces are modules over the $C^*$-algebra $C(X)$ of continuous
complex-valued functions on $X$, and $D$ almost commutes with the module
structure; i.e.,
$$
[D,f]\in \mathcal{K}
$$
for all $f\in C(X)$, where $\mathcal{K}$ is the space of compact operators
from $H_1$ to $H_2$. }
\end{definition}
Note that the compactness of the commutator in this definition originates from
the fundamental property of differential operators on manifolds: for a smooth
function $f$, the order of the commutator $[D,f]$ is at most the order of $D$
minus one.

Atiyah showed that generalized elliptic operators on a compact space $X$ define
elements in the $K$-homology group $K_0(X)$ of the space, where $K_0$ is the
generalized homology theory dual to  topological $K$-theory. Atiyah also
conjectured that the $K$-homology groups can be defined solely in terms of
generalized elliptic operators. This conjecture was proved by Kasparov
\cite{Kas1} and independently by Brown--Douglas--Fillmore \cite{BDF1}. It
turned out more convenient to work with a certain modification of Atiyah's
original definition.

We start from the description of the even group $K_0(X)$. (The odd groups
$K_1(X)$ will be introduced later.)\vspace{1mm}

{\bf 2. Even cycles and the group $K_0(X)$.} This group is generated by
so-called \emph{even cycles}.
\begin{definition}
\label{d8} An \emph{even cycle} in the $K$-homology of a space $X$ is a
pair $(F,H)$ given by a $\mathbb{Z}_2$-graded Hilbert space
$$
H=H_0\oplus H_1,
$$
where the components $H_{0,1}$ are $*$-modules over the $C^*$-algebra $C(X)$,
and an odd (with respect to the grading) bounded operator
$$
F:H\to H.
$$
The operator and the module structure have the properties
\begin{equation}
\label{kompa} f(F-F^*)\sim 0, \quad f(F^2-1)\sim 0,\quad[F,f]\sim 0 \quad
\text{for all } f\in C(X),
\end{equation}
where $\sim$ means equality modulo compact operators.
\end{definition}

The even $K$-homology group $K_0(X)$ can be obtained by introducing some
equivalence relation on the cycles. The simplest one is the \emph{stable
operator homotopy}. (See \cite{Bla1}, chapter VIII, where one can find a
number of other equivalence relations on cycles. These relations are
pairwise equivalent.) Let us define this relation.

Two cycles are said to be \emph{isomorphic} if the corresponding $C(X)$-modules
$H$ are isomorphic and the operators $F$ coincide under this isomorphism. Two
cycles are \emph{homotopic} if they become isomorphic after some homotopy of
the operators $F$. A \emph{trivial cycle} is a cycle for which all relations in
\eqref{kompa} are satisfied exactly. Finally, two cycles are \emph{stably
operator homotopic} (for short, stably homotopic) if they become homotopic
after adding some trivial cycles to each of them.

One can show that the set of stable homotopy classes of even cycles is an
abelian group with respect to the direct sum. This group is denoted by
$K_0(X)$. It is called the $K$-\emph{homology group} of space $X$.
\begin{remark}
The $K$-homology groups behave covariantly for continuous mappings. Consider a
continuous mapping $f:X\to Y$. Then a cycle $(F,H)$ over $X$ can be treated as
a cycle over $Y$ if the $C(Y)$-module structure on $H$ is defined as the
composition of the induced mapping $f^*:C(Y)\to C(X)$ with the original
$C(X)$-module structure on $H$. This mapping of cycles induces a homomorphism
$$
f_!:K_0(X)\longrightarrow K_0(Y).
$$
\end{remark}

After these definitions, it is almost a tautology to say that ordinary elliptic
operators on closed manifolds define $K$-homology elements. However, because of
its importance, we describe the corresponding construction in detail.

Namely, consider an elliptic pseudodifferential operator
$$
D:C^\infty(M,E)\longrightarrow C^\infty(M,F)
$$
acting on sections of some bundles $E$ and $F$ on a closed manifold $M$. Then
its partial isomorphism part
$$
D'=(1+D^*D)^{-m/2}D
$$
in the polar decomposition defines a  Fredholm operator
\[
D':L^2\left( M,E\right) \longrightarrow L^2\left( M,F\right) ,
\]
where both spaces are $C\left( M\right) $-modules. (The module structure is the
pointwise product of functions.) Moreover, $D'$ commutes with the module
structure up to compact operators. For smooth functions, this follows from the
composition formulas for pseudodifferential operators. Then the general case
follows by continuity. We define the matrix operator
\[
F=\left(
\begin{array}{cc}
0 & D'^{*} \\
D' & 0
\end{array}
\right) .
\]
Then $F$ is a self-adjoint odd operator on the naturally
$\Bbb{Z}_2$-graded $C\left( M\right) $-module $H=L^2\left(
M,E\right)\oplus L^2\left( M,F\right) $. Thus, the pair $(F,H)$ is an
even cycle. The corresponding element in $K$-homology is denoted by
\begin{equation}
\label{khom} [D]\in K_0(M).
\end{equation}

{\bf 3. The quantization mapping.} This construction can be interpreted
in quite a different way. Namely, the principal symbol of $D$
$$
\sigma(D):\pi^*E\longrightarrow \pi^* F, \quad \pi:T^*M\to M,
$$
is a vector bundle isomorphism over $T^*M$ except for the zero section.
The corresponding difference element in $K$-theory with compact
supports is denoted by
\begin{equation}
\label{diff}[\sigma(D)]\in K^0_c(T^*M).
\end{equation}
Considering two elements \eqref{khom} and \eqref{diff} together, one can
readily show that there is a well-defined mapping
\begin{equation}
\label{quant}
\begin{array}{rcl}
Q:K_c^0(T^*M)& \longrightarrow & K_0(M), \vspace{1mm}\\
\left[\sigma (D)\right] & \mapsto & [D],
\end{array}
\end{equation}
which sends symbols to the corresponding operators. This mapping is
naturally called the \emph{quantization mapping}. It is well defined on
homotopy classes, since operators with the same symbol are homotopic.

{\bf 4. Odd cycles.} We have defined the group $K_0(X)$ in terms of even
cycles. The odd group $K_1(X)$ is generated by odd cycles.

\begin{definition}
An {\em odd cycle} in the $K$-homology of a space $X$ is given by a
Hilbert space $H$ that is a $*$-module over $C(X)$ and a bounded operator
$$
F:H\to H
$$
such that
\begin{equation}
\label{kompa1} f(F-F^*)\sim 0, \quad f(F^2-1)\sim 0,\quad[F,f]\sim 0
\end{equation}
for an arbitrary function $f$.
\end{definition}
The only difference between odd and even cycles is in the
$\mathbb{Z}_2$-grading structure.

The set of stable homotopy classes of odd cycles on $X$ is denoted by
$K_1(X)$. Let us show how odd cycles arise from elliptic \emph{self-adjoint}
operators on closed manifolds. Consider an elliptic self-adjoint operator
$$
A:C^\infty(M,E)\longrightarrow C^\infty(M,F).
$$
We can extend it to the $L^2$-spaces of sections (if the order of $A$ is zero).
This gives the odd cycle $((1+A^2)^{-m/2}A,L^2(M,E))$. The corresponding
$K$-homology class is denoted by
$$
[A]\in K_1(M).
$$
On the other hand, the principal symbol of $A$ defines a difference element
$$
[\sigma(A)]\in K^1_c(T^*M).
$$
Let us recall its definition (see \cite{APS3}). The principal symbol of
an elliptic self-adjoint operator is an invertible Hermitian endomorphism of a bundle
over the cosphere bundle $S^*M$. Hence, over $S^*M$ we can consider the
vector bundle denoted by $\im\sigma(\Pi_+(A))$ and generated at a point
$(x,\xi)\in S^*M$ by the eigenvectors of the symbol $\sigma(A)(x,\xi)$
corresponding to positive eigenvalues. Then the difference element is
defined by the formula
$$
[\sigma(A)]=\partial[\im \Pi_+(A)]\in K^1_c(T^*M),
$$
where
$$
\partial:K(S^*M)\to K^1_c(T^*M)
$$
is the coboundary mapping in the $K$-theory of the pair $S^*M\subset
B^*M$. (Here $B^*M$ is the bundle of unit balls in $T^*M$ with respect
to some Riemannian metric, and the cosphere bundle $S^*M$ is realized as
its boundary.)

By analogy with the even case, there is a well-defined quantization
mapping
\begin{equation}
\label{quant1}
\begin{array}{rcl}
Q:K_c^1(T^*M)& \longrightarrow & K_1(M), \\
\left[\sigma (A)\right] & \mapsto & [A].
\end{array}
\end{equation}
This is well defined, since the coboundary $\partial$ induces an isomorphism
\begin{equation*}
 K(S^*M)/K(M)\simeq K^1_c(T^*M).
\end{equation*}

A classical theorem (see \cite{Kas3}) claims that the quantization mapping is
an isomorphism. This isomorphism is called the \emph{Poincar\'e isomorphism}.

\begin{theorem}
\label{pua} \emph{(The Poincar\'e isomorphism.)} The quantization mappings
given by \eqref{quant} and \eqref{quant1} define isomorphisms
$$
Q: K_c^*(T^*M)\longrightarrow K_*(M),
$$
where the index $*$ is either $0$ or $1$.
\end{theorem}
The proof can be found in \cite{Kas3}.

\subsection{Duality and the topological index}
As was already mentioned, $K$-homology is dual to topological $K$-theory. The
most important manifestation of this duality is the pairing
\begin{equation}
\langle,\rangle: K_i(X)\times K^i(X)\longrightarrow \mathbb{Z}, \quad i=0,1,
\label{po1}
\end{equation}
which we define in this subsection. The pairing \eqref{po1} will be defined as
the composition of the product
\begin{equation}
K_i(X)\times K^i(X)\longrightarrow K_{0}(X) \label{proda3}
\end{equation}
with the mapping
$$
p_!:K_0(X)\longrightarrow K_0(pt)=\mathbb{Z}
$$
induced by the projection of $X$ into a one-point space.

Let us start by saying that $p_!$ takes a Fredholm operator to its
analytic index. Indeed, an even cycle over a point is just a Fredholm
operator, and Fredholm operators have only one stable homotopy invariant,
namely, the index. This shows that $K_0(pt)=\mathbb{Z}$.

It remains to define the product \eqref{proda3}. Let us consider it for
$i=j=0$.

The formula for the product (e.g., see \cite{Ati4}) for a general element
of the group $K_0(X)$ mimics the construction of an operator with
coefficients in a vector bundle. Namely, this mapping takes a cycle
$(F,H)$ and a vector bundle $G$ to the even cycle
$$
(F\otimes 1_{\mathbb{C}^N}, 1_{H}\otimes P(H\otimes \mathbb{C}^N)),
$$
where $P:\mathbb{C}^N\to \mathbb{C}^N$ is some projection over $X$
defining $G$ and
$$
1_H\otimes P(H\otimes \mathbb{C}^N)\subset H\otimes \mathbb{C}^N
$$
is the range of the projection $1_H\otimes P$. Thus, the pairing
$ \langle[D],[G]\rangle$ of an elliptic operator $D$ with a bundle $G$
can be computed by applying the Atiyah--Singer formula to $D\otimes 1_G$:
$$
\langle[D],[G]\rangle=\ind_t (D\otimes 1_G).
$$

Let us also mention that the product of odd groups
$$
K_1(X)\times K^1(X)\longrightarrow K_0(X)
$$
is defined in the spirit of the theory of Toeplitz type operators
(see~\cite{BaDo2}).

Comparing the constructions of the present section with the de Rham
theory, one can make the following glossary of similar terms:
$$\renewcommand\arraystretch{1.5}
\begin{tabular}{|c|c|}\hline
de Rham theory & Elliptic theory \\ \hline \hline
cohomology $H^*(M)$ & topological $K$-theory $K^*(M)$ \\ \hline
homology $H_*(M)$ & $K$-homology $K_*(M)$ \\ \hline
cocycle $\omega$, $d\omega=0$ & vector bundle $E$ \\ \hline
cycle $\gamma,$ $\partial\gamma=0$ & elliptic operator $D$ \\ \hline
integral $\int_\gamma\omega$ & index $\ind(D\otimes 1_E)$ \\ \hline
\end{tabular}
$$

\subsection{Poincar\'e duality on manifolds with boundary. Absolute and relative
cycles}

Let us now consider duality for manifolds with boundary.

Let $M$ be a compact smooth manifold with boundary $\partial M$. In this
case, one also has Poincar\'e duality, frequently called
Poincar\'e--Lefschetz duality. In (co)homology (for oriented $M$), the
duality amounts to two group isomorphisms
$$
H^i(M)\longrightarrow H_{n-i}(M,\partial M), \quad H^i(M,\partial
M)\longrightarrow H_{n-i}(M), \quad n=\dim M.
$$
Thus, duality relates the ordinary groups to the so-called relative
groups.

Let us define similar isomorphisms in $K$-theory. For a manifold $M$ with
boundary, one has two natural $K$-homology groups: the ordinary group
$K_*(M)$ and the so-called \emph{relative $K$-homology group}
$K_*(M,\partial M)$. The relative groups are defined as follows.

It is clear that in Definition \ref{d8} we have used only the algebra $C(M)$
of functions on $M$. At the same time, a manifold with boundary has another
natural $C^*$-algebra, namely, the algebra $C_0(M\setminus\partial M)$ of
functions vanishing at the boundary. We denote the group generated by cycles
over $C_0(M\setminus\partial M)$ by $K_*(M,\partial M)$ and call it the
\emph{relative $K$-homology} of the manifold.

Let us show how elements of these groups arise from elliptic operators on
manifolds with boundary. On the analogy with homology theory, we refer to
cycles over the algebra $C(M)$ as \emph{absolute cycles} and cycles over
$C_0(M\setminus\partial M)$ as \emph{relative cycles}.

We first describe the ordinary (absolute) cycles. \vspace{1mm}

\textbf{1. Elliptic operators and absolute cycles}. In this case, we
consider elliptic operators $D$ induced near the boundary by vector
bundle isomorphisms. Just as in the case of closed manifolds, such
operators almost commute with functions $f\in C(M)$:
$$
[D,f]\sim 0.
$$
Therefore, the construction of the previous subsection can be carried out
word for word in this case; i.e., $D$ defines a $K$-homology class
$$
[D]\in K_*(M).
$$
On the other hand, the principal symbol of such an operator is an
isomorphism on $T^*M$ \emph{everywhere near the boundary}. Hence, the
operator $D$ has a difference element
$$
[\sigma(D)]\in K^*_c(T^*(M\setminus\partial M)).
$$
\begin{remark}
{This class of zero-order operators induced in a neighborhood of the boundary
by vector bundle isomorphisms naturally arises in index theory of classical
boundary value problems. Namely, the classical procedure of order reduction
(see \cite{Hor3} or \cite{SaScS4}) reduces a boundary value problem for a
differential operator to a pseudodifferential operator of order zero of
precisely the form considered here. Moreover, the reduction preserves the
index.}
\end{remark}

\textbf{2. Elliptic operators and relative cycles.} Relative cycles arise
as follows. An elliptic operator $D$ of order one on $M$ defines an
element
\[
\left[ D\right] \in K_*\left( M, \partial M\right),
\]
which can be constructed as follows. Consider an embedding
$M\subset\widetilde{M}$ of $M$ in some closed manifold $\widetilde{M}$ of the
same dimension. Let $\widetilde{D}$ be an arbitrary extension of $D$ to this
closed manifold as an elliptic operator. On $\widetilde{M}$ we consider the
zero-order operator
\[
\widetilde{F}=\left( 1+\widetilde{D}^{*}\widetilde{D}\right) ^{-1/2}
\widetilde{D}.
\]
The restriction of this operator to $M$ is defined as
\begin{equation}
F=i^{*}\widetilde{F}i_{*}:L^2\left( M,E\right) \longrightarrow L^2\left(
M,F\right) , \label{relo}
\end{equation}
where $i_{*}:L^2\left( M\right) \rightarrow L^2\left( \widetilde{M}
\right) $ is the extension by zero and $i^{*}:L^2\left(
\widetilde{M}\right) \rightarrow L^2\left( M\right) $ is the restriction
operator.

If $D$ is symmetric, then one can readily verify that $F$ satisfies the
properties
\[
\begin{array}{c}
F-F^{*}\sim 0,\quad f\left( F^2-1\right) \sim 0, \quad \left[ F,f\right] \sim 0
\end{array}
\]
for a function $f\in C_0\left(M\backslash \partial M\right)$ vanishing on the
boundary. These relations show that $F$ defines an element in $K_1\left(
M,\partial M \right) $.

If $D$ is not symmetric, then instead of the operator $F$ one considers the
corresponding matrix operator as in Subsection \ref{ss51}. In this case, the
pair $(F,H)$ defines an element of the group $K_0\left(M,\partial
M\right).$\vspace{1mm}

\textbf{3. The quantization mapping and the Poincar\'e isomorphism.} For
elliptic operators defining absolute cycles on manifolds with boundary,
we obtain two elements
$$
[\sigma(D)]\in K^*_c(T^*(M\setminus\partial M)), \quad [D]\in K_*(M).
$$
Similarly, operators corresponding to relative cycles also define a pair
of elements
$$
[\sigma(D)]\in K^*_c(T^*M), \quad [D]\in K_*(M,\partial M).
$$
(Here $[\sigma(D)]$ is the Atiyah--Singer difference element of $D$.) One can
show that the corresponding \emph{quantization mappings}
\begin{equation}
\label{quant3}
\begin{array}{ccc}
Q:K_c^*(T^*(M\setminus \partial M)& \longrightarrow & K_*(M), \vspace{1mm}\\
K_c^*(T^*M)& \longrightarrow & K_*(M,\partial M), \vspace{1mm}\\
\left[\sigma (D)\right] & \mapsto & [D]
\end{array}
\end{equation}
are well defined. It was proved in \cite{Kas3} that both quantization
mappings are isomorphisms. They are called the \emph{Poincar\'e
isomorphisms in $K$-theory on manifolds with boundary.} \vspace{1mm}

\textbf{4. The exact sequence of a pair in $K$-homology.} It is well
known that in $K$-homology for an embedding $\partial M\subset M$ there
is a six-term exact sequence
\begin{equation}
\xymatrix{ & K_0\left(\partial M\right)\ar[r] & K_0\left( M\right)
\ar[rd] \\
\limfunc K_1\left(M,\partial M \right)\ar[ur] & & & K_0\left(M,\partial M \right).\ar[dl] \\
& K_1\left( M\right)\ar[lu] & K_1\left( \partial M\right)\ar[l] } \label{ino}
\end{equation}
It follows from the results of~\cite{BDT1} that this sequence is
isomorphic to the exact sequence of the pair $T^*M|_{\partial M}\subset
T^*M $ in topological $K$-theory:
\begin{equation}
\xymatrix{ &K_c\left( T^{*}\partial M\right) \ar[r] & K_c\left( T^{*}\left(
M\backslash \partial M\right) \right)
\ar[rd] \\
K^1_c\left( {T^{*}M} \right)\ar[ur] & & & K_c\left(
{T^{*}M} \right) .\ar[dl] \\
& K^1_c\left( T^{*}\left( M\backslash \partial M\right) \right)\ar[lu] &
K^1_c\left( T^{*}\partial M\right)\ar[l] } \label{ino1}
\end{equation}
The isomorphism can be obtained by applying the quantization mappings
term by term to the sequence~\eqref{ino1}.

\section{Appendix C. Poincar\'e duality on $\mathbb{Z}_n$-manifolds}
\label{ch5}

In this section, we construct Poincar\'e duality for $\mathbb{Z}_n$-manifolds
(see \cite{SaSt10}). Unfortunately, these manifolds with singularities have no
duality in the framework of the usual topological $K$-groups. Quite remarkably,
however, the duality can be restored if one applies the approach of
noncommutative geometry \cite{Con1} to this problem and states duality in
terms of $K$-theory of some noncommutative algebras of functions on the
corresponding spaces. We show that the approach of noncommutative geometry, at
least for $\mathbb{Z}_n$-manifolds, can be completely described in terms of
elliptic operator theory.

Before we proceed to the description of duality for $\mathbb{Z}_n$-manifolds,
we recall the main features of the passage from topological $K$-theory to
$K$-theory of algebras. This passage goes as follows.
\begin{enumerate}
\item The topological group $K^0(X)$ of a compact space $X$ can be identified
with the Grothendieck group of homotopy classes of projections in matrix
algebras over $C(X)$. Clearly, this definition makes sense for an arbitrary
unital $C^*$-algebra $\mathcal{A}$. The corresponding group is denoted by
$K_0(\mathcal{A})$. Similarly, the odd group $K^1(X)$ is identified with the
group of stable homotopy classes of unitary operators in matrix algebras over
$C(X)$. The similar group for a $C^*$-algebra $\mathcal{A}$ is denoted by
$K_1(\mathcal{A})$. Summarizing, the new groups are reduced in the commutative
case to the topological $K$-groups
$$
K^*(X)\simeq K_*(C(X)).
$$

\item $K$-homology. It is clear from Definition \ref{d8} that the definition of
the $K$-homology group of a space uses only the $C^*$-algebra of functions on
it. Therefore, the same definition for an arbitrary $C^*$-algebra $\mathcal{A}$
gives two groups $K^0(\mathcal{A})$ and $K^1(\mathcal{A})$. These groups are
generated, respectively, by even and odd cycles over $\mathcal{A}$.
\end{enumerate}

\begin{remark}
We note the change of variance of the functors (and the corresponding change in
the position of indices). For example, an analog of the topological $K$-groups
is the group $K_*(\mathcal{A})$, which is a covariant functor with respect to
algebra homomorphisms. Similarly, the $K$-homology of spaces translates into
the $K$-cohomology (a contravariant functor) of algebras.
\end{remark}

For a nonunital algebra $\mathcal{A}$, the $K_0$-group is defined as the kernel
of the natural surjective mapping
$$
K_0(\mathcal{A})\stackrel{\operatorname{def}}=\ker\left(K_0(\mathcal{A}^+)\longrightarrow
K_0(\mathbb{C})\right),
$$
where $\mathcal{A}^+=\mathcal{A}\oplus \mathbb{C}$ is the algebra with attached
unit. The mapping is induced by the algebra homomorphism
$\mathcal{A}\oplus\mathbb{C}\to \mathbb{C}$. The $K_1$-group is defined
slightly simpler:
$$
K_0(\mathcal{A})\stackrel{\operatorname{def}}=K_1(\mathcal{A}^+).
$$
The construction is similar to the definition of topological $K$-theory for
locally compact spaces (see \cite{Ati2}).

\subsection{Relative cycles. The $C^*$-algebra of a $\mathbb{Z}_n$-manifold}
\label{ss61} Since a $\mathbb{Z}_n$-manifold for $n=1$ is just a smooth
manifold with boundary, it is natural to construct duality separately for
relative and absolute cycles.

Consider a twisted $\mathbb{Z}_n$-manifold $(M,\pi)$ (see Definition \ref{d6}),
or a $\mathbb{Z}_n$-manifold for short. On $M$, we take some elliptic operator
$$
D:C^\infty(M,E)\longrightarrow C^\infty(M,F)
$$
satisfying Assumption~1. Its principal symbol defines a difference element
\begin{equation}
\label{difa}[\sigma(D)]\in K_c(\overline{T^*M}^\pi)
\end{equation}
of the topological $K$-group of the $\mathbb{Z}_n$-manifold $(T^*M,\pi)$. Note
that \eqref{difa} includes additional information (compared with the image of
$[\sigma(D)]$ in $K_c(T^*M)$) related to the structure of the space
$\overline{T^*M}^\pi$. Let us define a similar element in $K$-homology. Note
that by assumption, the operator $D$ almost commutes both with functions $f\in
C_0(M\setminus\partial M)$ and with transpositions of the sheets of $\pi$ in a
neighborhood of the boundary.

By $\mathcal{A}_{M,\pi}$ we denote the $C^*$-subalgebra of operators in
$L^2(M)$ generated by multiplications by functions $f\in C_0(M\setminus
\partial M)$ and transposition operators in a neighborhood of the boundary. If
we take a neighborhood of the boundary where transpositions of the sheets are
allowed and denote the complement of this neighborhood by $M'$, then the
algebra can be explicitly described in the form
\[
\mathcal{A}_{M,\pi }=\left\{ (u,v)\; \left|\;
\begin{array}{c}
u\in C_0\left( M^{\prime }\right), v\in C_0\left(
X\times \left( 0,1\right] ,\limfunc{End}\pi _{!}1\right)\vspace{1mm}\\
\beta \left( \left. u\right| _{\partial M^{\prime }}\right)\beta ^{-1}=\left.
v\right| _{t=1}
\end{array}
\right. \right\}
\]
as a subalgebra in the direct sum $C_0\left( M^{\prime }\right)\oplus C_0\left(
X\times \left( 0,1\right] ,\limfunc{End}\pi _{!}1\right)$. Here ${\rm
End}\pi_!1$ is the vector bundle of endomorphisms of $\pi_!1\in\Vect(X)$.

\begin{remark}
{ The algebra $\mathcal{A}_{M,\pi }$ can also be defined as the $C^*$-algebra
of the etale-groupoid \cite{BrNi1} corresponding to the extension of the
equivalence relation \eqref{eka} to the neighborhood of the boundary}.
\end{remark}

Under Assumption 1, one can prove the following result.

\begin{lemma}
The spaces $L^2(M,E)$ and $L^2(M,F)$ have a natural module structure over the
algebra $\mathcal{A}_{M,\pi}$ such that $D$ almost commutes with this module
structure. Thus, $D$ defines a class $[D]\in K^*(\mathcal{A}_{M,\pi}).$
\end{lemma}

\begin{proof}
By Assumption 1, in a neighborhood of the boundary one has vector bundle
isomorphisms
$$
E\simeq \pi^* E_0, \quad E_0\in \Vect(X\times [0,1]).
$$
In a similar way, for the spaces of sections we have
$$
C(\partial M\times [0,1],E)\simeq C(X\times [0,1],\pi_!1\otimes E_0).
$$
The latter space has a natural $C_0\left( X\times \left( 0,1\right]
,\limfunc{End}\pi _{!}1\right)$-module structure. This module structure can be
glued with the $C(M)$-module structure on the entire manifold, thus defining
the desired module structure over the algebra $\mathcal{A}_{M,\pi}$. A similar
construction applies to $F$. The desired almost commutativity of $D$ with the
elements of $\mathcal{A}_{M,\pi}$ follows from the equivariance of $D$ with
respect to transpositions of sheets of the covering.
\end{proof}

\subsection{Absolute cycles. Nonlocal operators}\label{noni}

Clearly, the relative cycles on a $\mathbb{Z}_n$-manifold were obtained in the
previous subsection by restricting the class of operators to those almost
commuting with the algebra $\mathcal{A}_{M,\pi}$. In this subsection, we
enlarge the class of absolute cycles on $M$ requiring only that the operators
in question almost commute  with elements of the algebra of continuous
functions on the singular space $\overline{M}^\pi$ but not of the larger
algebra
$$
C(M)\supset C(\overline{M}^\pi) .
$$

{\bf 1. Definition of the class of operators.} Needless to say, ordinary
absolute cycles on $M$ (represented by elliptic operators of order zero induced
by vector bundle isomorphisms near the boundary) are cycles for the algebra
$C(\overline{M}^\pi)$ as well. However, there are also more general cycles on
the space $\overline{M}^\pi$. To find them, let us replace
$C(\overline{M}^\pi)$ by a homotopic algebra of functions on $M$ that are
lifted from the base of the covering in a given neighborhood of the boundary.
We denote this algebra by the same symbol. Here is an ``heuristic derivation"
of relative cycles.

Consider some operator $D$ on $M$. In a neighborhood of the boundary, its
direct image $(\pi\times 1)_!D$ acts on the base of the covering $\pi\times
1:\partial M\times[0,1)\to X\times [0,1)$. The spaces in which the operator
acts have the usual $C(X)$-module structure (determined by the pointwise
product of functions). Therefore, $D$ almost commutes with functions in
$C(\overline{M}^\pi)$ if (a) it is a pseudodifferential operator far from the
boundary and (b) its direct image in a neighborhood of the boundary is a
pseudodifferential operator on the base. Thus, we arrive at the following
definition.

\begin{definition}
An \emph{admissible operator} on $M$ is an operator of the form
\begin{equation}
\label{alef1} D=\psi'D^{\prime }\varphi'+\psi'' \left( \pi ^{!}D^{\prime \prime
}\right) \varphi'' ,
\end{equation}
where
$$
D^{\prime }:C^\infty(M,E)\to C^\infty(M,F)
$$
is a pseudodifferential operator of order zero on $M$,
$$
D^{\prime \prime}:C^\infty(X\times[0,1],(\pi\times 1)_!E)\longrightarrow
C^\infty(X\times[0,1],(\pi\times 1)_!F)
$$
is a pseudodifferential operator on $X\times[0,1]$, the cutoff functions
$\varphi'$ and $\psi' $ vanish in a neighborhood of the boundary, and
$\varphi''$ and $\psi''$ vanish in a neighborhood of the subset
$X\times\{1\}\subset X\times[0,1]$. Finally, the operator $D''$ is induced by a
vector bundle isomorphism in the neighborhood $X\times[0,1/2)\subset
X\times[0,1]$ on the cylinder.
\end{definition}

\begin{remark}
Inverse images like $\pi^!D''$ are in general \emph{nonlocal} operators. By
way of illustration, consider the trivial covering
\[
\partial M=\stackunder{n\;\;{{\text{copies}}}}{\underbrace{X\sqcup X\sqcup ...\sqcup
X}}\longrightarrow X
\]
of $n$ copies of $X$. On $\partial M$, we choose a trivial bundle $E=\Bbb{C}$.
Then $\pi _{!}E=\Bbb{C}^n$; the direct image of a differential operator on the
total space is always a diagonal operator
\[
\pi _{!}P=\limfunc{diag}\left( \left. P\right| _{X_1},...,\left. P\right|
_{X_n}\right) :C^\infty \left( X,\Bbb{C}^n\right) \longrightarrow C^\infty
\left( X,\Bbb{C}^n\right) ,
\]
whereas the inverse image of an operator on the base is a matrix operator of
general form. Therefore, near the boundary we admit usual pseudodifferential
operators as well as operators corresponding to transpositions of values of a
function on the sheets of the covering. In a sense, admissible operators are
obtained as an extension of the algebra of pseudodifferential operators by the
nonlocal operators corresponding to elements of $\mathcal{A}_{M,\pi}$.
\end{remark}

{\bf 2. Symbols of admissible operators.} One can readily define the symbol of
an admissible operator and state the ellipticity condition. Namely, the symbol
of an admissible operator $D$ is a pair $\left( \sigma_M ,\sigma _X\right) $ of
usual symbols, where
\[
\sigma _M:p_M^{*}\left. E\right| _{M'}\longrightarrow p_M^{*}\left. F\right|
_{M'},\quad p_M:S^{*}M\longrightarrow M,
\]
is a symbol on the part of the manifold where the operator is
pseudodifferential and
\[
\sigma _X:p_X^{*}\left. \left(\pi_!E\right)\right|
_{{X\times[0,1]}}\longrightarrow p_X^{*}\left. \left(\pi_!F\right)\right|
_{X\times[0,1]},\quad p_X:S^{*}(X\times[0,1])\longrightarrow X\times[0,1],
\]
is a symbol on the cylinder $X\times[0,1]$. The symbols satisfy the
compatibility condition
\[
\pi _{!}\left( \left. \sigma _M\right| _{\partial M'}\right) =\sigma _X.
\]
An operator is said to be \emph{elliptic} if both components of its symbol are
invertible. An elliptic operator is Fredholm in appropriate Sobolev spaces.

Summarizing, we see that admissible elliptic operators $D$ on a
$\mathbb{Z}_n$-manifold define absolute cycles
$$
[D]\in K_*(\overline{M}^\pi).
$$
Let us show that the symbol of an elliptic admissible operator defines an
element in $K$-theory. In other words, let us define an analog of the
Atiyah--Singer difference construction for admissible elliptic operators.
\vspace{1mm}

\textbf{3. The difference construction for admissible operators.} Let us cut
$M$ into two parts
\[
M^{\prime }=M\backslash \left\{ \partial M\times \left[ 0,1\right) \right\}
\;\qquad \text{and\qquad }\partial M\times \left[ 0,1\right] .
\]
Then the symbol $\sigma \left( D\right) $ of an admissible elliptic operator
$D$ is naturally represented as a pair $( \sigma_M,\sigma_X)$. Each symbol in
this pair has the corresponding difference element
\[
\left[ \sigma_M\right] \in K^0_c\left( T^{*}M^{\prime }\right) ,\qquad \left[
\sigma_X\right] \in K^0_c\left( T^{*}\left( X\times \left( 0,1\right] \right)
\right).
\]
In the latter case, we use the difference construction for an absolute cycle
with symbol $\sigma_X$ of order zero.

However, a single element of some topological $K$-group cannot be constructed
from these data, since the manifolds $T^{*}M^{\prime }$ and $T^{*}\left(
X\times \left( 0,1\right] \right) $ cannot be glued smoothly (their boundaries
may be nondiffeomorphic). Nevertheless, the pasting can be done if we glue the
algebras of these manifolds rather than the manifolds themselves.

Let $\mathcal{A}_{T^*M,\pi}$ be the $C^*$-algebra of the
$\mathbb{Z}_n$-manifold $(T^*M,\pi)$. It is the subalgebra
\[
\mathcal{A}_{T^{*}M,\pi }\subset C_0\left( T^{*}M^{\prime }\right) \oplus
C_0\left( T^{*}\left( X\times \left( 0,1\right] \right) ,\limfunc{End}
p^{*}\pi _{!}1\right) ,
\]
defined by the compatibility condition
\[
\mathcal{A}_{T^{*}M,\pi }=\left\{ u\oplus v\;|\;\beta \left. u\right|
_{\partial M^{\prime }}\beta ^{-1}=\left. v\right| _{t=1}\right\}
\]
on the boundaries.

The \emph{difference construction} for admissible operators is a mapping
\begin{equation}
\chi :\limfunc{Ell}\left( M,\pi \right) \longrightarrow K_0\left(
\mathcal{A}_{T^{*}M,\pi }\right) \label{stk}
\end{equation}
of the group $\limfunc{Ell}\left( M,\pi \right)$ of stable homotopy classes
of elliptic admissible operators on $M$ into the $K_0$-group of the algebra
$\mathcal{A}_{T^{*}M,\pi }$.

Let us explicitly describe the element $\chi \left[ D\right] $ corresponding to
an elliptic admissible operator
\[
D:C^\infty \left( M,E\right) \longrightarrow C^\infty \left( M,{F}\right)
\]
with symbol $\sigma(D)=(\sigma_M,\sigma_X)$. By stabilization, one can always
assume that the bundle $F$ is trivial: $F=\mathbb{C}^k$. Now let us choose
some embeddings of $E$ and ${\Bbb{C}^k}$ in trivial vector bundles and some
projections
\[
P_E,P_{\Bbb{C}^k}:\Bbb{C}^{N+L}\longrightarrow \Bbb{C}^{N+L}
\]
defining $E$ and ${\Bbb{C}^k}$ as
\[
E\simeq \func{Im}P_E\subset \Bbb{C}^N\oplus 0,\qquad {\Bbb{C}^k}\simeq \func{Im}
P_{\Bbb{C}^k}\subset 0\oplus \Bbb{C}^L.
\]
Let $P_{\pi _{!}E},P_{\pi _{!}{\Bbb{C}^k}}$ be the direct images of these
projections in a neighborhood of the boundary. Clearly, these projections
define the direct images of the corresponding subbundles.

Then the difference construction of $D$ is defined by the formula
\[
\chi \left[ D\right] =\left[ P_1\oplus P_2\right] -\left[ P_{\Bbb{C}^k}\oplus P_{\pi
_{!}{\Bbb{C}^k}}\right] ,
\]
where the projection $P_1$ over $T^*M^{\prime }$ is defined by
\begin{equation}
\left\{
\begin{array}{cc}
P_E\cos ^2\left| \xi \right| + P_{\Bbb{C}^k}\sin ^2\left| \xi \right| +
\left(\sigma^{-1}_M P_{\Bbb{C}^k} +\sigma_M P_E \right) \sin \left| \xi
\right| \cos \left| \xi \right|
, & \left| \xi
\right| \leq \pi /2,\vspace{1mm} \\
P_{\Bbb{C}^k}, & \left| \xi \right| >\pi /2
\end{array}
\right. \label{strr}
\end{equation}
(we assume that the symbol $\sigma_M $ is homogeneous of order zero in the
covariable $\xi$), and the projection $P_2$ over $T^*(X\times \left[
0,1\right])$ is defined by one of the expressions
\[
\begin{array}{cc}
P_{\pi _{!}E}\cos ^2\left| \xi \right| +P_{\pi _{!}{\Bbb{C}^k}}\sin ^2\left| \xi
\right| +1/2 \left({\sigma }_X^{-1} P_{\pi _{!}{\Bbb{C}^k}}+ {\sigma }_X P_{\pi
_{!}E} \right)\sin 2\left| \xi \right|
, \vspace{1mm} \\
P_{\pi _{!}E}\cos ^2\varphi +P_{\pi _{!}{\Bbb{C}^k}}\sin ^2\varphi +
1/2\left({\sigma }^{-1}_X P_{\pi _{!}{\Bbb{C}^k}}+ {\sigma }_X P_{\pi _{!}E}
\right) \sin 2\varphi,
\vspace{1mm} \\
P_{\pi _{!}{\Bbb{C}^k}}.
\end{array}
\]
Here the first case is taken  for the parameter values $x^{\prime }\in X\times
\left[ 1/2,1\right] ,\left| \xi \right| \leq \pi /2 $, the second case for
$x^{\prime }\in X\times \left[ 0,1/2\right] ,\left| \xi \right| <\pi t$, and
the third expression applies to the remaining points. Here for brevity we write
\[
\varphi =\left| \xi \right| +\pi /2\left( 1-2t\right) .
\]

\begin{remark}
Geometrically, these projections define a subbundle that coincides over the
zero section (i.e. for $\left| \xi \right| =0$) with $E\subset \Bbb{C}^{N+L}$,
coincides for $\left| \xi \right| \geq \pi /2$ with the orthogonal subbundle
${\Bbb{C}^k}\subset \Bbb{C}^{N+L}$, and is a rotation of one of the subbundles
towards the other via the isomorphisms defined by the elliptic symbol $\sigma
\left( D\right)$ at the intermediate points.
\end{remark}

\begin{proposition}
\label{sem1}The difference construction $\emph{(\ref{stk})} $ is a well-defined
group isomorphism.
\end{proposition}

\begin{proof}
The difference construction $\chi $ preserves the equivalence relations on the
groups $\limfunc{Ell}^0\left( M,\pi \right) $ and $K\left(
\mathcal{A}_{T^{*}M,\pi }\right) $. This can be proved by observing that a
homotopy of operators gives a homotopy of symbols and hence a homotopy of the
corresponding projections $P_{1,2}$. Furthermore, $\chi \left[ D\right] $ is
independent of the choice of embeddings in a trivial bundle, since all such
embeddings are pairwise homotopic. The proof of the fact that $\chi $ defines a
one-to-one mapping can be obtained along the same lines.
\end{proof}

\subsection{The Poincar\'e isomorphism}
Using the difference construction of the previous subsection, we define
quantization mappings on $\mathbb{Z}_n$-manifolds: one for absolute cycles
$$
K_{*}\left( \mathcal{A}_{T^{*}M,\pi }\right) \longrightarrow K_{*}\bigl(
\overline{M}^\pi \bigr) ,
$$
corresponding to the theory of admissible nonlocal operators, and another for
relative cycles
$$
K_c^{*}\bigl( \overline{T^{*}M}^\pi \bigr) \longrightarrow K^{*}\left(
\mathcal{A}_{M,\pi }\right),
$$
corresponding to the operators discussed in Section \ref{seczn} and
Subsection \ref{ss61}.

\begin{theorem}
{\em (The Poincar\'e isomorphism.)} The quantization mappings on
$\mathbb{Z}_n$-manifolds are isomorphisms.
\end{theorem}

\begin{proof}
1) The algebra $\mathcal{A}_{T^{*}M,\pi }$ has the ideal
\[
I=C_0\left( T^{*}\left( X\times \left( 0,1\right) \right) ,\limfunc{End}
p^{*}\pi _{!}1\right)
\]
with quotient $\mathcal{A}_{T^{*}M,\pi }/I\simeq C_0\left( T^{*}M\right).$
Consider the exact sequence
\begin{equation}
\rightarrow K_c\left( T^{*}X\right) \stackrel{ }{\rightarrow }K_0\left(
\mathcal{A}_{T^{*}M,\pi }\right) \rightarrow K_c\left( T^{*}M\right)
\rightarrow K_c^1\left( T^{*}X\right) \rightarrow  \label{exac}
\end{equation}
of this pair. (To obtain this sequence, we use the natural isomorphisms
$$
K_{*}\left( C_0\left( Y,\limfunc{End}G\right) \right) \simeq K_{*}\left(
C_0\left( Y\right) \right) \simeq K^{*}\left( Y\right)
$$
for a vector bundle $G\in\limfunc{Vect}\left(Y\right)$.) Now consider the
commutative diagram
\[
\begin{array}{ccccccccc}
\rightarrow & K_c^*\left( T^{*}X\right) & \rightarrow & K_*\left( \mathcal{A}
_{T^{*}M,\pi }\right) & \rightarrow & K^*_c\left( T^{*}M\right) & \rightarrow
& K_c^{*+1}\left( T^{*}X\right) & \rightarrow \\
& \downarrow & & \downarrow & & \downarrow & & \downarrow & \\
\rightarrow & K_*\left( X\right) & \rightarrow & K_*\bigl( \overline{M}^\pi
\bigr) & \rightarrow & K_*\left( M,\partial M\right) & \rightarrow &
K_{*+1}\left( X\right) & \rightarrow
\end{array}
\]
The second row is the exact sequence of the pair $X\subset \overline{M}^\pi $
in $K$-homology. The vertical arrows (except for the second) are Poincar\'e
isomorphisms on closed manifolds and manifolds with boundary. Thus, by the
$5$-lemma, the mapping
\[
K_{*}\left( \mathcal{A}_{T^{*}M,\pi }\right) \longrightarrow K_{*}\bigl(
\overline{M}^\pi \bigr)
\]
is an isomorphism as well.

2) The proof of the second isomorphism can be carried out in a similar way. In
this case, the following diagram is relevant:
\[
\begin{array}{ccccccccc}
\!\!\leftarrow\!\! & K^1\left( T^{*}X\right) & \!\!\leftarrow\!\! & K^0\bigl(
\overline{ T^{*}M}^\pi \bigr) & \!\!\leftarrow\!\! & K^0\left( T^{*}\left(
M\backslash
\partial M\right) \right) & \!\!\leftarrow\!\! & K^0\left( T^{*}X\right) & \leftarrow \\
& \downarrow & & \downarrow & & \downarrow & & \downarrow & \\
\!\!\leftarrow\!\! & K_1\left( X\right) & \!\!\leftarrow\!\! & K^0\left(
\mathcal{A}_{M,\pi }\right) & \!\!\leftarrow\!\! & K_0\left( M\right) &
\!\!\leftarrow\!\! & K_0\left( X\right) & \leftarrow
\end{array}
\]
(Here the first row is the exact sequence of the pair $\Bbb{R}\times
T^{*}X\subset \overline{T^{*}M}^\pi $.)
\end{proof}

\subsection{Poincar\'e duality}

\hspace{1mm}

{\bf 1. Definition of the pairing.} Let us now define  Poincar\'e--Lefschetz
duality on $\mathbb{Z}_n$-manifolds following the same scheme as in the case of
smooth manifolds. Since we have relative and absolute cycles, two dualities are
expected:
\begin{equation}
K_c^i\bigl( \overline{T^{*}M}^\pi \bigr) \times K_i\left( \mathcal{A}_{M,\pi
}\right) \longrightarrow \Bbb{Z} \label{aa}
\end{equation}
and
\begin{equation}
K_i\left( \mathcal{A}_{T^*M,\pi }\right) \times K^i\bigl(\overline{M}^\pi
\bigr) \longrightarrow \Bbb{Z}. \label{aa2}
\end{equation}
To save space, we consider only the first duality. The second can be considered
in a similar way.

Let us define the pairing \eqref{aa} as follows: we act on the first argument
by the Poincar\'e isomorphism
$$
Q: K_c^{*}\bigl( \overline{T^{*}M}^\pi \bigr) \longrightarrow K^{*}\left(
\mathcal{A}_{M,\pi }\right)
$$
and then apply the index pairing
$$
K^{*}\left(
\mathcal{A}_{M,\pi }\right)\times K_*\left( \mathcal{A}_{M,\pi }\right)
\longrightarrow \mathbb{Z}
$$
of $K$-groups of opposite variance.
\begin{theorem}
\label{pda} {\em (Poincar\'e duality.)}\cite{SaSt10} On a
$\mathbb{Z}_n$-manifold, the pairings
\begin{equation}
K_c^i\bigl( \overline{T^{*}M}^\pi \bigr) \times K_i\left( \mathcal{A}_{M,\pi
}\right) \longrightarrow \Bbb{Z}, \quad i=0,1, \label{aa3}
\end{equation}
are nonsingular on the free parts of the groups.
\end{theorem}

\begin{proof} Fixing the first argument of the pairing, we obtain a mapping
\[
K^i_c\bigl( \overline{T^{*}M}^\pi \bigr) \otimes \Bbb{Q}\longrightarrow
K_i^{\prime }\left( \mathcal{A}_{M,\pi }\right) ,
\]
where we write $G^{\prime }=\limfunc{Hom}\left( G,\Bbb{Q}\right) $. This
mapping occurs in the commutative diagram
\begin{equation*}
\begin{CD}
K^1_c\left( T^{*}\!X, \Bbb{Q}\right)&\,\longleftarrow\,&
K_c^0\bigl(\overline{T^{*}\!M}^\pi ,\Bbb{Q}\bigr)&\,\longleftarrow\,&
K^0_c\left(T^{*}\!\left( M\backslash \partial M\right)
,\Bbb{Q}\right)&\,\longleftarrow\,& K_c^0\left( T^{*}\!X,\Bbb{Q}\right)
\\
@VVV @VVV @VVV @VVV
 \\
K^{1\prime }\left( X\right)&\longleftarrow&
K_0^{\prime }\left( \mathcal{A}_{M,\pi }\right)&\longleftarrow&
K^{0\prime }\left( M\right)&\longleftarrow&
K^{0\prime }\left( X\right).
\end{CD}
\end{equation*}
All vertical arrows except for the second are isomorphisms (by Poincar\'e
duality on smooth manifolds). Thus, by the 5-lemma the remaining mapping is
also an isomorphism. Hence, (\ref{aa3}) is nondegenerate in the first factor.

The nondegeneracy in the second factor can be proved in a similar way.
\end{proof}

{\bf 2. An application to ${\rm spin}^c$-manifolds.} Consider a pair $(M,\pi)$,
where $M$ has a $spin^c$-structure induced over the boundary by a
$spin^c$-structure on the base $X$ of the covering $\pi$. Then the group
$K_c^{*}\bigl( \overline{T^{*}M}^\pi \bigr) $ is a free $K^{*+n}\bigl(
\overline{M}^\pi \bigr) $-module with one generator (where $n=\dim M$). A
generator is given by the difference element
\[
\left[ \sigma \left( D\right) \right] \in K^n_c\bigl( \overline{T^{*}M}^\pi
\bigr)
\]
of the Dirac operator on $M$. (This can be proved by analogy with the case of
closed manifolds; e.g., see \cite{LaMi1}.) Therefore, one can define the
Poincar\'e duality pairing
\[
K^{*+n}\bigl( \overline{M}^\pi \bigr) \times K_{*}\left( \mathcal{A}_{M,\pi
}\right) \longrightarrow \Bbb{Z}
\]
using the composition $K^{*+n}\bigl( \overline{M}^\pi \bigr) \rightarrow
K_c\bigl( \overline{T^{*}M}^\pi \bigr)$. It follows from Theorem \ref{pda} that
this pairing is nonsingular on the free parts of the groups.

{\bf 3. Computation of the pairing.} Our aim is to find a computable formula
for the Poincar\'e duality pairing. To be definite, we consider only the case
$i=0$ in \eqref{aa}. To this end, we start from an explicit geometric
realization of the groups. Let us first give a realization of the group
$K_0\left( \mathcal{A}_{M,\pi }\right)$.

\begin{lemma}
\label{lemc2} The group $K_0\left( \mathcal{A}_{M,\pi }\right) $ is isomorphic
to the group of stable homotopy classes of triples
\[
\left( E,F,\sigma \right) ,\qquad E,F\in \limfunc{Vect}\left( M\right)
,\;\sigma :\pi _{!}\left. E\right| _{\partial M}\longrightarrow \pi _{!}\left.
F\right| _{\partial M},
\]
where $\sigma $ is a vector bundle isomorphism. Here trivial triples are those
induced by a global vector bundle isomorphism over $M$.
\end{lemma}

\begin{proof}
This lemma is similar to Proposition~\ref{sem1}: both give a topological
realization of the $K_0$-group of the $C^*$-algebra of the
$\mathbb{Z}_n$-manifolds. Hence, a triple $\left(E,\Bbb{C}^k,\sigma\right)$
(note that an arbitrary triple can be reduced to this form) defines the element
\[
\left[ P_E\oplus P_2\right] -\left[ P_{\Bbb{C}^k}\oplus P_{\pi
_{!}{\Bbb{C}^k}}\right] \in K_0\left( \mathcal{A}_{M,\pi }\right)
\]
of the $K_0$--group, where the projection $P_2$ over $X\times \left[
0,1\right]$ is defined as
\[
P_2=P_{\pi _{!}E}\cos ^2\varphi +P_{\pi _{!}{\Bbb{C}^k}}\sin ^2\varphi +
(\sigma^{-1} \left( x\right)P_{\pi_{!}{\Bbb{C}^k}}+ \sigma \left( x\right)
P_{\pi _{!}E})\frac{\sin 2\varphi}2 ,\; \varphi =\frac \pi 2\left( 1-t\right),
\]
and $P_E$ and $P_{\Bbb{C}^k}\subset \Bbb{C}^N$ are some projections on
subbundles isomorphic to $E$ and ${\Bbb{C}^k}$, respectively. One must also
assume that these subbundles are mutually orthogonal.

The proof of the fact that this construction gives a well-defined isomorphism
with $K_0\left( \mathcal{A}_{M,\pi }\right) $ can be carried out by analogy
with the proof of Proposition \ref{sem1}.
\end{proof}

This realization enables one to define the product
\[
K_c^0\bigl( \overline{T^{*}M}^\pi \bigr) \times K_0\left( \mathcal{A}_{M,\pi
}\right) \longrightarrow K_0\left( \mathcal{A}_{T^{*}M,\pi }\right)
\]
geometrically in terms of operators with coefficients in vector bundles. More
precisely, for two elements
\[
\left[ \sigma \right] \in K^0_c\bigl( \overline{T^{*}M}^\pi \bigr) ,\qquad {\rm
and} \quad \left[ E,F,\sigma ^{\prime }\right] \in K_0\left( \mathcal{A}_{M,\pi
}\right)
\]
we consider the symbol
\begin{equation}
\sigma \otimes 1_E\oplus \sigma ^{-1}\otimes 1_F \label{st2}
\end{equation}
on $M$. The direct image of the restriction of this symbol to the boundary is
\begin{multline*}
\pi _{!}\left( \left. \sigma \otimes 1_E\oplus \sigma ^{-1}\otimes 1_F\right|
_{\partial M}\right)=  \pi _{!}\sigma \otimes 1_{\pi _{!}\left. E\right|
_{\partial M}}\oplus \pi _{!}\sigma ^{-1}\otimes 1_{\pi _{!}\left. F\right|
_{\partial M}} \\
\simeq \left( \pi _{!}\sigma \oplus \pi _{!}\sigma ^{-1}\right) \otimes 1_{\pi
_{!}\left. E\right| _{\partial M}}.
\end{multline*}
The last isomorphism is induced by the vector bundle isomorphism
\[
\pi _{!}\left. E\right| _{\partial M}\stackrel{\sigma ^{\prime }}{\simeq }
\pi _{!}\left. F\right| _{\partial M}.
\]
There is a standard homotopy of the symbol $\left( \pi _{!}\sigma \oplus \pi
_{!}\sigma ^{-1}\right) \otimes 1_{\pi _{!}\left. E\right| _{\partial M}}$ to
the identity:
\[
\begin{pmatrix}
\pi _{!}\sigma & 0 \\
0 & 1
\end{pmatrix}
\begin{pmatrix}
\cos \tau & \sin \tau \\
-\sin \tau & \cos \tau
\end{pmatrix}
\begin{pmatrix}
1 & 0 \\
0 & \pi _{!}\sigma ^{-1}
\end{pmatrix}
\begin{pmatrix}
\cos \tau & -\sin \tau \\
\sin \tau & \cos \tau
\end{pmatrix} ,\quad \tau \in \left[ 0,\pi /2\right] .
\]
Therefore, we have constructed an extension of the symbol (\ref{st2}) to an
elliptic symbol of some admissible operator on $M$. Finally, the desired
product
\[
\left[ \sigma \right] \times \left[ E,F,\sigma ^{\prime }\right] \in K_0\left(
\mathcal{A}_{T^{*}M,\pi }\right)
\]
can be defined as the difference element of this symbol. One can prove that the
pairing
\begin{equation}\label{sta}
\left\langle ,\right\rangle :K_c^0\bigl( \overline{T^{*}M}^\pi \bigr) \times
K_0\left( \mathcal{A}_{M,\pi }\right) \longrightarrow K_0\left( \mathcal{A}
_{T^{*}M,\pi }\right) \stackrel{\limfunc{ind}}{\longrightarrow }\Bbb{Z}
\end{equation}
defined as a composition of this product with the index mapping coincides with
Poincar\'e duality defined earlier.

\subsection{A topological index for $\mathbb{Z}_n$-manifolds}
Clearly, definition~\eqref{sta} of duality in the previous subsection  still
contains one component defined analytically. This is the index mapping
$$
\ind: K_0(\mathcal{A}_{T^*M,\pi})\longrightarrow \mathbb{Z}.
$$
A topological formula for it was obtained in \cite{SaSt10}. Let us briefly
recall this result.

We suppose that the following geometric condition is satisfied: there is a free
action of a finite group $G$ on $\partial M$ by diffeomorphisms such that
$\pi$ is the projection to the quotient space $\partial M/G$.

Consider two pairs $\left(M,\pi\right)$ and
$\left(M^{\prime},\pi^{\prime}\right) $.

\begin{definition}{\em
\emph{An }embedding \emph{\ }$f$\emph{\ of the pair }$\left( M,\pi \right)
$\emph{\ in }$\left( M^{\prime },\pi ^{\prime }\right) $\emph{\ is an embedding
of manifolds with boundary such that }$f:M\rightarrow M^{\prime },$\emph{\
}$f\left(\partial M\right) \subset \partial M^{\prime }$\emph{\ and the
restriction of $f$ to the boundary is a $G$-equivariant mapping. }}
\end{definition}
An embedding induces the direct image mapping
$$
f_!:K_*(\mathcal{A}_{T^*M,\pi})\longrightarrow K_*(\mathcal{A}_{T^*M',\pi'})
$$
of the $K$-groups.

For these $\mathbb{Z}_n$-manifolds, one can give a universal space in which an
arbitrary $\mathbb{Z}_n$-manifold can be embedded, the embedding being unique
up to homotopy. To this end, by $\pi _N:EG_N\longrightarrow BG_N$ we denote
the $N$-universal bundle for the group $G$ such that the spaces $EG_N$ and
$BG_N$ are smooth compact manifolds without boundary. (For the existence of
such models, e.g., see \cite{LuMi1}.)

\begin{proposition}
For a sufficiently large $N$, there exists an embedding $f$ of $\left( M,\pi
\right) $ in $\left( EG_N\times [0,\infty),\pi _N\right) .$ The embedding is
unique up to homotopy.
\end{proposition}

\begin{proof}[Sketch of Proof] By the $N$-universality of the covering
$\pi _N$, there exists an equivariant mapping $\partial M\rightarrow EG_N.$ If
the dimension of the space $EG_N$ is sufficiently large, then a general
position argument shows that slightly deforming this mapping one obtains a
smooth embedding. Then, by virtue of the $N$-connectedness of $EG_N$, this
extends to a mapping $M\rightarrow EG_N.$ In turn, by a small deformation
outside the boundary $\partial M$, this mapping can be made an embedding
globally.
\end{proof}

The index theorem can be stated with the use of an embedding in the universal
space. To this end, we introduce the $K$-group of the (infinite-dimensional)
classifying space for $\mathbb{Z}_n$-manifolds as the direct limit
\begin{equation}\label{clas}
K_*(\mathcal{A}_{ T^*EG\times [0,\infty),\pi_\infty})=
\lim\limits_{\longrightarrow} K_*(\mathcal{A}_{ T^*EG_N\times
[0,\infty),\pi_N})
\end{equation}
of the $K$-groups corresponding to the filtration of $EG$ by skeletons.

\begin{theorem} {\em \cite{SaSt10}}
The even $K$-group of the classifying space is isomorphic to $\mathbb{Z}$, and
the index can be computed in terms of the direct image mapping:
$$
\ind D=f_![\sigma(D)],
$$
$$
f_![\sigma(D)] \in K_0(\mathcal{A}_{ T^*EG\times [0,\infty),\pi_\infty})\simeq
\mathbb{Z}.
$$
\end{theorem}

\providecommand{\bysame}{\leavevmode\hbox to3em{\hrulefill}\thinspace}


\begin{thebibliography}{APS76b}

\bibitem[AB64]{AtBo2}
M.F. Atiyah and R.~Bott, \emph{The index problem for manifolds with boundary},
Bombay Colloquium on Differential Analysis (Oxford), Oxford University Press,
1964, pp.~175--186.

\bibitem[ABP73]{ABP1}
M.F. Atiyah, R.~Bott, and V.K. Patodi, \emph{On the heat equation and the index
theorem}, Invent. Math. \textbf{19} (1973), 279--330.

\bibitem[AD62]{AgDy1}
M.~S. Agranovich and A.~S. Dynin, \emph{General boundary-value problems for
elliptic systems in higher-dimensional regions}, Dokl. Akad. Nauk SSSR
\textbf{146} (1962), 511--514.

\bibitem[ADS83]{ADS1}
M.F. Atiyah, H.~Donnelly, and I.M. Singer, \emph{Eta invariants, signature
defects of cusps and values of {$L$-functions}}, Annals Math. \textbf{118}
(1983), 131--177.

\bibitem[Agr94]{Agr3}
M.~S. Agranovich, \emph{{Elliptic operators on closed manifolds.}}, Partial
differential equations. VI., {Encycl. Math. Sci.}, vol.~63, 1994, pp.~1--130
(English).

\bibitem[APS75]{APS1}
M.~Atiyah, V.~Patodi, and I.~Singer, \emph{Spectral asymmetry and {Riemannian}
geometry {I}}, Math. Proc. Cambridge Philos. Soc. \textbf{77} (1975), 43--69.

\bibitem[APS76a]{APS2}
\bysame, \emph{Spectral asymmetry and {Riemannian} geometry {II}}, Math. Proc.
Cambridge Philos. Soc. \textbf{78} (1976), 405--432.

\bibitem[APS76b]{APS3}
\bysame, \emph{Spectral asymmetry and {Riemannian} geometry {III}}, Math. Proc.
Cambridge Philos. Soc. \textbf{79} (1976), 71--99.

\bibitem[AS63]{AtSi0}
M.F. Atiyah and I.M. Singer, \emph{The index of elliptic operators on compact
manifolds}, Bull. Amer. Math. Soc. \textbf{69} (1963), 422--433.

\bibitem[AS68]{AtSi1}
\bysame, \emph{The index of elliptic operators {I}}, Ann. of Math. \textbf{87}
(1968), 484--530.

\bibitem[Ati61]{Ati6}
M.~F. Atiyah, \emph{{{C}haracters and cohomology of finite groups}}, Publ.
Math. IHES \textbf{9} (1961), 23--64.

\bibitem[Ati69]{Ati4}
\bysame, \emph{Global theory of elliptic operators}, Proc. of the Int.
Symposium on Functional Analysis (Tokyo), University of Tokyo Press, 1969,
pp.~21--30.

\bibitem[Ati89]{Ati2}
\bysame, \emph{K-theory}, second ed., The Advanced Book Program,
Addison--Wesley, Inc., 1989.

\bibitem[BBW93]{BBW1}
B.~Boo\ss-Bavnbek and K.~Wojciechowski, \emph{Elliptic boundary problems for
{D}irac operators}, Birkh{\"a}user, Boston--Basel--Berlin, 1993.

\bibitem[BD82]{BaDo2}
P.~Baum and R.~G. Douglas, \emph{{{T}oeplitz operators and {P}oincare
duality}}, Toeplitz centennial ({Toeplitz Mem. Conf. Tel Aviv 1981}),
Operator Theory, Adv. Appl., vol.~4, 1982, pp.~137--166.

\bibitem[BDF77]{BDF1}
L.~Brown, R.~Douglas, and P.~Fillmore, \emph{{E}xtensions of
$\mathbb{C}^*$-algebras and ${K}$-homology}, Ann. Math. II \textbf{105}
(1977), 265--324.

\bibitem[BDT89]{BDT1}
P.~Baum, R.~G. Douglas, and M.~E. Taylor, \emph{{{C}ycles and relative cycles
in analytic ${K}$-homology}}, J. Differ. Geom. \textbf{30} (1989), no.~3,
761--804 (English).

\bibitem[Bla98]{Bla1}
B.~Blackadar, \emph{{$K$-theory for operator algebras}}, Mathematical Sciences
Research Institute Publications, no.~5, Cambridge University Press, 1998,
Second edition.

\bibitem[BN94]{BrNi1}
J.-L. Brylinski and V.~Nistor, \emph{{C}yclic cohomology of etale groupoids},
K-theory \textbf{8} (1994), 341--365.

\bibitem[Bot92]{Bot2}
B.~Botvinnik, \emph{Manifolds with singularities and the {A}dams-{N}ovikov
spectral sequence}, London Mathematical Society Lecture Note Series, vol.
170, Cambridge University Press, Cambridge, 1992.

\bibitem[BT82]{BoTu1}
R.~Bott and L.~Tu, \emph{{D}ifferential forms in algebraic topology}, Graduate
Texts in Mathematics, vol.~82, Springer--Verlag, Berlin--Heidelberg--New
York, 1982.

\bibitem[Con94]{Con1}
A.~Connes, \emph{Noncommutative geometry}, Academic Press Inc., San Diego, CA,
1994.

\bibitem[Don78]{Don3}
H.~Donnelly, \emph{{{E}ta-invariants for ${G}$-spaces}}, Indiana Univ. Math. J.
\textbf{27} (1978), 889--918.

\bibitem[DZ98]{DaZh3}
X.~Dai and W.~Zhang, \emph{Higher spectral flow}, J. Funct. Anal. \textbf{157}
(1998), no.~2, 432--469.

\bibitem[FM92]{FrMe1}
D.~Freed and R.~Melrose, \emph{A mod $k$ index theorem}, Invent. Math.
\textbf{107} (1992), no.~2, 283--299.

\bibitem[Gil81]{Gil4}
P.B. Gilkey, \emph{The residue of the global eta function at the origin}, Adv.
in Math. \textbf{40} (1981), 290--307.

\bibitem[Gil89a]{Gil2}
\bysame, \emph{The geometry of spherical space form groups}, Series in Pure
Mathematics (Singapore), vol.~7, World Sientific Publ. Co. Pte. Ltd., 1989.

\bibitem[Gil89b]{Gil7}
\bysame, \emph{{T}he eta invariant of even order operators}, Lecture Notes in
Mathematics \textbf{1410} (1989), 202--211.

\bibitem[Gil95]{Gil1}
\bysame, \emph{{I}nvariance theory, the heat equation, and the
{A}tiyah-{S}inger index theorem}, second ed., Studies in Advanced Mathematics,
CRC Press, Boca Raton, FL, 1995.

\bibitem[GS95]{GrSe2}
G.~Grubb and R.~Seeley, \emph{Weakly parametric pseudodifferential operators
and {Atiyah--Patodi--Singer} boundary problems}, Invent. Math. \textbf{121}
(1995), 481--529.

\bibitem[Hir73]{Hirz2}
F.~Hirzebruch, \emph{Hilbert modular surfaces}, Enseignement Math. \textbf{19}
(1973), 183--281.

\bibitem[H{\"o}r85]{Hor3}
L.~H{\"o}rmander, \emph{The analysis of linear partial differential operators.
{III}}, Springer--Verlag, Berlin Heidelberg New York Tokyo, 1985.

\bibitem[HR00]{HiRo1}
N. Higson and J. Roe, \emph{Analytic {$K$}-homology}, Oxford University Press,
Oxford, 2000.

\bibitem[Kas73]{Kas1}
G.G. Kasparov, \emph{{{T}he generalized index of elliptic operators}}, Funct.
Anal. Appl. \textbf{7} (1973), 238--240 (English. Russian original).

\bibitem[Kas88]{Kas3}
\bysame, \emph{Equivariant ${KK}$-theory and the {N}ovikov conjecture}, Inv.
Math. \textbf{91} (1988), no.~1, 147--201.

\bibitem[LM89]{LaMi1}
H.B. Lawson and M.L. Michelsohn, \emph{{S}pin geometry}, Princeton Univ. Press,
Princeton, 1989.

\bibitem[LM98]{LuMi1}
G.~Luke and A.S. Mishchenko, \emph{{Vector bundles and their applications}},
Mathematics and its Applications, vol. 447, Kluwer Academic Publishers,
Dordrecht, 1998 (English).

\bibitem[Mel93]{Mel2}
R.~Melrose, \emph{The Atiyah--Patodi--Singer index theorem}, Research Notes in
Mathematics, A. K. Peters, Boston, 1993.

\bibitem[M{\"u}l84]{Mul5}
W.~M{\"u}ller, \emph{Signature defects of cusps of {Hilbert} modular varieties
and values of {$L$-series} at $s=1$}, J. Diff. Geometry \textbf{20} (1984),
55--119.

\bibitem[NSS99]{NScS5}
V.~Nazaikinskii, B.-W. Schulze, and B.~Sternin, \emph{On the homotopy
classification of elliptic operators on manifolds with singularities}, Univ.
Potsdam, Institut f{\"u}r Mathematik, Potsdam, Oktober 1999, Preprint N
99/21.

\bibitem[NSSS98]{NScSS3}
V.~Nazaikinskii, B.-W. Schulze, B.~Sternin, and V.~Shatalov, \emph{Spectral
boundary value problems and elliptic equations on singular manifolds},
Differential Equations, {\bf 34}, N 5 (1998), pp 696--710.

\bibitem[Pal65]{Pal1}
R.S. Palais, \emph{Seminar on the {A}tiyah--{Singer} index theorem}, Princeton
Univ. Press, Princeton, NJ, 1965.

\bibitem[Phi96]{Phil1}
J.~Phillips, \emph{{{S}elf-adjoint Fredholm operators and spectral flow}},
Canad. Math. Bull. \textbf{39} (1996), no.~4, 460--467.

\bibitem[Sal95]{Sal1}
D.~Salamon, \emph{{{T}he spectral flow and the {M}aslov index}}, Bull. LMS
\textbf{27} (1995), no.~1, 1--33.

\bibitem[Sav99]{Sav1}
A.Yu. Savin, \emph{{O}n operators that allow the decomposition of index
formulas for spectral boundary value problems}, Doklady Mathematics
\textbf{60} (1999), no.~2, 220--222.

\bibitem[See67]{See5}
R.T. Seeley, \emph{Complex powers of an elliptic operator}, Proc. Sympos. Pure
Math. \textbf{10} (1967), 288--307.

\bibitem[SS99]{SaSt1}
A.Yu. Savin and B.Yu. Sternin, \emph{Elliptic operators in even subspaces},
Matem. sbornik \textbf{190} (1999), no.~8, 125--160, English transl.:
Sbornik: Mathematics {\bf 190}, N 8 (1999), p. 1195--1228; arXiv:
math/9907027.

\bibitem[SS00a]{SaSt7}
\bysame, \emph{{E}ta {I}nvariant and {P}arity {C}onditions},
Univ. Potsdam, Institut f{\"u}r Mathematik, Potsdam, Oktober 2000, {Preprint
00/21, available at
http://www.math.uni-potsdam.de/a\_partdiff/prepr/2000\_21.zip}.

\bibitem[SS00b]{SaSt2}
\bysame, \emph{{E}lliptic operators in odd subspaces}, Sbornik: Mathematics
{\bf 191}, N 8 (2000), arXiv: math/9907039.

\bibitem[SS01]{SaSt10}
\bysame, \emph{Index {D}efects in the {T}heory of {N}onlocal
{B}oundary {V}alue {P}roblems and the $\eta$-invariant}, Univ. Potsdam,
Institut f{\"u}r Mathematik, Potsdam, November 2001, Preprint 01/31, arXiv:
math/0108107.

\bibitem[SS02]{SaSt6}
\bysame, \emph{{T}he {E}ta-invariant and {P}ontryagin {D}uality
in ${K}$-{T}heory}, Math. Notes \textbf{71} (2002), no.~2, 245--261, arXiv:
math/0006046.

\bibitem[SSS98]{ScSS18}
B.-W. Schulze, B.~Sternin, and V.~Shatalov, \emph{On general boundary value
problems for elliptic equations},  Sbornik: Mathematics {\bf 189}, N 10 (1998),
p. 1573--1586.

\bibitem[SSS99a]{SaScS1}
A.~Savin, B.-W. Schulze, and B.~Sternin, \emph{On invariant index formulas for
spectral boundary value problems}, Differentsial'nye uravnenija \textbf{35}
(1999), no.~5, 709--718.

\bibitem[SSS99b]{SaScS4}
\bysame, \emph{The homotopy classification and
the index of boundary value problems for general elliptic operators}, Univ.
Potsdam, Institut f{\"u}r Mathematik, Oktober 1999, Preprint N 99/20, arXiv:
math/9911055.

\bibitem[SSS01]{SaScS8}
\bysame, \emph{On the homotopy classification
of elliptic boundary value problems}, Partial differential equations and
spectral theory: PDE2000 Conference in Clausthal, Germany (Basel, Boston,
Berlin), Operator theory: advances and applications, vol. 126,
Birkh{\"a}user, 2001, pp.~299--306.

\bibitem[Sul70]{Sul1}
D.~Sullivan, \emph{{G}eometric topology. {L}ocalization, periodicity and
{G}alois symmetry}, MIT, Cambridge, Massachusets, 1970.

\end{thebibliography}
\end{document}